\documentclass[smallextended]{svjour3}      % onecolumn (second format)
\smartqed  % flush right qed marks, e.g. at end of proof
\usepackage{lineno}
\modulolinenumbers[5]

\usepackage{mathptmx}      % use Times fonts if available on your TeX system

\usepackage[colorlinks=true]{hyperref}
\usepackage{xparse}
\usepackage{amsmath}
\usepackage{amsfonts}

\newcommand{\bxi}{\boldsymbol{\xi}}

\newcommand{\dVblank}{\, {\rm d}V}
\newcommand{\dvblank}{\, {\rm d}v}

\newcommand{\eps}{\varepsilon}

\NewDocumentCommand \dV{ o }{%
    \IfNoValueTF{#1}{\dVblank}%
    {
        \dV_{#1}
    }
}

\NewDocumentCommand \dv{ o }{%
    \IfNoValueTF{#1}{\dvblank}%
    {
        \dv_{#1}
    }
}

\NewDocumentCommand \mc{ m }{%
    \mathcal{#1}%
}

\NewDocumentCommand \pN{ m m m }{%
    \frac{\partial^{#3} #1}{\partial #2^{#3}}%
}

\NewDocumentCommand \an{ m }{%
    \langle {#1} \rangle%
}

\NewDocumentCommand \mbf{ m }{%
    \mathbf{#1}
}

\NewDocumentCommand \mbt{ m }{%
    \mbf{\tilde{#1}}
}

\NewDocumentCommand \bs{ m }{%
    \boldsymbol{#1}
}

\NewDocumentCommand \stateOne{ m }{
    \underline{\mbf{#1}}%
}

\NewDocumentCommand \stateTwo{ m o }{%
    \IfNoValueTF{#2}{\stateOne{#1}}
    {
    \stateOne{#1}\an{#2}%
    }
}

\NewDocumentCommand \s{ m o o }{%
    \IfNoValueTF{#3}{\stateTwo{#1}[#2]}%
    {
    \underline{\mbf{#1}}(#3)\an{#2}%
    }
}

\NewDocumentCommand \stateOneI{ m m }{
    \underline{#1}_{#2}%
}

\NewDocumentCommand \stateTwoI{ m m o }{%
    \IfNoValueTF{#3}{\stateOneI{#1}{#2}}
    {
    \stateOneI{#1}{#2}\an{#3}%
    }
}

\NewDocumentCommand \sI{ m m o o }{%
    \IfNoValueTF{#4}{\stateTwoI{#1}{#2}[#3]}%
    {
        \stateOneI{#1}{#2}(\mbf{#4})\an{#3}
    }
}

\NewDocumentCommand \snsI{ m m o }{%
    \IfNoValueTF{#3}{\underline{\mbf{#1}}_{#2}}%
    {
      \underline{\mbf{#1}}_{#2}(#3)%
    }
}

\NewDocumentCommand \stateOneN{ m }{
    \underline{#1}%
}

\NewDocumentCommand \stateTwoN{ m o }{%
    \IfNoValueTF{#2}{\stateOneN{#1}}
    {
    \stateOneN{#1}\an{#2}%
    }
}

\NewDocumentCommand \sN{ m o o }{%
    \IfNoValueTF{#3}{\stateTwoN{#1}[#2]}%
    {
    \underline{#1}(#3)\an{#2}%
    }
}

\NewDocumentCommand\presuper{ m m }{%
  {\mathop{}%
   \mathopen{\vphantom{#2}}^{#1}%
   \kern-1\scriptspace%
   #2}
}

\NewDocumentCommand\presuperz{ m m m }{%
  {\mathop{}%
   \mathopen{\vphantom{#2}}^{#1}%
   \kern-3\scriptspace%
   #2\kern-8\scriptspace_{#3}%
}
}

\usepackage{amssymb}
\usepackage[numbers]{natbib}

\usepackage{pgfplots}
\usetikzlibrary{arrows}
\usetikzlibrary{calc}

\usepackage[english]{babel}
\usepackage[utf8x]{inputenc}
\usepackage[T1]{fontenc}

\usepackage{amsmath}
\usepackage{graphicx}
\usepackage{subfig}
\usepackage{bm} % bold math symbols

\makeatletter
\let\cl@chapter\undefined
\makeatletter

\usepackage[nameinlink]{cleveref}
\crefname{equation}{Eq.}{Eqs.}
\crefname{figure}{Fig.}{Figs.}
\crefname{section}{Section}{Sections}
\crefname{remark}{Remark}{Remarks}
\crefname{example}{Example}{Examples}
\crefname{appendix}{}{}

\hypersetup{allcolors=blue}

\makeatletter
\newcommand{\oset}[3][0ex]{%
  \mathrel{\mathop{#3}\limits^{
    \vbox to#1{\kern-2\ex@
    \hbox{$\scriptstyle#2$}\vss}}}}
\makeatother

\begin{document}

% Remove spacing around align
\setlength{\abovedisplayskip}{6pt}
\setlength{\belowdisplayskip}{6pt}

\title{A unified, stable and accurate meshfree framework for peridynamic correspondence modeling. Part II: wave propagation and enforcement of stress boundary conditions}

\author{Masoud Behzadinasab \and
        John T. Foster \and
        Yuri Bazilevs
}

%\authorrunning{Short form of author list} % if too long for running head

\institute{M. Behzadinasab \at
              School of Engineering, Brown University, 184 Hope St., Providence, RI 02912, USA \\
              \email{masoud\_behzadinasab@brown.edu}           
           \and
           John T. Foster \at
              Hildebrand Department of Petroleum \& Geosystems Engineering, The University of Texas at Austin, Austin, TX 78712, USA
           \and
           Yuri Bazilevs \at
              School of Engineering, Brown University, 184 Hope St., Providence, RI 02912, USA
}

\date{Received: date / Accepted: date}
% The correct dates will be entered by the editor

\maketitle

\begin{abstract}
The overarching goal of this work is to develop an accurate, robust, and stable methodology for finite deformation modeling using strong-form peridynamics (PD) and the correspondence modeling framework. We adopt recently developed methods that make use of higher-order corrections to improve the computation of integrals in the correspondence formulation. A unified approach is presented that incorporates the reproducing kernel (RK) and generalized moving least square (GMLS) approximations in PD to obtain higher-order non-local gradients. We show, however, that the improved quadrature rule does not suffice to handle instability issues that have proven problematic for the correspondence-model based PD. In Part I of this paper, a bond-associative, higher-order core formulation is developed that naturally provides stability without introducing artificial stabilization parameters. Numerical examples are provided to study the convergence of RK-PD, GMLS-PD, and their bond-associated versions to a local counterpart, as the degree of non-locality (i.e., the horizon) approaches zero. Problems from linear elastostatics are utilized to verify the accuracy and stability of our approach. It is shown that the bond-associative approach improves the robustness of RK-PD and GMLS-PD formulations, which is essential for practical applications. The higher-order, bond-associated model can obtain second-order convergence for smooth problems and first-order convergence for problems involving field discontinuities, such as curvilinear free surfaces. In Part II of this paper we use our unified PD framework to: (a) study wave propagation phenomena, which have proven problematic for the state-based correspondence PD framework; (b) propose a new methodology to enforce natural boundary conditions in correspondence PD formulations, which should be particularly appealing to coupled problems. Our results indicate that bond-associative formulations accompanied by higher-order gradient correction provide the key ingredients to obtain the necessary accuracy, stability, and robustness characteristics needed for engineering-scale simulations.

\end{abstract}

%\linenumbers

\section{Introduction}
\label{sec:intro}

In Part I of this series, we provided a unified meshfree approach to higher-order peridynamic (PD) correspondence modeling, namely the reproducing kernel peridynamics (RK-PD) \citep{hillman2020generalized} (in which the non-local gradient operator is equivalent to the PD differential operator \citep{madenci2016peridynamic}) and GMLS-PD based on the generalized moving least squares operator \citep{trask2019asymptotically} were considered in our unified approach. Both of these algorithms can compute non-local gradients with $n^\text{th}$-order accuracy. We showed that the higher-order corrections cannot eliminate the well-known issue of numerical instability in correspondence PD theory, which manifest itself as zero-energy mode oscillations in simulations that involve inhomogeneous deformations \citep{littlewood2010simulation,breitenfeld2014non,tupek2014extended,silling2017stability,behzadinasab2020stability}, and may lead to large prediction errors in practical applications \citep{kramer2019sandia}. To provide stability for the higher-order framework, we developed a continuum-based, bond-associated model, inspired by \citep{breitzman2018bond,chen2018bond,behzadinasab2020semi}, which takes into account non-uniform deformations and maintains stability without additional stabilization mechanisms employed. It was shown that the higher-order, bond-associated model can obtain second-order accuracy for smooth problems and first-order accuracy for non-smooth ones (e.g., involving solution-derivative discontinuity).

In Part II of this work, we aim to study how our unified, stabilized PD formulation handles wave propagation phenomena and to develop a convenient methodology for incorporating natural boundary conditions in PD. 

Wave propagation phenomena have proven problematic for the correspondence formulation of PD \citep{silling2017stability,bazant2016wave}, as the formulation becomes significantly dispersive for relatively high wave numbers, imposing practical difficulties in a group of common PD dynamic problems such as impact and penetration problems (or general fracture mechanics problems involving rapid crack growth). Here we show how the bond-associative, higher-order approaches improve the behavior of the correspondence-based PD setting in capturing wave propagation phenomena. 

Traditionally, stress boundary conditions have been represented as volumetric body forces in PD, which may work for simple loading scenarios (see, e.g., \citep{ha2010studies,bobaru2015cracks} involving constant stress on flat boundaries). It is, however, not clear how to represent general stress conditions by body forces, which considerably limits the general applicability of PD-based approaches. While Galerkin PD formulations can be developed (see, e.g. \citep{chen2011continuous,madenci2018weak}) and used to include natural boundary conditions in the weak form, the weak-form strategy involves a six-dimensional (double) integration \citep{littlewood2015roadmap,bobaru2016handbook}, imposing significant computational cost. Thus, strong-form-based approaches, commonly employed in meshfree methods \citep{silling2005meshfree}, are appealing for practical applications. We introduce a convenient technique, while maintaining robustness, to incorporate natural boundary conditions in the strong-form PD, based on the stabilized, higher-order correspondence framework. In the presented approach, we assume that limited information is available on boundaries. That is, either essential or natural boundary conditions are known at each boundary point. The methodology developed here should be particularly appealing to coupled problems, e.g., fluid-structure interaction modeling where pressure and viscous stresses from the surrounding fluid are imposed as natural boundary conditions on the structure surface.

The remainder of the paper is organized as follows. The governing equations of the bond-associated, higher-order PD formulation are summarized in \cref{sec:formulation}. Wave propagation is studied in \cref{sec:wave}. A novel approach for non-local modeling of natural boundary conditions in the strong-form PD is presented in \cref{sec:natural}, including the numerical examples to show the performance of the proposed methodology. Concluding remarks and future directions are summarized in \cref{sec:conclusions}. 

In what follows, all vectors are column vectors unless otherwise noted. Bold symbols indicate tensors of rank 1 (i.e., vectors) or higher.

\section{Stabilized, Higher-Order Correspondence Model of PD}
\label{sec:formulation}

In this section, the approach detailed in Part I of this paper to unify correspondence material modeling, incorporating bond-associative and higher-order corrections, is summarized.

The goal of these methods is to use non-local integration techniques, rather than differentiation, to approximate classical continuum theory. That is, for a problem in the Lagrangian frame, {\em dynamic equilibrium} exists if the internal forces are in balance with the inertial and externally applied forces, i.e., 
\begin{equation}
   \rho \pN{\mbf{u}}{t}{2}\Big|_{(\mbf{X},t)} - \nabla \cdot \mbf{P}(\mbf{X}) - \mbf{b}(\mbf{X}) = 0 , 
  \label{eqn:classical}
\end{equation}
where $\mbf{X}$ is a material point in the reference configuration, $\mbf{u}$ is the displacement field, $\rho$ is the mass density in the undeformed configuration, and $\mbf{b}$ is the body force per uint volume. The static case is obtained by setting the first term in Eq.~(\ref{eqn:classical}) to zero.

The second term in Eq.~(\ref{eqn:classical}) is the divergence of the first Piola--Kirchhoff stress $\mbf{P}$, which is the energy conjugate of the deformation gradient $\mbf{F}$ defined as
$$\mbf{F} = \mbf{I} + \nabla \mbf{u}, $$
where $\mbf{I}$ is the Identity Matrix. While local models evaluate the gradient operator using partial differentiation, the noted PD approach uses the following integration technique to compute a gradients at the discrete level:
\begin{align}
  \left(\nabla_h \mbf{f}\right)_I = \sum_{J\in\mc{H}_I} \left[\mbf{f}_{J} - \mbf{f}_I\right] \, \bs\gamma_{IJ}^\intercal , 
  \label{eqn:grad}
\end{align}
where $I$ and $J$ denote discrete nodes, and $\mc{H}_I$ is the {\em family} of $I$ defined as 
$$ \mc{H}_I = \left\{ J \ | \ J \in \tilde{\mc{B}}_0, \, 0 < | \mbf{X}_J - \mbf{X}_I | \leq \delta \right\} , $$
where $\delta$ is called the {\em horizon} for the body $\mc{B}_0$ in the reference configuration. The coefficients $\bs\gamma_{IJ}$ are obtained by RK and GMLS methods such that higher-order accurate gradients may be computed provided the polynomial unisolvency is maintained (i.e., enough neighbors exist for representing the space of $n^\text{th}$-order polynomials; see Part I for more details). 

\begin{itemize}
  \item RK-PD: 
  $$\bs\gamma_{IJ} = \bs\Phi_{IJ} \, V_J ,$$
  where $V_J$ is the volume of node $J$, and 
  $$\bs\Phi_{IJ} = \alpha_{IJ} \, \oset[.2ex]{\nabla}{\mbf{Q}}^\intercal \, \mbt{M}^{-1}_I \, \mbf{Q}_{IJ} .$$
  $\alpha$ is a weighting function, known as the {\em influence function} in PD, which is a function of the relative distance between material points with respect to $\delta$. The moment matrix $\mbt{M}$ is defined by 
  $$\mbt{M}_I = \sum_{J\in\mc{H}_I} \alpha_{IJ} \, \mbf{Q}_{IJ} \, \mbf{Q}_{IJ}^\intercal \, V_J .$$
  $\mbf{Q}(\bxi)$ is a vector of the set of monomials $\{\bxi^{\, \beta}\}_{|\beta|=1}^{n}$, with a free parameter $n$ that determines the degree of accuracy of the non-local gradient operator. $\oset[.2ex]{\nabla}{\mbf{Q}}$ is given by
  $$\oset[.2ex]{\nabla}{\mbf{Q}} \; \equiv \mbf{Q}^{(\delta_{j1}, \, \delta_{j2}, \, \delta_{j3})} = [0 \ , \ \dots \ , \ 0 \ , \ \underset{jth \; {\rm entry}}{1} \ , \ 0 \ , \ \dots \ , \ 0]^\intercal.$$

  \item GMLS-PD: 
  $$\bs\gamma_{IJ} = \bs\omega_{IJ} \, \dfrac{\left[\mbf{X}_J-\mbf{X}_I\right]}{\left|\mbf{X}_J-\mbf{X}_I\right|^2} , $$
  where $\bs\omega_{IJ}$ is a set of quadrature weights associated with the family of $I$. The weights are obtained by solving the following constrained optimization problem: 
  \begin{align*}
    & \underset{\bs\omega_{IJ}}{\rm argmin} \sum_{J\in\mc{H}_I} [ \bs\omega_{IJ}:\bs\omega_{IJ} ] \notag \\
    & {\rm such \ that,} ~~~~ \left(\nabla_h [p]\right)_I = \left(\nabla [p]\right)_I ~~~ \forall p \in \mbf{V}_h , 
  \end{align*}
  in which $\mbf{V}_h$ denotes a Banach space of functions whose gradients should be reproduced exactly. Selecting $\mbf{V}_h$ as the space of $n^\text{th}$-order polynomials would be beneficial to obtain higher-order gradients. 
\end{itemize}
For both RK-PD and GMLS-PD (with the case of polynomials as the reproducing condition), it can be shown that $\Big[ \dfrac{(n+d)!}{n!d!}-1 \Big]$ number of non-colinear (in 2D) or non-coplanar (in 3D) bonds are required to obtain polynomial unisolvency \citep{liu1997moving}, where $d$ is the space dimension. The horizon size, for uniform discretizations, must be chosen greater than $n \times h$ with $h$ the nodal spacing.

Using \cref{eqn:grad}, a higher-order, non-local, deformation gradient is given by
\begin{align}
  \mbt{F}_I = \mbf{I} + \sum_{J\in\mc{H}_I} \left[\mbf{u}_{J} - \mbf{u}_I\right] \, \bs\gamma_{IJ}^\intercal . 
  \label{eqn:F}
\end{align}
Similarly, by employing the gradient operator on the stress field and applying the trace operator, a higher-order, non-local divergence of the stress is obtained by
\begin{align}
  \left(\nabla_h \cdot \mbf{P}\right)_I
  &\approx \left(\nabla_h \cdot \mbt{P}\right)_I \notag \\
  &= \sum_{J\in\mc{H}_I} \left[\mbt{P}_{J} - \mbt{P}_I\right] \, \bs\gamma_{IJ} ,
  \label{eqn:div}
\end{align}
where $\mbt{P} = \mbf{P}(\mbt{F})$.

To naturally stabilize the model, a bond-associative correction is implemented by modifying \cref{eqn:div} as
\begin{align}
  \left(\nabla_h \cdot \mbt{P}\right)_I = \sum_{J\in\mc{H}_I} \left[\mbt{P}_{JI} - \mbt{P}_I\right] \, \bs\gamma_{IJ} ,
  \label{eqn:BA-div}
\end{align}
where $\mbt{P}_{JI}$ is called the {\em bond-associated} first Piola--Kirchhoff stress that is the energy conjugate of the bond-associated deformation gradient, i.e., $\mbt{P}_{JI} = \mbf{P}(\mbt{F}_{JI})$. Here, $\mbt{F}_{JI}$ is defined as
\begin{align}
  \mbt{F}_{JI} 
  &= \mbt{F}_J + \Delta \mbt{F}_{JI}^{\, nh} \notag \\
  &= \mbt{F}_J + \left[ \mbf{x}_J - \mbf{x}_I - \frac{\mbt{F}_I+\mbt{F}_J}{2}[\mbf{X}_J-\mbf{X}_I] \right] \frac{[\mbf{X}_J-\mbf{X}_I]^\intercal}{|\mbf{X}_J-\mbf{X}_I|^2} ,
  \label{eqn:BA-F}
\end{align}
where $\mbt{F}_J$ is the neighbor deformation gradient, evaluated by \cref{eqn:F}, and $\mbf{x}$ is the current position vector field, i.e., $\mbf{x} = \mbf{X} + \mbf{u}$. To account for non-uniform deformations, $\Delta \mbt{F}_{JI}^{\, nh}$ is included in the bond-level deformation gradient. $\Delta \mbt{F}_{JI}^{\, nh} = \mbf{0}$ for homogeneous deformations. Neglecting $\Delta \mbt{F}_{JI}^{\, nh}$ recovers the base divergence operator, i.e., \cref{eqn:div}.

\section{Wave Propagation}
\label{sec:wave}

In section, we use a 1D linear elastic wave propagation problem to study the dispersion relations of the proposed RK-PD and GMLS-PD models and their bond-associative counterparts. This study allows us to gain insight into the effectiveness of these formulations in dynamic problems. For an infinite 1D bar, the following plane wave solution may be written explicitly
\begin{equation}
  u(X,t) = u_0 {\rm e}^{{\rm i}(kX - \omega t)} , 
  \label{eqn:wave}
\end{equation}
where $u_0$ is the wave amplitude, $k$ is the wave number, and $\omega$ is the angular frequency. In addition, the following relation holds:
$$ \omega / k = c , $$
where $c$ is the elastic wave propagation speed given by
$$ c = \sqrt{E/\rho} . $$

To set up this problem with the methods proposed, a bar of size $[-2\delta, 2\delta]$ (long enough to remove the boundary effects for the center node) is discretized using a uniform and a non-uniform grid. The material properties are set to $E=1$ and $\rho=1$. The nodal displacements are prescribed as per \cref{eqn:wave}, and the angular frequency is computed by solving the equation of motion for the center node ($X=0$) at $t=0$:
%\begin{equation}
$$ \rho \pN{u}{t}{2}\Big|_{(I,0)} = \frac{\partial_h P}{\partial X}\Big|_{(I,0)} , $$
where $I$ refers to the node placed at $X=0$. Using \cref{eqn:wave,eqn:div}, and rearranging terms, we obtain
$$ - \rho \omega^2 u_0 = \sum_{J\in\mc{H}_I} \left[\tilde{P}_{J}(0) - \tilde{P}_I(0)\right] \, \gamma_{IJ} . $$
A 1D linear material model, $\tilde{P} = E (\tilde{F} - 1)$, is adopted here. Quadratic ($n=2$) discretizations are employed with the average nodal spacing $h=1$. For the non-uniform case, a 15\% random normal perturbation is applied to the nodal positions.  In addition, to study the horizon-size effect, $\delta=3h$ and $\delta=4h$ are considered. Note that if $\delta$ is chosen to scale with $h$, $kh$ (or $k\delta$) becomes an invariant in the mechanical response.

\begin{figure*}[!ht]
  \centering
  \subfloat[][$\delta = 3 \, h$]{\includegraphics[height=0.4\textwidth]{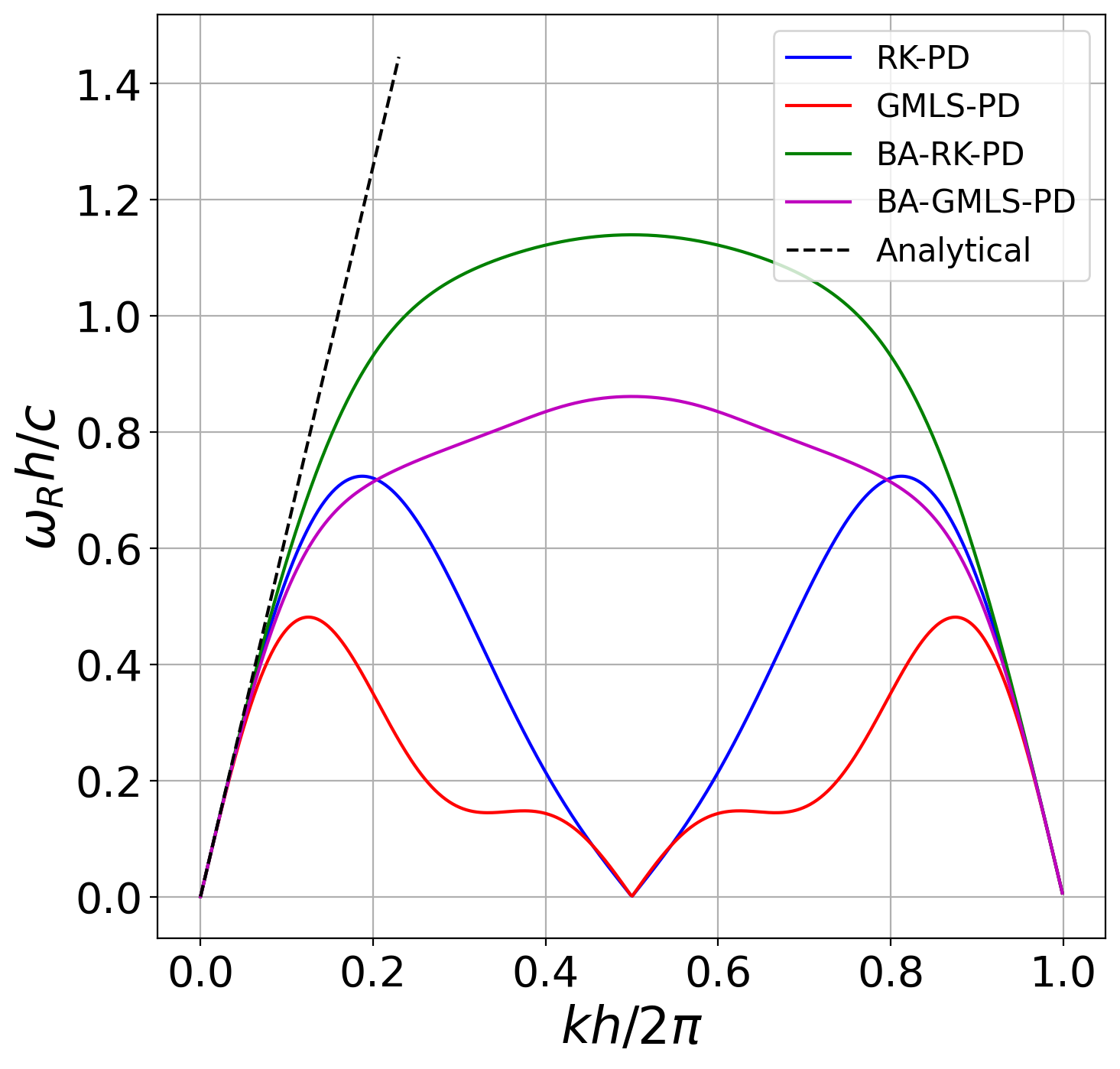}}
  \hspace*{0.5cm}
  \subfloat[][$\delta = 4 \, h$]{\includegraphics[height=0.4\textwidth]{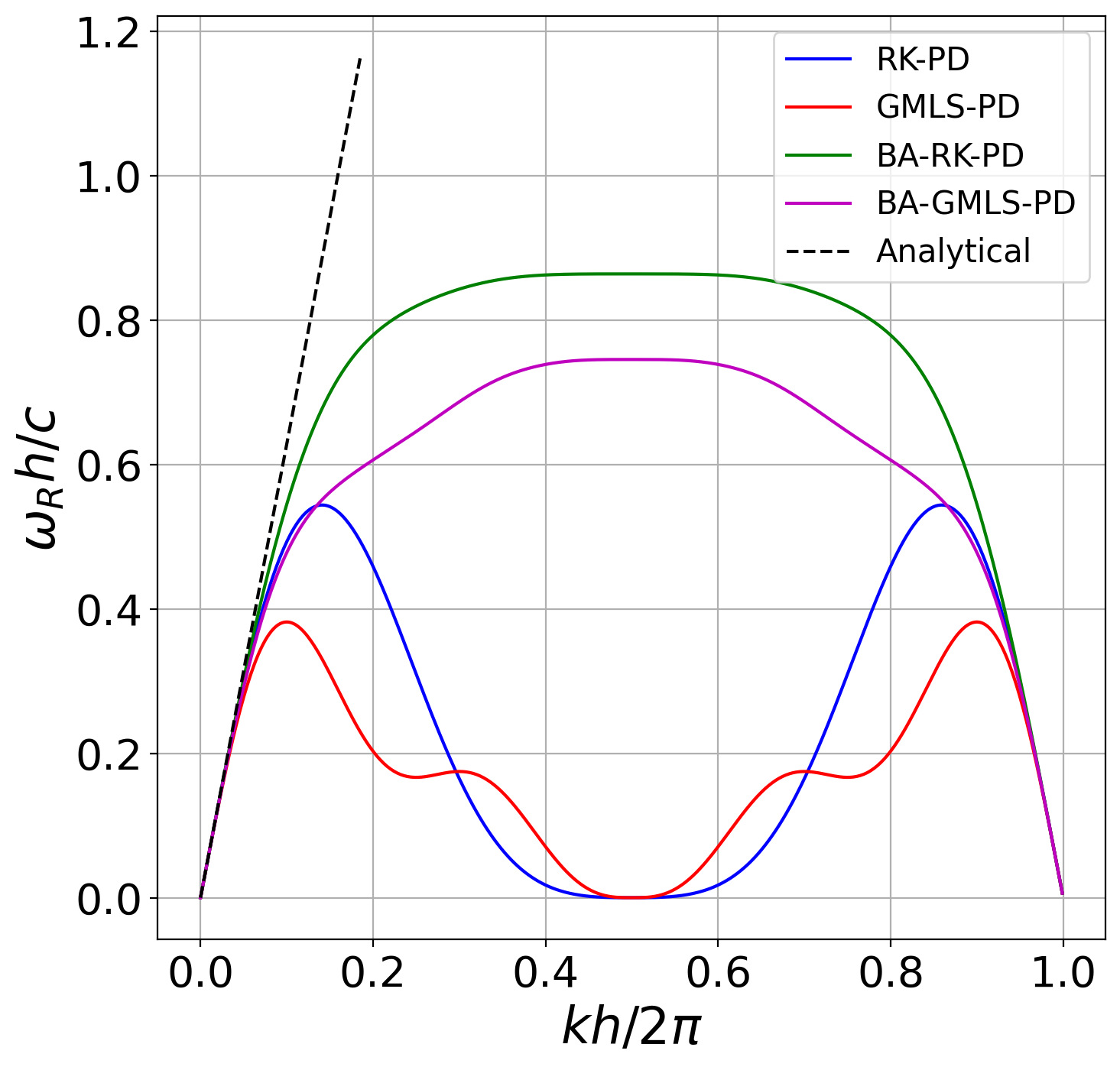}}
  \caption{Wave dispersion diagrams for 1D wave propagation on a uniform grid for two horizon sizes. Dashed lines correspond to the classical (local) linear elastodynamics solution.}
  \label{fig:wave-uniform}
\end{figure*}

\cref{fig:wave-uniform} shows the normalized real part of $\omega$ plotted versus the normalized $k$ (number of wavelengths per $h$) for the uniform grid. RK-PD and BA-RK-PD denote the reproducing kernel peridynamics and its bond-associated variant, respectively. GMLS-PD and BA-GMLS-PD refer to the generalized moving least square peridynamic method and its bond-associated version, respectively. For the uniform grid, due to symmetry, the imaginary part of $\omega$ is zero. All the models capture the analytical response as $kh\rightarrow0$. Wave dispersion occurs for higher wave numbers, as expected. It is observed that the bond-associated variants are significantly less dispersive than the base formulations. Moreover, the number of zero angular frequencies (for a non-zero wavelength), which is indicative of instability \citep{belytschko2000unified}, is decreased using the bond-associative models. Comparing the results between the two horizon sizes, it is seen that using a larger $\delta$ for a fixed $h$ leads to more dispersion, which is also expected.

\begin{figure*}[!ht]
  \centering
  \subfloat[][$\delta = 3 \, h$]{\includegraphics[height=0.4\textwidth]{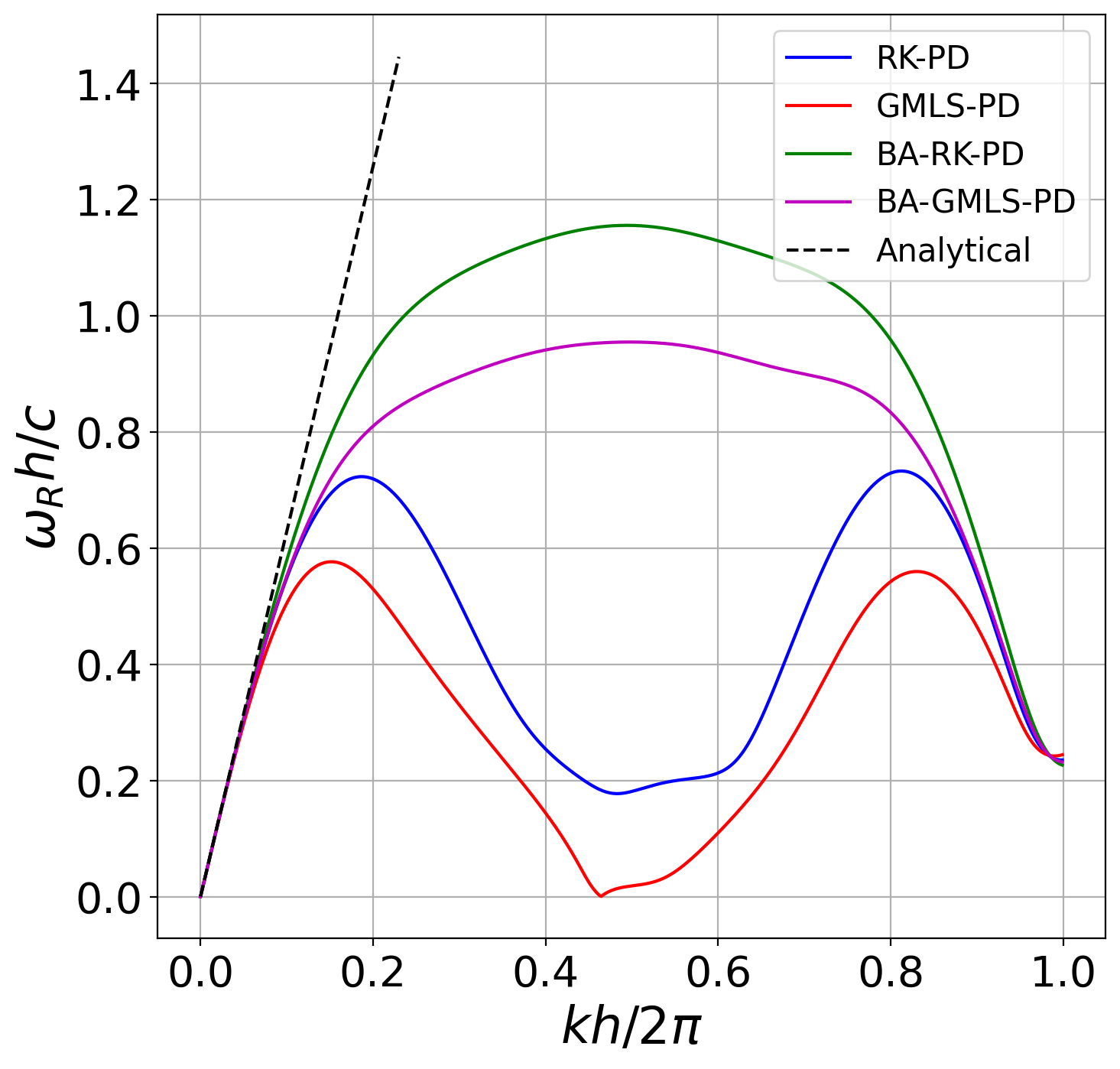}}
  \hspace*{0.5cm}
  \subfloat[][$\delta = 3 \, h$]{\includegraphics[height=0.4\textwidth]{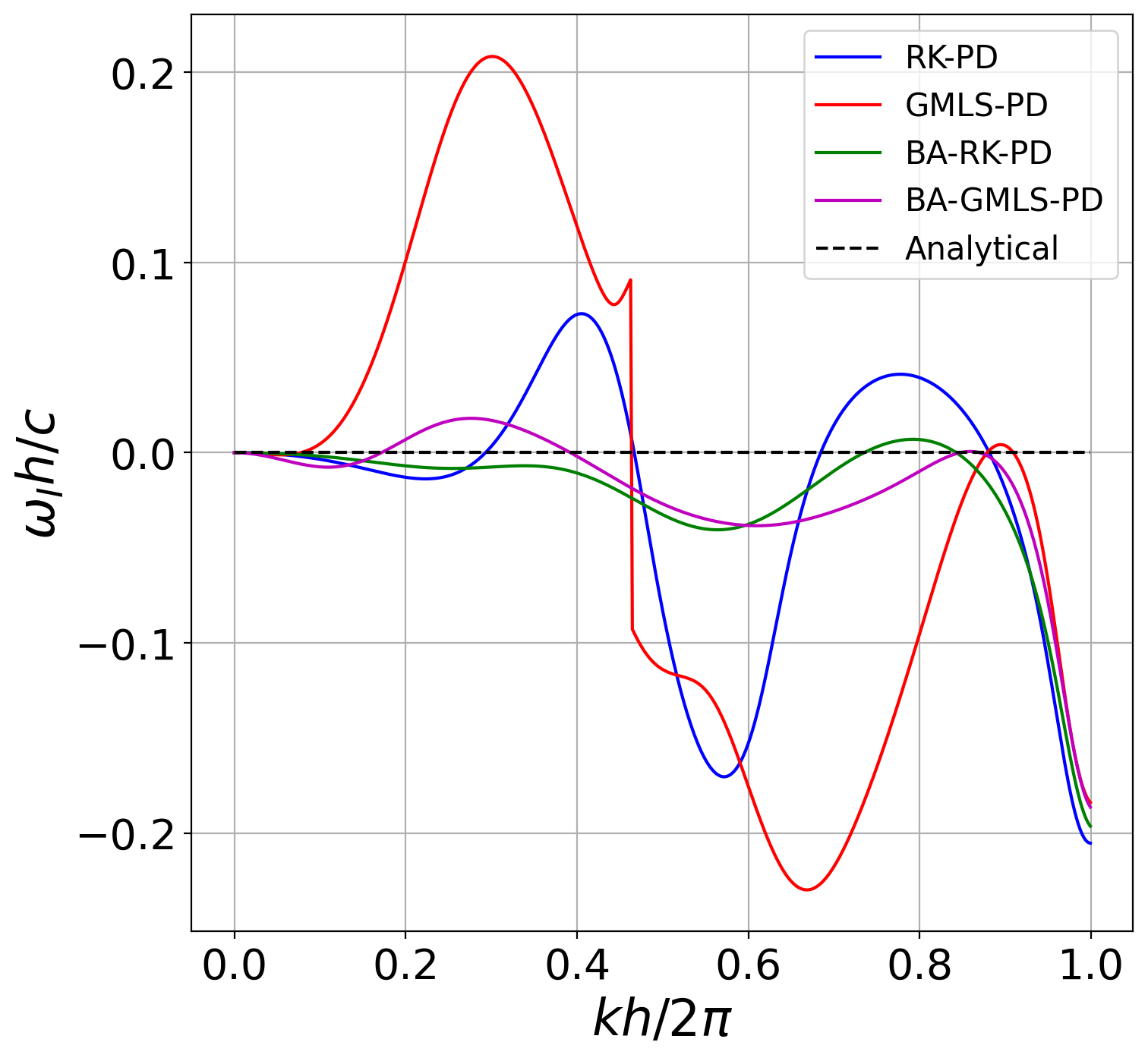}}

  \subfloat[][$\delta = 4 \, h$]{\includegraphics[height=0.4\textwidth]{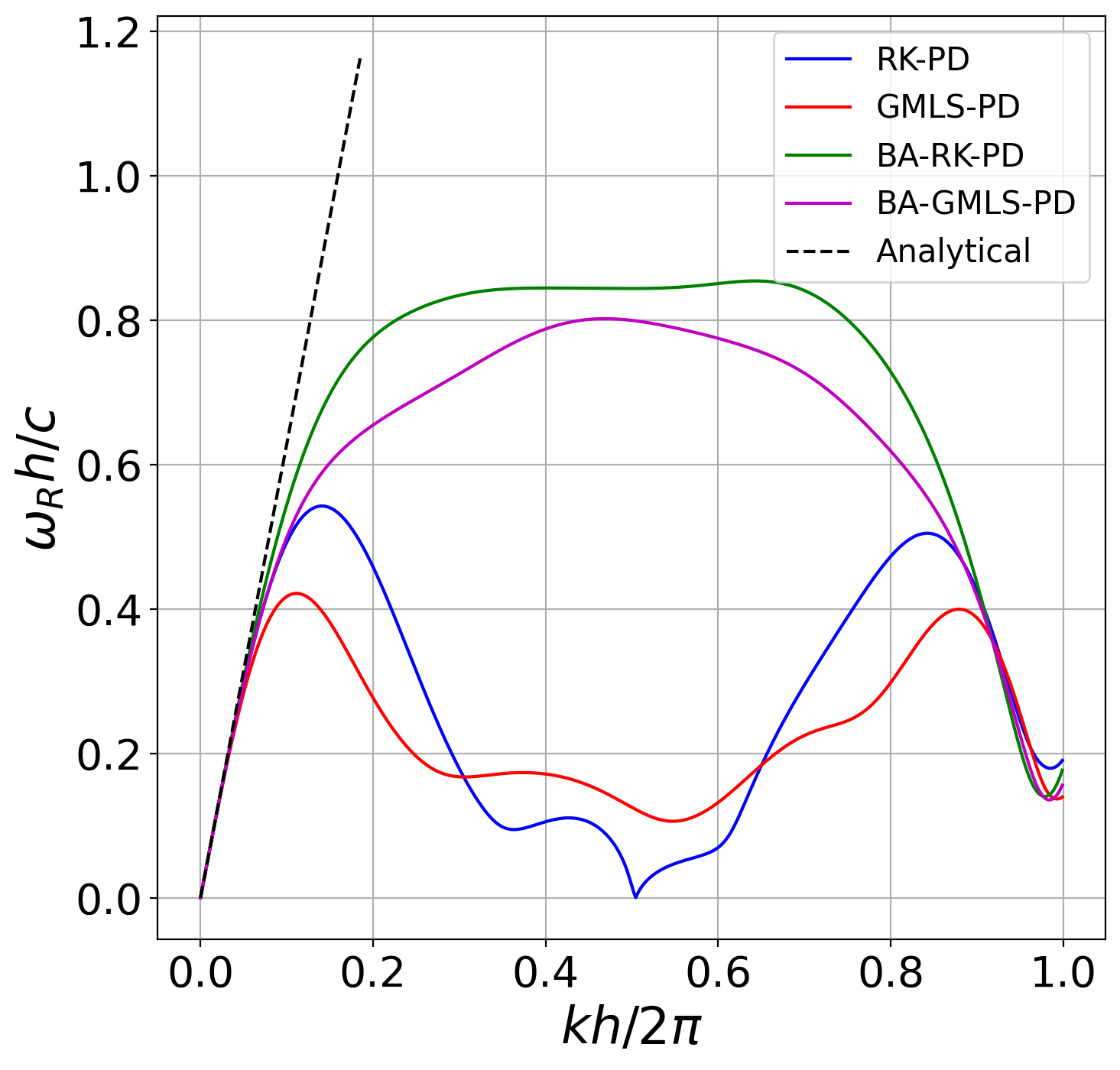}}
  \hspace*{0.5cm}
  \subfloat[][$\delta = 4 \, h$]{\includegraphics[height=0.4\textwidth]{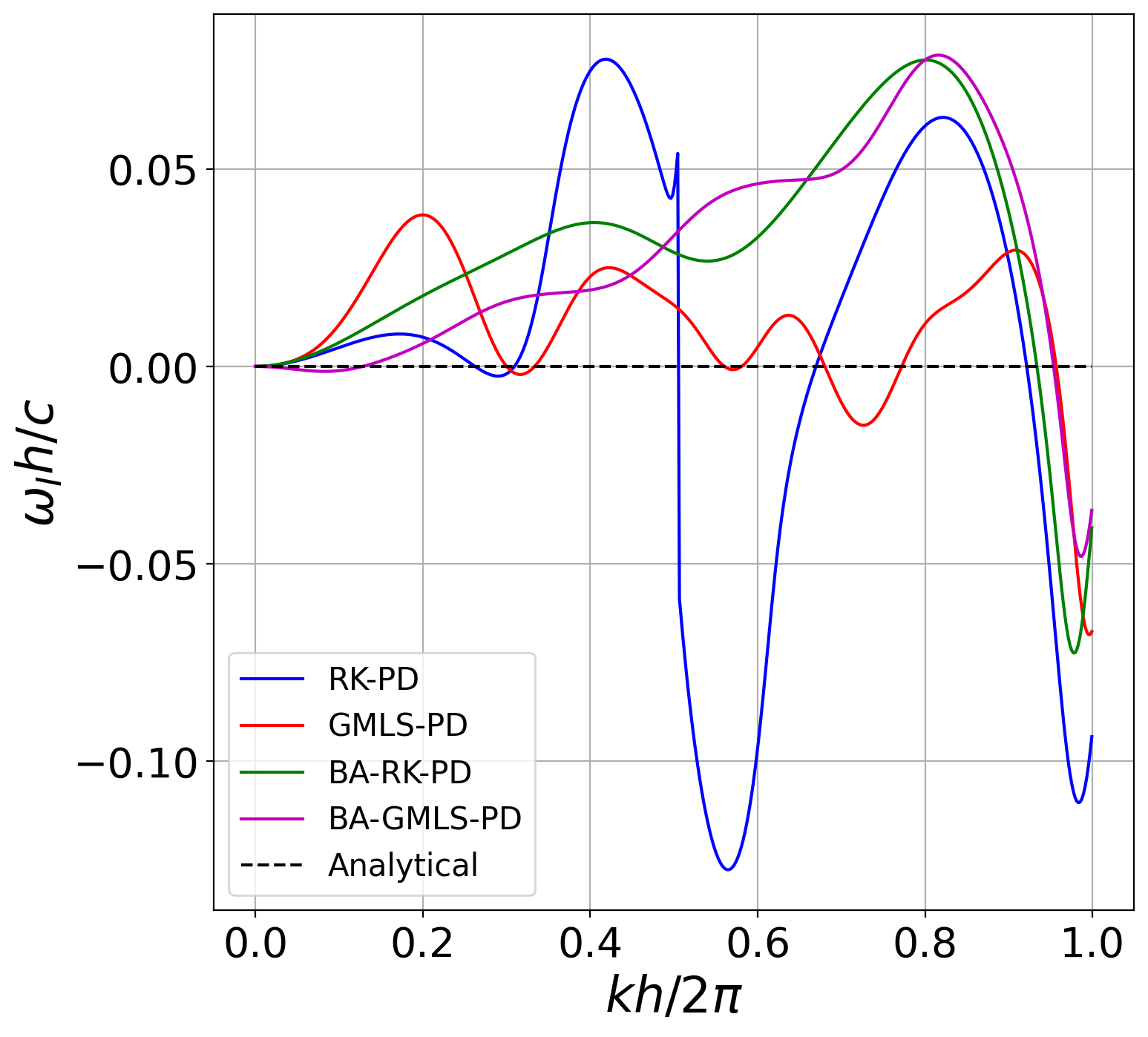}}
  \caption{Wave dispersion diagrams for 1D wave propagation on a non-uniform grid for two horizon sizes. Real component is shown in (a) and (c). Imaginary component is plotted in (b) and (d). Dashed lines correspond to the classical (local) linear elastodynamics solution.}
  \label{fig:wave-non-uniform}
\end{figure*}

\cref{fig:wave-non-uniform} shows the normalized real and imaginary components of $\omega$ plotted versus the normalized $k$ for the non-uniform grid. As for the uniform-grid case, a local solution is captured for small $kh$ regardless of the formulation chosen. However, for larger $kh$, the bond-associated formulations are significantly less dispersive for both the real and imaginary components of $\omega$. Using a larger $\delta/h$ leads to more wave dispersion. We also observed that for the bond-associated variants, for small to moderate $kh$ ($kh/2\pi < 0.5$), the ratio of Im($\omega$) to Re($\omega$) is small, i.e., Im($\omega$)/Re($\omega$) < 0.05, which is not at all the case for the base models. This is yet another manifestation of increased stability and robustness with respect to mesh distortion of the bond-associative approaches relative to the base formulations. 

Only 1D studies were presented in this section. We expect to obtain similar results in higher dimensions because \citet{bazant2016wave} showed that 1D wave dispersion results are comparable to their 2D counterparts for the correspondence PD framework. Although the present results are very encouraging, further studies involving dynamic and nonlinear calculations, in multiple dimensions, using bond-associated formulations are needed to make definitive conclusions about their performance in this regime.

In addition, \citet{bazant2016wave} compared the bond-based PD wave dispersion with more classical non-local counterparts, and concluded that bond-based PD, and not state-based PD, behaves more similarly to the classical non-local case showing (significantly) less wave dispersion. The present results clearly show that it is not the choice of modeling, but rather an inherent numerical instability of the state-based framework that led to the undesired wave dispersion characteristics. This instability is circumvented by the proposed bond-associated formulation, rendering the more versatile state-based models stable and accurate for the approximation of wave dispersion phenomena.

\section{Natural Boundary Conditions}
\label{sec:natural}

A novel approach is given in this section for non-local modeling of natural boundary conditions in peridynamics. Correspondence-model based PD is employed, which enables a direct use of classical stress measures in PD. 

\begin{figure*}[!htbp]
  \centering
  \includegraphics[width=0.9\textwidth,trim={0.2cm 1.5cm 0.5cm 1.5cm },clip]{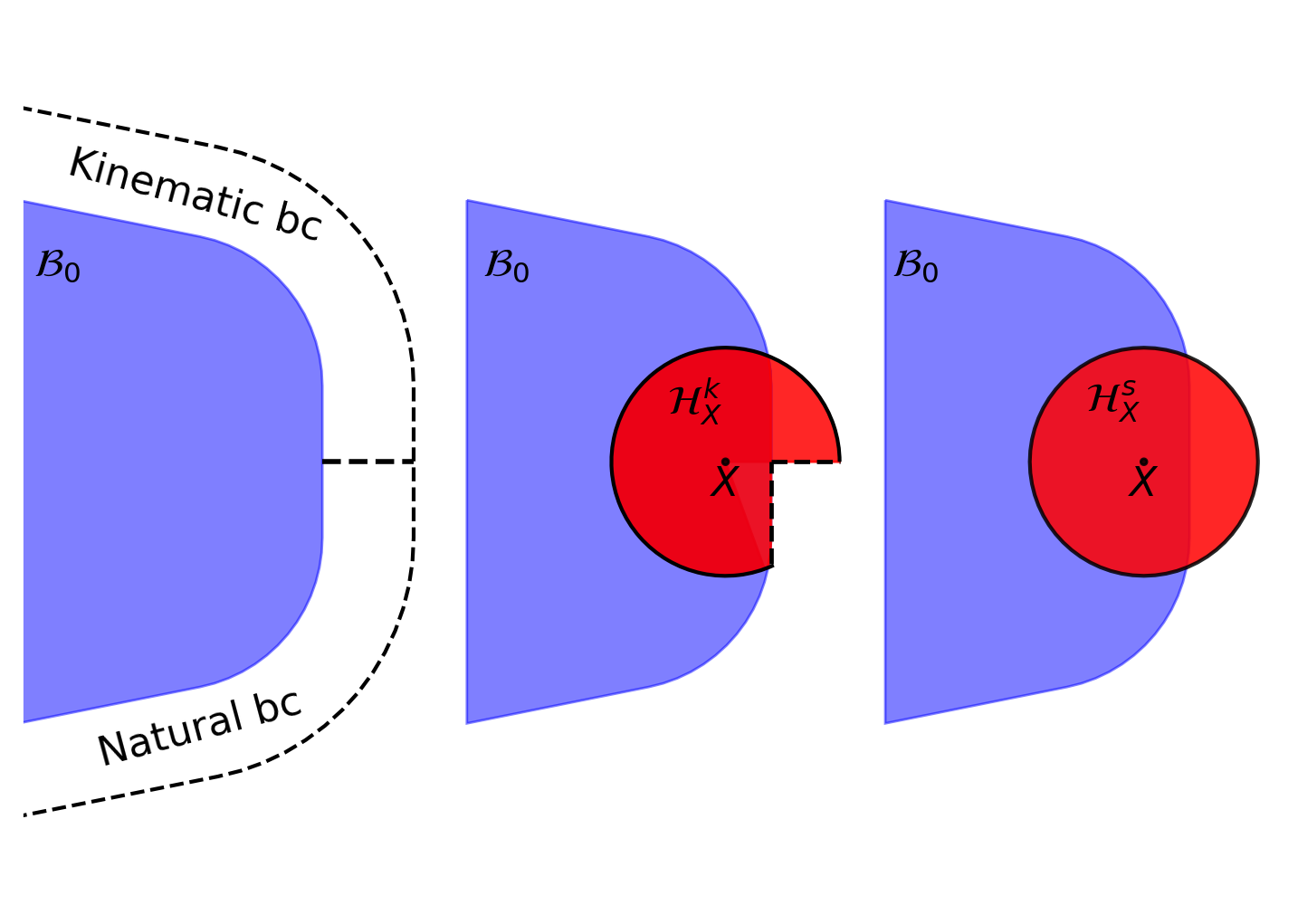}
  \caption{Two different neighborhoods used to evaluate the deformation gradient and stress divergence. Natural-bc neighbors are only considered for evaluating the divergence term.}
  \label{fig:neighborhoods}
\end{figure*}
The common PD practice involves applying boundary conditions on a non-local collar (with size $\delta$) surrounding the body, to be consistent with the non-local nature of the theory \citep{zhou2010mathematical}. Here we assume that either essential (kinematic) or natural (stress) boundary conditions are prescribed at each boundary point. This also means there is no kinematic information for the boundary points with natural conditions (called natural-bc nodes in the remainder of the text). The natural-bc nodes, therefore, cannot be involved in the evaluation of their neighbors' deformation gradient, for which the knowledge of the kinematics (i.e., displacement) is required. By this argument, we use two different set of families for material points in the vicinity of boundaries (including damaged surfaces): one set, denoted as $\mc{H}^k$, includes only the non-broken bonds from the bulk and the essential-bc nodes to compute the deformation gradient; another set denoted as $\mc{H}^s$ involves all the surrounding nodes to compute the stress divergence. \cref{fig:neighborhoods} provides a schematic of the modeling approach.

To incorporate this two-neighborhood approach, \cref{eqn:F} is modified as
\begin{align}
  \mbt{F}_I = \mbf{I} + \sum_{J\in\mc{H}^k_I} \left[\mbf{u}_{J} - \mbf{u}_I\right] \, \bs\gamma_{IJ}^\intercal ,
  \label{eqn:Hk-F}
\end{align}
and \cref{eqn:BA-div} is adapted as 
\begin{align}
  \left(\nabla_h \cdot \mbt{P}\right)_I = \sum_{J\in\mc{H}^s_I} \left[\mbt{P}_{JI} - \mbt{P}_I\right] \, \bs\gamma_{IJ} ,
  \label{eqn:Hs-div}
\end{align}
where
\begin{align*}
  \mbt{P}_{JI} = 
  \begin{cases}
    \mbf{P}(\mbt{F}_{JI}) ~~~ & \text{if $JI$ is non-broken bond with kinematics,} \\
    \hat{\mbf{P}}_J ~~~ & \text{if $J$ is natural-bc node,} \\
    \mbf{0} ~~~ & \text{if $JI$ is broken bond.} 
  \end{cases}
\end{align*}
Here, $\hat{\mbf{P}}_J$ is the stress tensor at $J$ that is directly supplied through boundary conditions. Computation of $\mbt{F}_{JI}$, using \cref{eqn:BA-F}, is only needed if $JI$ is a non-broken bond with a kinematic condition.

\subsection{Numerical results}
\label{sec:numerics}

In this section, the reproducing kernel peridynamics and its bond-associated counterpart are called RK-PD and BA-RK-PD, respectively. The generalized moving least square peridynamic method and its bond-associated variant are denoted as GMLS-PD and BA-GMLS-PD, respectively. A linear bond-associated model (without the higher-order corrections), equivalent to the linear BA-RK-PD, is called BA-PD. The following numerical examples are provided for a threefold goal: (1) to show the suitability of the developed methodology for modeling natural boundary conditions, (2) to test the convergence of the different versions of the correspondence formulation to linear elastics, (3) to demonstrate the necessity of both elements, i.e., bond-associative and higher-order corrections, to obtain a robust formulation.

We consider linear ($n=1$), quadratic ($n=2$), and cubic ($n=3$) models in a set of two-dimensional (plane-strain) static problems. The horizon $\delta$ is chosen according to the nodal spacing and the order of the formulation ($n$), in each problem, to ensure having sufficient number of neighbors to obtain polynomial unisolvency. Boundary conditions (essential or natural) are prescribed on a fictitious layer of size $\delta$, and elastostatic solutions are computed.

Convergence of the models to linear elastostatics is studied in the following examples. To recover a local solution, the size of the horizon and the mesh approach zero concurrently. As the non-local deformation gradient, i.e. \cref{eqn:Hk-F}, is a linear function of the displacement, the following constitutive law governs a linear formulation
$$ \bs\eps = (\mbf{F}+\mbf{F}^T)/2 - \mbf{I} , $$
$$ \mbf{P} = \lambda \ {\rm tr}(\bs\eps) \, \mbf{I} + 2 \, \mu \, \bs\eps, $$
where $\lambda$ and $\mu$ are the Lam\'e parameters.

A cubic B-spline kernel, commonly used as a smoothing function in meshfree methods, is utilized as the influence function for RK-PD (and its BA variant). This radial function is given as:
\begin{align*}
  \alpha(\hat\xi) = 
  \begin{cases}
    \dfrac{2}{3} - 4\,\hat\xi^2 + 4\hat\xi^3  ~~~ & \text{for} ~~~ 0 < \hat\xi \leq \dfrac{1}{2} \\
    \dfrac{4}{3} - 4\,\hat\xi + 4\,\hat\xi^2 - \dfrac{4}{3}\hat\xi^3  ~~~ & \text{for} ~~~ \dfrac{1}{2} < \hat\xi \leq 1 \\
    0  ~~~ & \text{otherwise} 
  \end{cases}
  , 
\end{align*}
where
$$\hat\xi \equiv \frac{|\mbf{X}_J-\mbf{X}_I|}{\delta} . $$
As noted in Part I of this paper, the influence function for GMLS-PD and BA-GMLS-PD is adopted as $\alpha(\bxi) = \frac{1}{|\bxi|^2}$, and is already included in the relevant equations, i.e., \cref{eqn:Hs-div,eqn:Hk-F}. 

Free-surface behavior is modeled as a natural boundary condition. Fictitious nodes are employed in the free space, and their associated bonds are broken, i.e., given zero stress values that directly contribute to the internal-force evaluation. A similar approach may be implemented for fracture modeling, i.e., broken bonds carry zero stress tensor. 

We want to emphasize that two set of quadrature weights are required for material nodes with natural-bc neighbors (including broken bonds, such as free-surface neighbors): (1) to obtain a correct gradient for computing the kinematic variable (i.e., deformation gradient), which must always pass a patch test, a set of $\bs\Phi_{IJ}$ (for RK) and $\bs\omega_{IJ}$ (for GMLS) is constructed using only the nodes from the bulk of material and Dirichlet boundary region; and (2) to evaluate the internal force (i.e., divergence of stress), a set of $\bs\Phi_{IJ}$ and $\bs\omega_{IJ}$ is obtained by considering the full neighborhood involving natural-bc nodes (including broken bonds to model surface effects). The natural-bc bonds are directly supplied with stress values (no kinematic variable is computed for the natural-bc bonds).

The error is computed using the root-mean-square (RMS) norm
$$ ||e||_{2} = \sqrt{\frac{\sum_i^N e_i^2}{N}} ,$$
which is equivalent to the full $L_2$ norm for quasi-uniform pointsets \citep{wendland2004scattered}.

\subsubsection{Manufactured solution}
\label{subsec:square}

The following 2D manufactured-solution example is computed on a square domain $[-1,1] \times [-1,1]$. The analytical solution is given by
\begin{align*}
  & u_1(x,y) = A \sin(\pi x/2) \cos(\pi y/2) + B \exp(x) \exp(y) , \notag \\
  & u_2(x,y) = C \cos(\pi x/2) \sin(\pi y/2) + D \exp(x) \exp(y) ,
\end{align*}
resulting in the stress components
\begin{align*}
  P_{11}(x,y) = & \frac{\pi}{2} \Big[ [A+C] \lambda + 2 A \mu \Big] \cos(\pi x/2) \cos(\pi y/2) \notag \\
  & + \Big[ [B+D] \lambda + 2 B \mu \Big] \exp(x) \exp(y) , \notag \\
  P_{12}(x,y) = & \, P_{21}(x,y) = - \frac{\pi}{2} [A+C] \mu \sin(\pi x/2) \sin(\pi y/2) \notag \\
  & + [B+D] \mu \exp(x) \exp(y) , \notag \\
  P_{22}(x,y) = & \frac{\pi}{2} \Big[ [A+C] \lambda + 2 C \mu \Big] \cos(\pi x/2) \cos(\pi y/2) \notag \\
  & + \Big[ [B+D] \lambda + 2 D \mu \Big] \exp(x) \exp(y) .
\end{align*}
%For equilibrium, a body force $\mbf{b}$ is required to satisfy $\mbf{b} + \nabla \cdot \mbf{P} = \mbf{0}$; therefore,
and the corresponding body force $\mbf{b}$ that satisfies static equilibrium, i.e. \cref{eqn:classical},
\begin{align*}
  b_1(x,y) = & - \frac{\pi^2}{4} \Big[ [A+C] \lambda + [3A+C] \mu \Big] \sin(\pi x/2) \cos(\pi y/2) \notag \\
      & + \Big[ [B+D] \lambda + [3B+D] \mu \Big] \exp(x) \exp(y) , \notag \\
  b_2(x,y) = & - \frac{\pi^2}{4} \Big[ [A+C] \lambda + [A+3C] \mu \Big] \cos(\pi x/2) \sin(\pi y/2) \notag \\
      & + \Big[ [B+D] \lambda + [B+3D] \mu \Big] \exp(x) \exp(y) .
\end{align*}
The set of constants is taken as $A=0.2$, $B=-0.15$, $C=-0.15$, and $D=0.1$. We consider a material with Young's modulus $E=100,000$ and Poisson's ratio $\nu=0.3$ in this problem.

\begin{figure*}[!ht]
  \centering
  \subfloat[][Base level (L0)]{\includegraphics[height=0.25\textwidth]{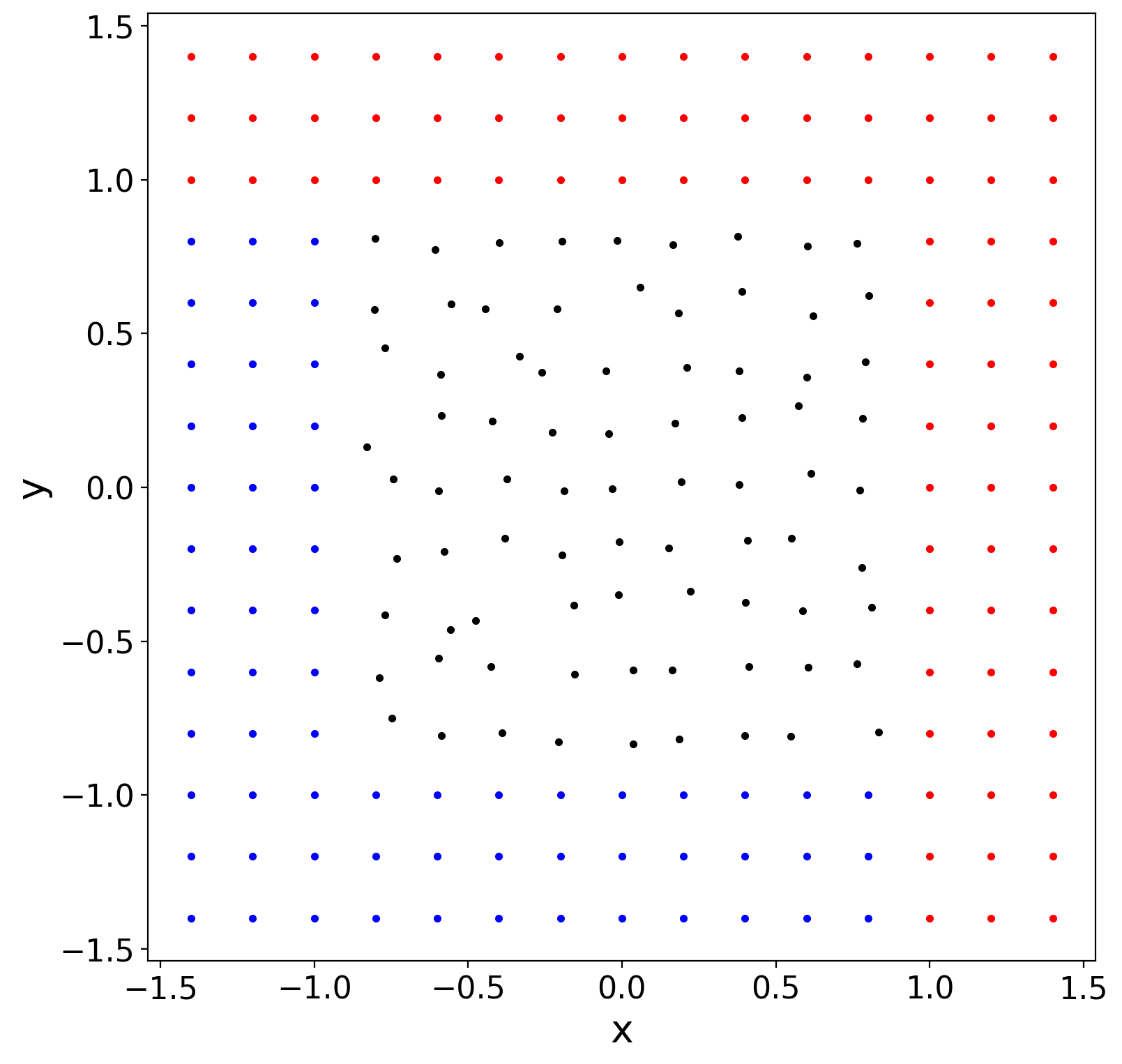}}
  \hspace*{0.1cm}
  \subfloat[][L1]{\includegraphics[height=0.25\textwidth,trim={2.5cm 0 0 0 },clip]{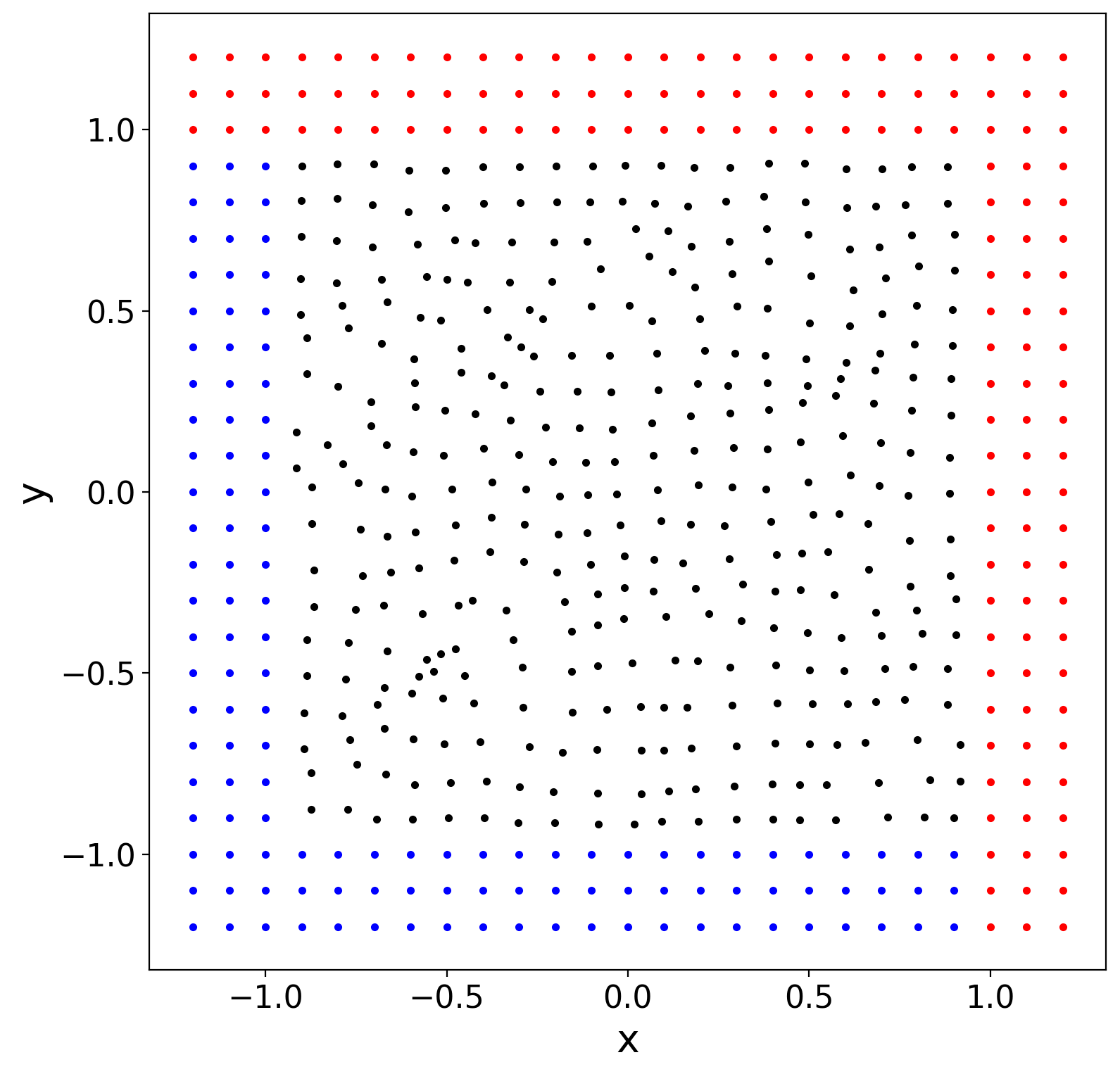}}
  \hspace*{0.1cm}
  \subfloat[][L2]{\includegraphics[height=0.25\textwidth,trim={2.5cm 0 0 0 },clip]{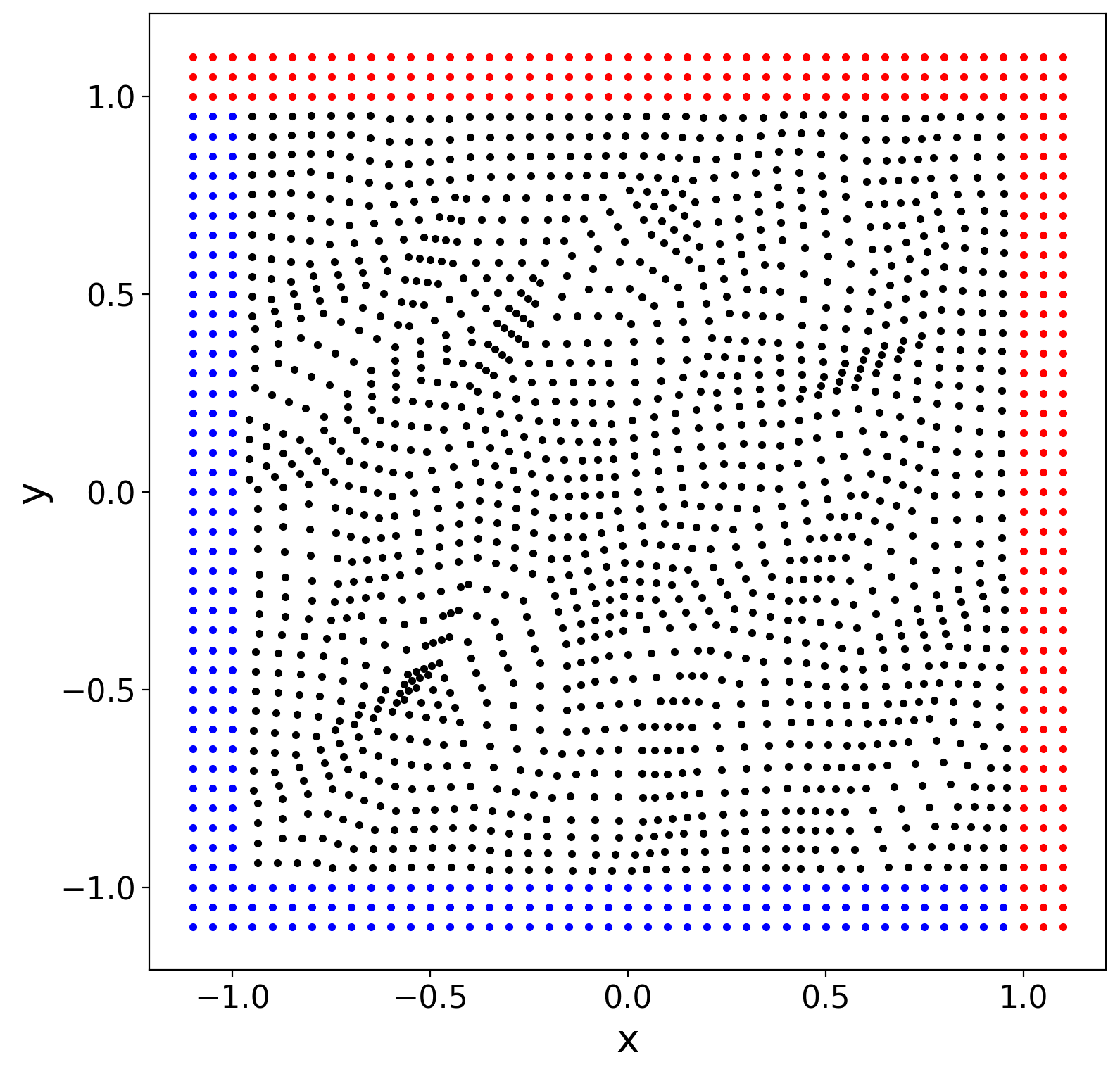}}
  \hspace*{0.1cm}
  \subfloat[][L3]{\includegraphics[height=0.25\textwidth,trim={2.5cm 0 0 0 },clip]{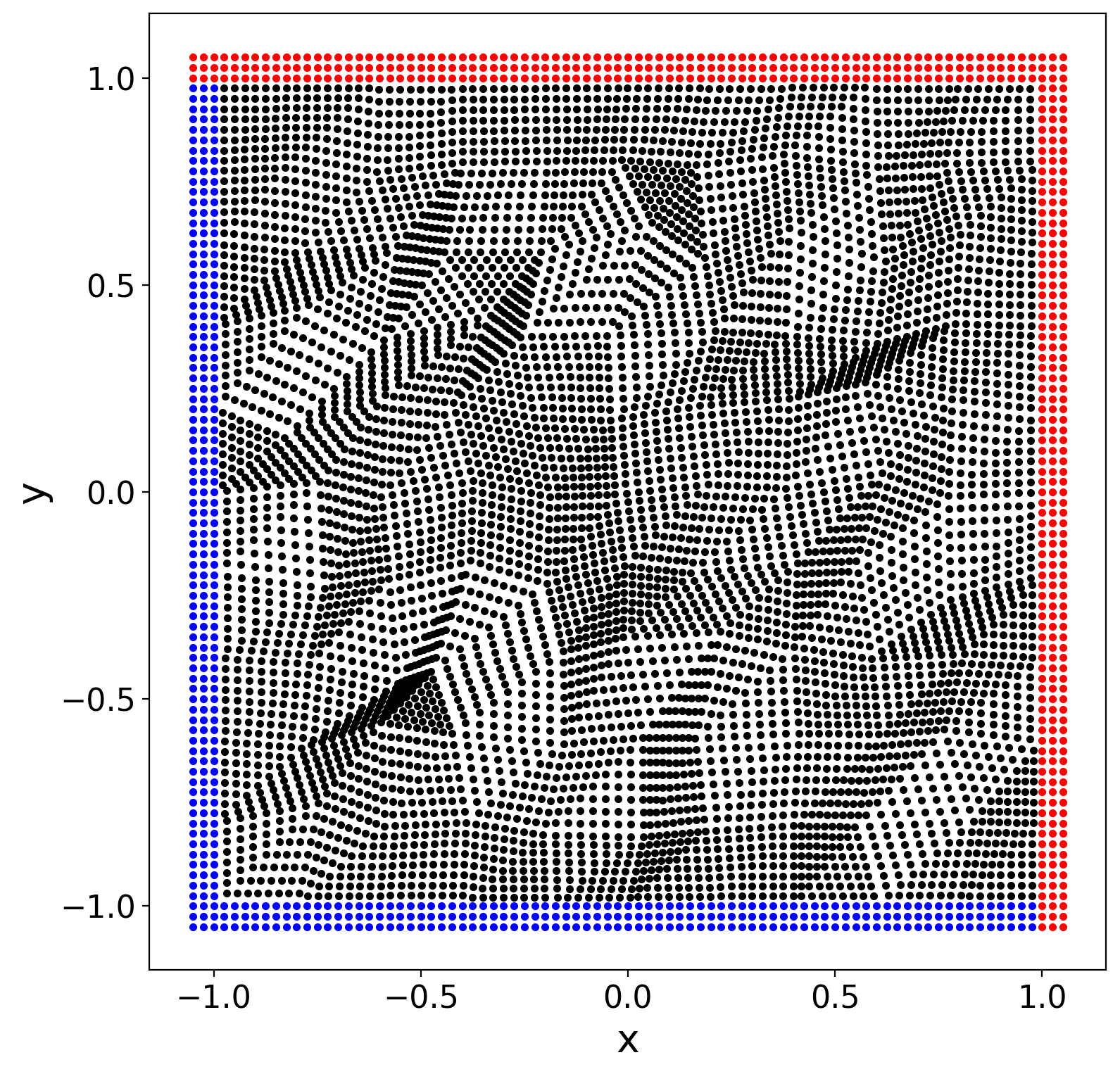}}
  \caption{Non-uniform nodal discretizations and refinements for the 2-D manufactured problem. Blue and red colors denote the nodes used to prescribe displacement-controlled and stress-controlled boundary conditions, respectively.}
  \label{fig:square-mesh}
\end{figure*}

Uniform and non-uniform nodal discretizations are employed in this example. A set of average nodal spacing $h=[0.2, \, 0.1, \, 0.05, \, 0.025]$ is considered for the convergence test. The non-uniform, coarse-level discretization is obtained by perturbing the nodes in the uniform, base level ($h=0.2$) with a random normal distribution with the standard deviation of 0.03 (15\% of the average nodal spacing). The refined levels are constructed in a systematic way by placing new nodes on the midpoints of the adjacent old nodes. The non-uniform nodal discretizations are shown in \cref{fig:square-mesh}. In this problem, natural boundary conditions are applied on the top and right boundary layers, and the bottom and left layers are subjected to essential boundary conditions. For the linear, quadratic, and cubic models, the horizon size is chosen as $\delta = 2.5 h$, $\delta = 3.5 h$, and $\delta = 4.5 h$, respectively. 

\begin{figure*}[!ht]
  \centering
  \subfloat[][Uniform discretization - Linear]{\includegraphics[height=0.39\textwidth]{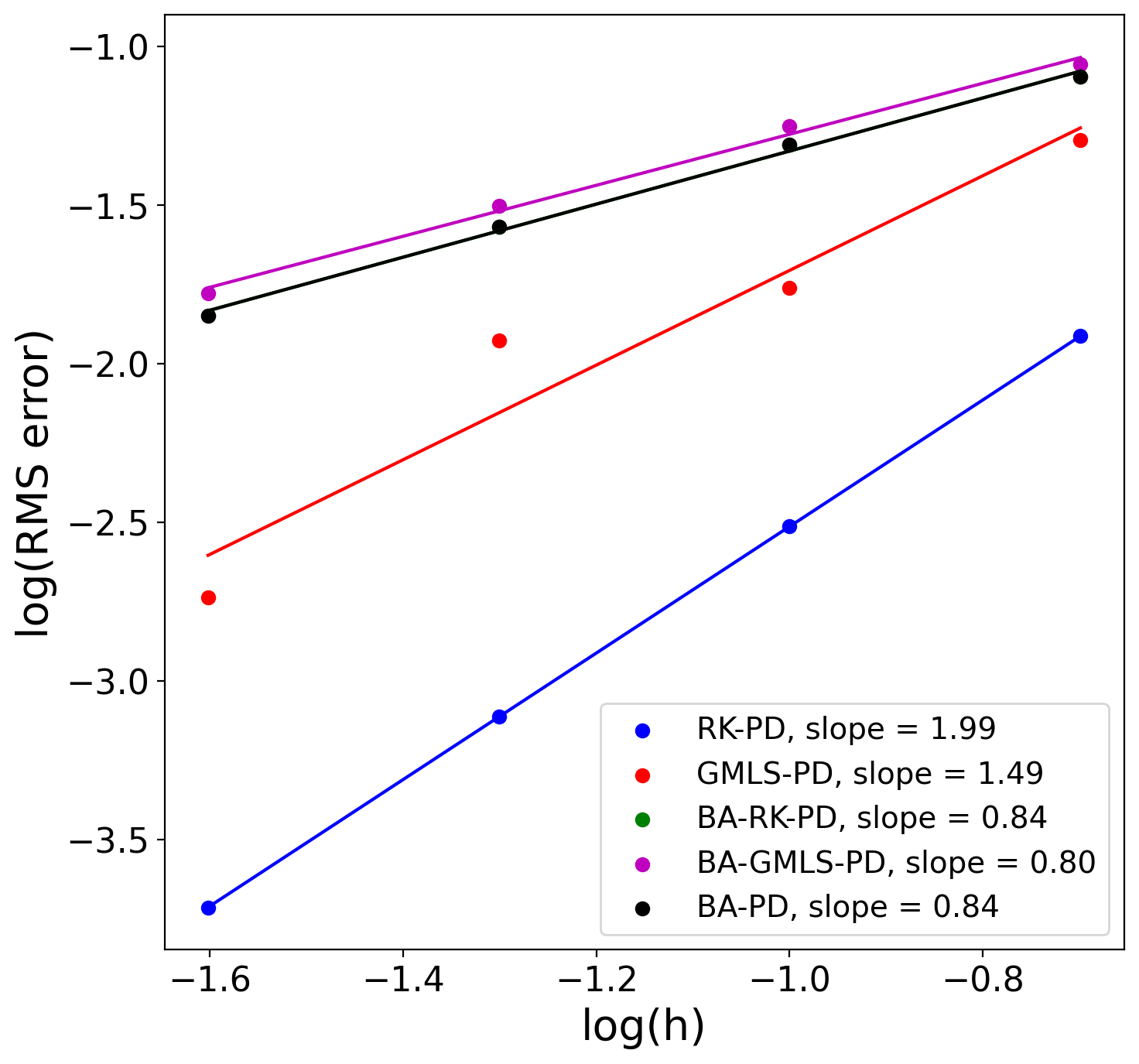}}
  \hspace*{1cm}
  \subfloat[][Non-uniform discretization - Linear]{\includegraphics[height=0.39\textwidth]{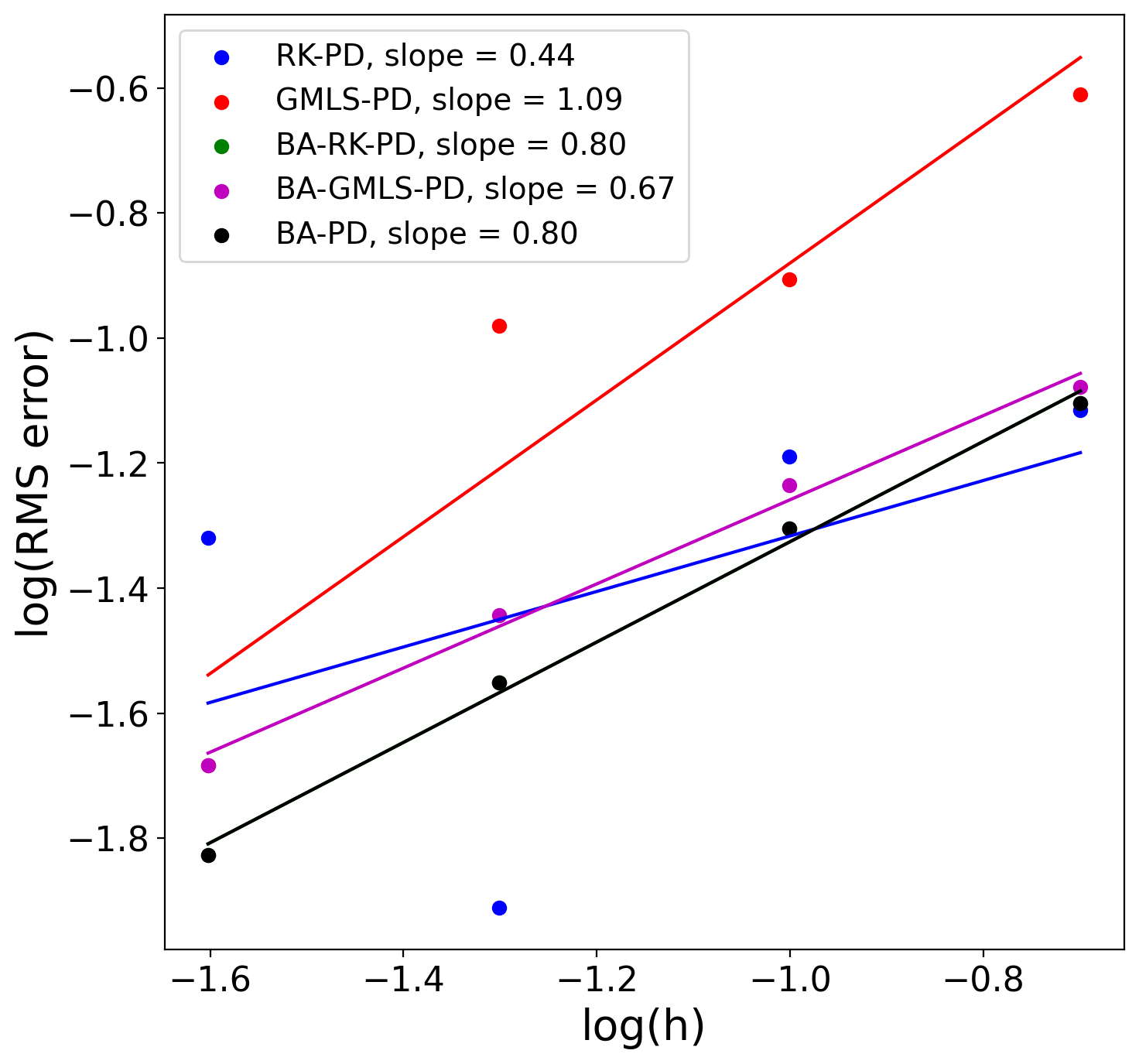}}

  \subfloat[][Uniform discretization - Quadratic]{\includegraphics[height=0.39\textwidth]{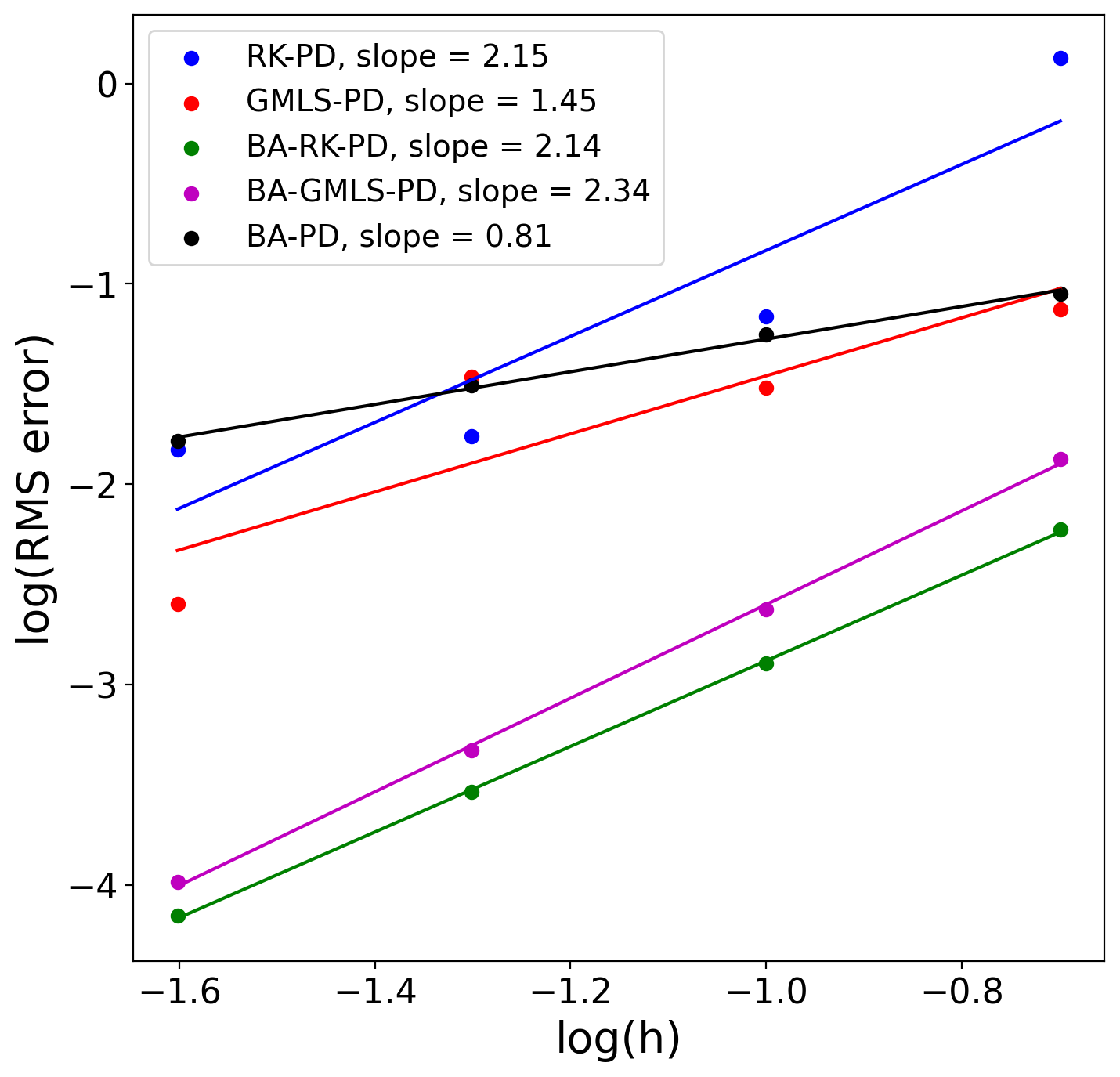}}
  \hspace*{1cm}
  \subfloat[][Non-uniform discretization - Quadratic]{\includegraphics[height=0.39\textwidth]{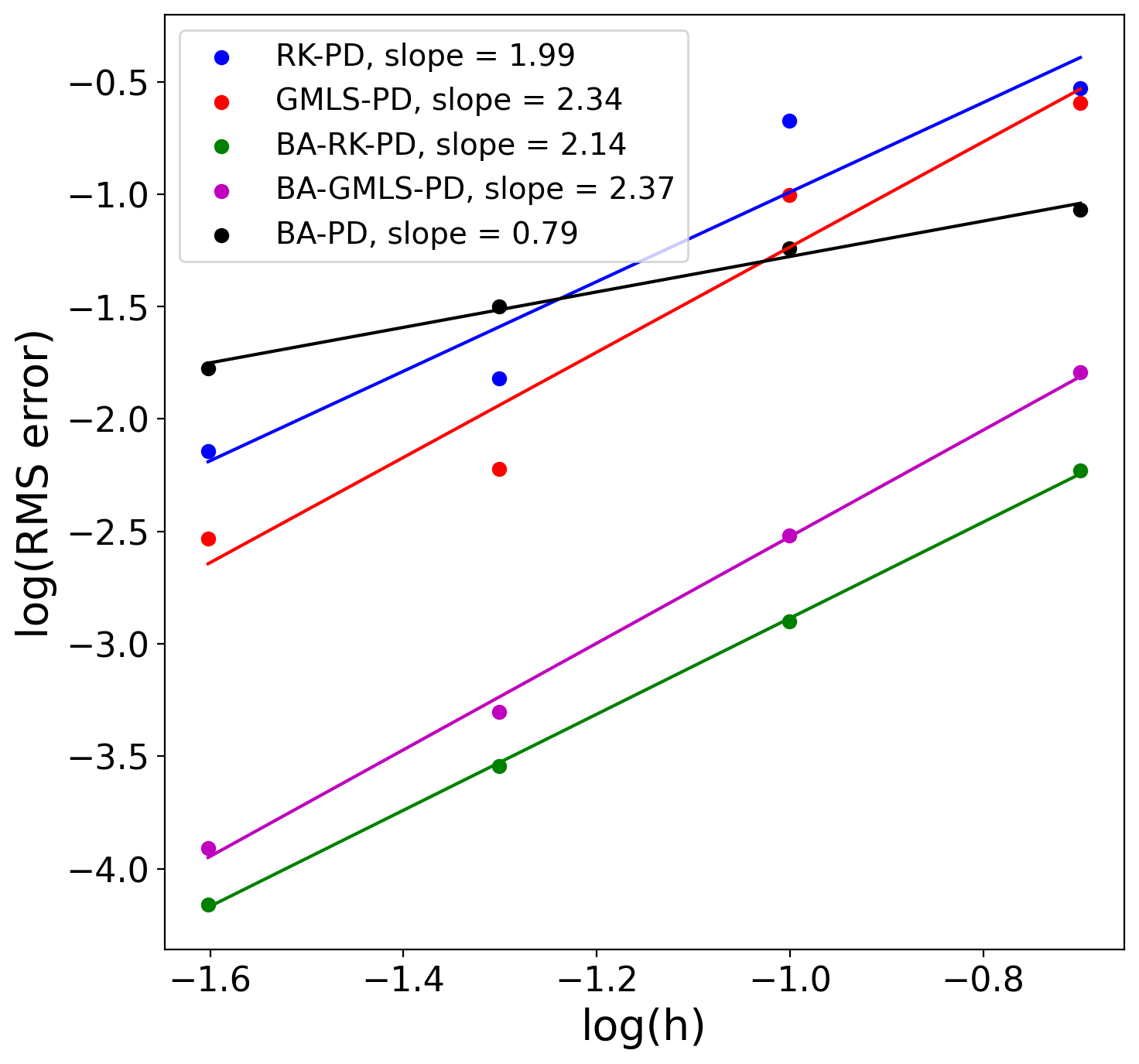}}

  \subfloat[][Uniform discretization - Cubic]{\includegraphics[height=0.39\textwidth]{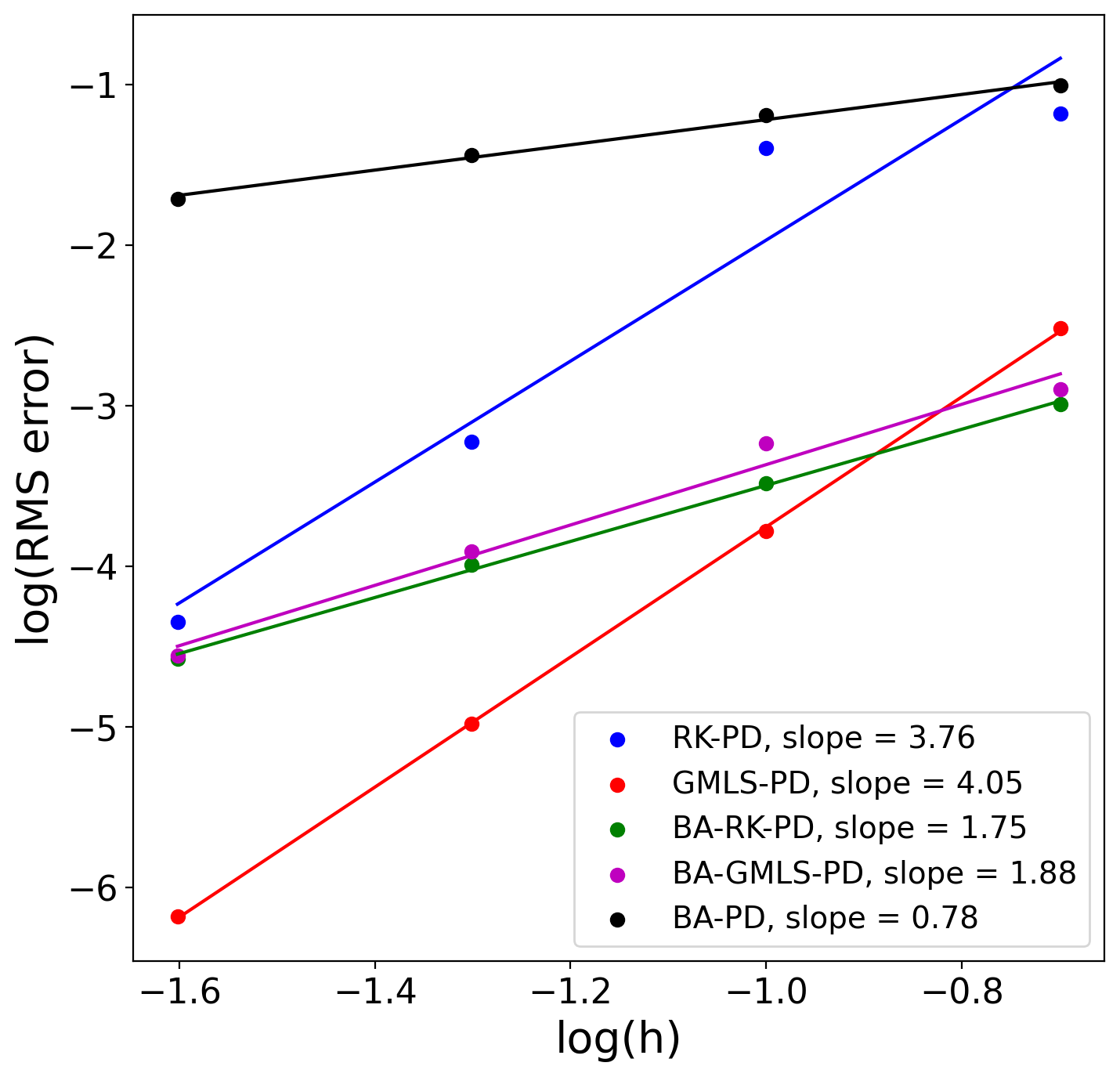}}
  \hspace*{1cm}
  \subfloat[][Non-uniform discretization - Cubic]{\includegraphics[height=0.39\textwidth]{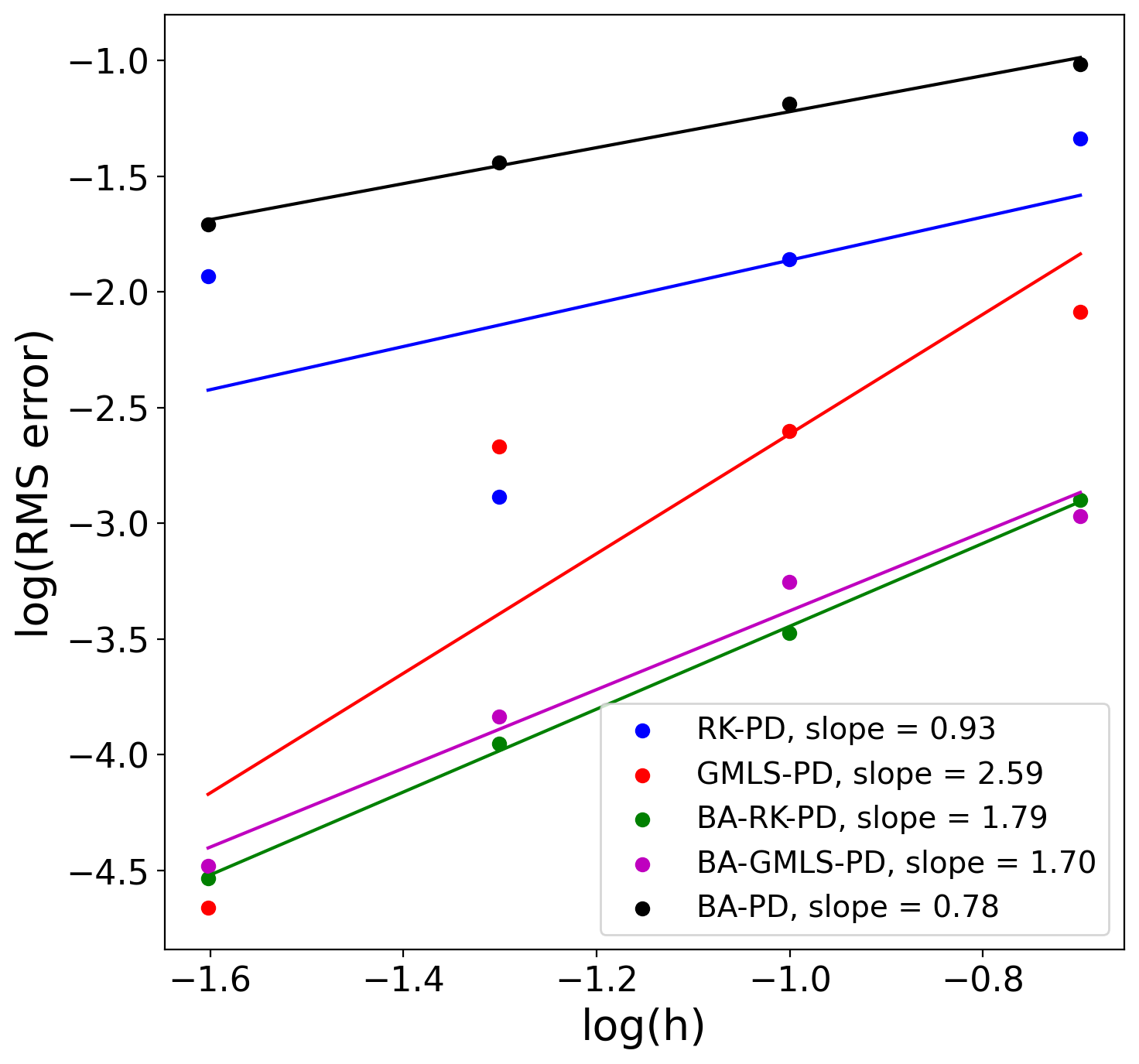}}
  \caption{Convergence of RK-PD, GMLS-PD, BA-RK-PD, BA-GMLS-PD, and BA-PD for the manufactured solution example. Linear, quadratic, and cubic formulations are tested on uniform and non-uniform discretizations. BA-PD and BA-RK-PD overlap for (a--b). Plots show displacement RMS error values.}
  \label{fig:square-convergence}
\end{figure*}

The displacement RMS errors are shown in \cref{fig:square-convergence}. For the uniform grids, the RK-PD and GMLS-PD results show a super-convergence behavior, where a rate of convergence of $n+1$ is empirically observed for odd-order models \citep{leng2019super}. For the noted models, an inconsistent behavior is observed for the non-uniform case, which is rooted in the lack of stability. BA-PD shows convergence with a near linear rate. It is indicated that the higher-order, bond-associated models resulted in a near first-order convergence rate for the linear case and a near second-order convergence rate for the quadratic and cubic formulations, in a problem that involves smooth functions (no discontinuity in the field variables). Only for the higher-order, bond-associated methods, increasing the order of accuracy consistently decreases error. Adding at least a second-order correction with the stabilized model seems necessary to achieve a robust formulation.

\begin{figure*}[!ht]
  \centering
  \subfloat[][RK-PD]{\includegraphics[height=0.21\textwidth]{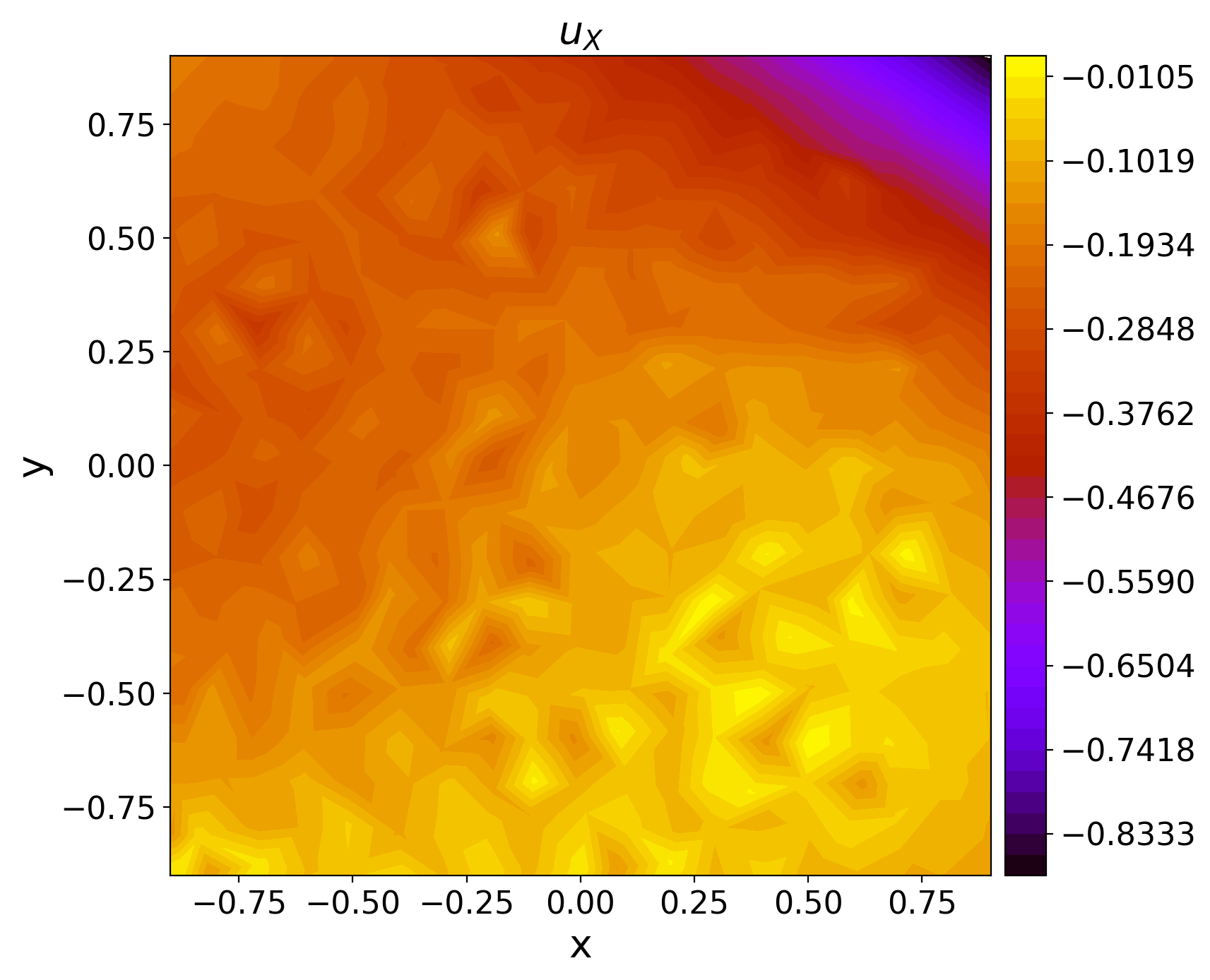}}
  \hspace*{0.2cm}
  \subfloat[][GMLS-PD]{\includegraphics[height=0.21\textwidth,trim={3cm 0 0 0 },clip]{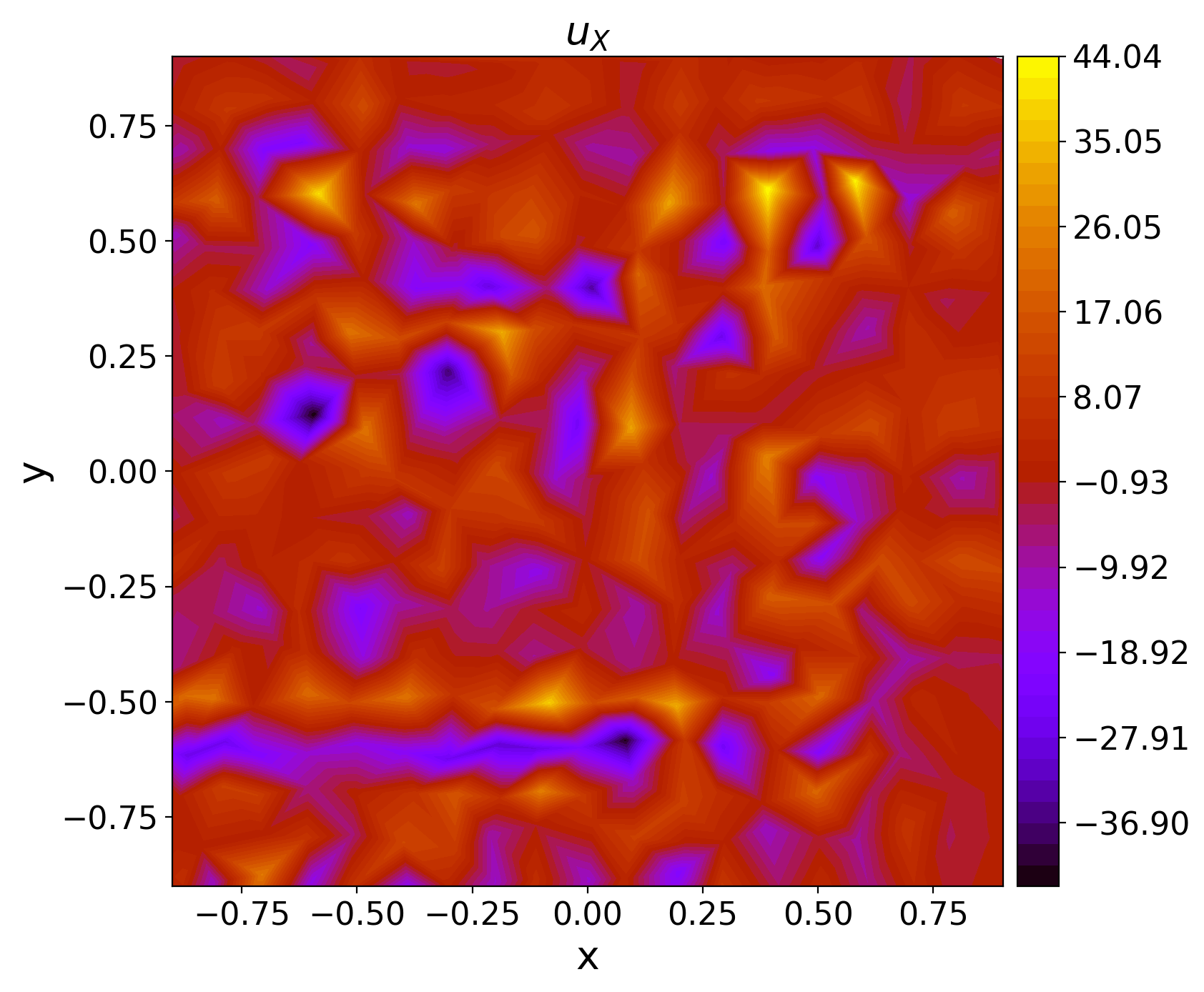}}
  \hspace*{0.2cm}
  \subfloat[][BA-RK-PD]{\includegraphics[height=0.21\textwidth,trim={3cm 0 0 0 },clip]{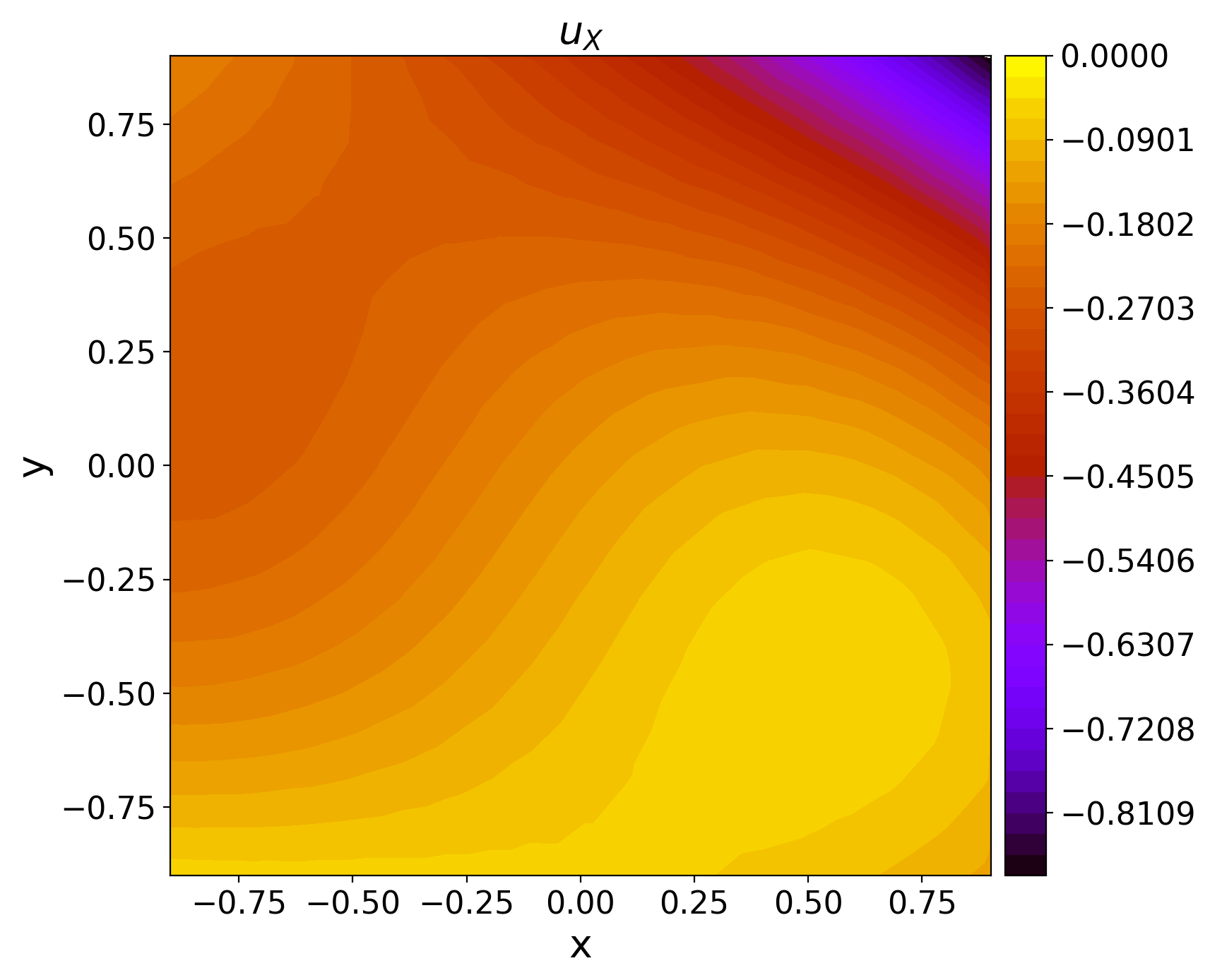}}
  \hspace*{0.2cm}
  \subfloat[][BA-GMLS-PD]{\includegraphics[height=0.21\textwidth,trim={3cm 0 0 0 },clip]{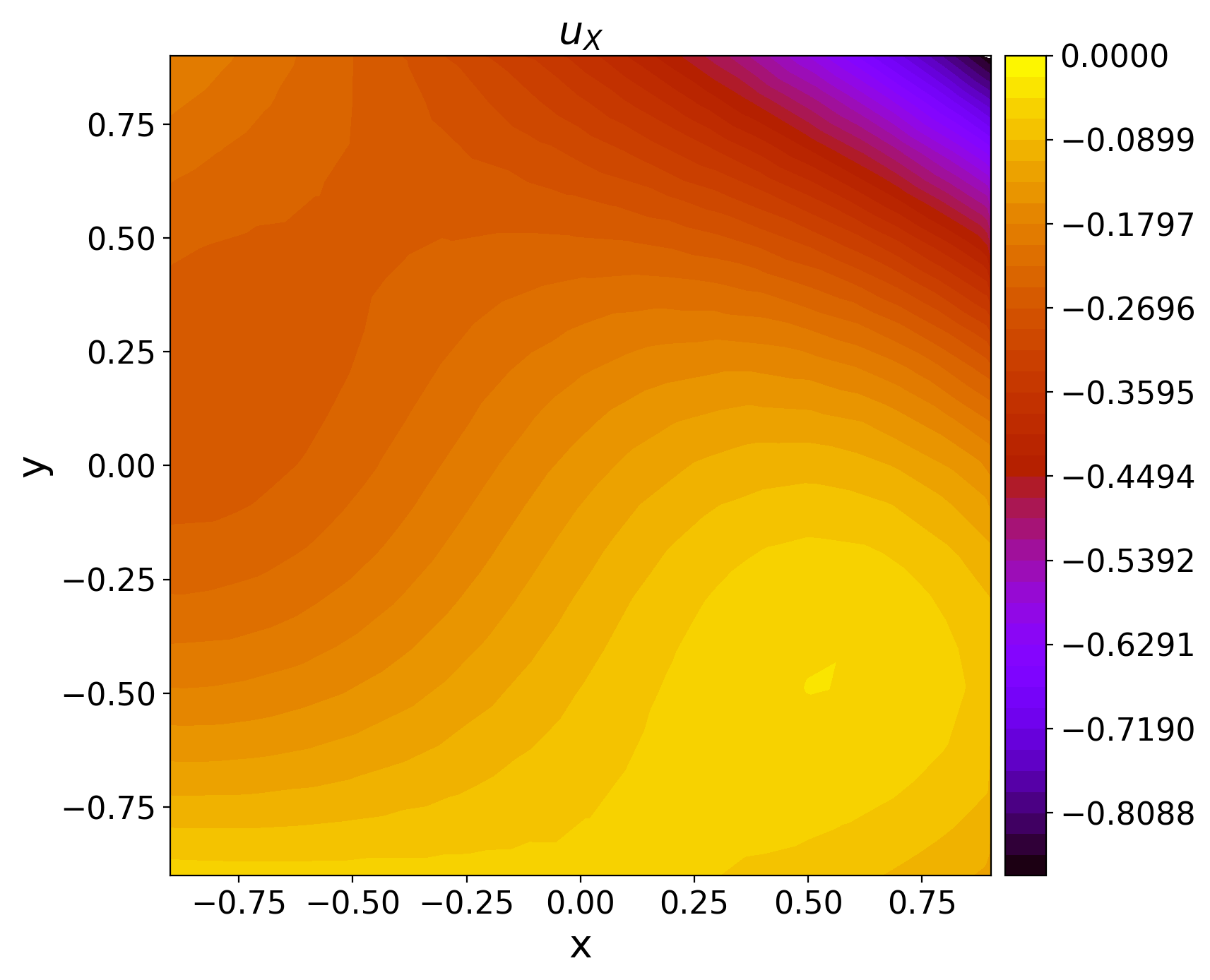}}

  \subfloat[][RK-PD]{\includegraphics[height=0.21\textwidth]{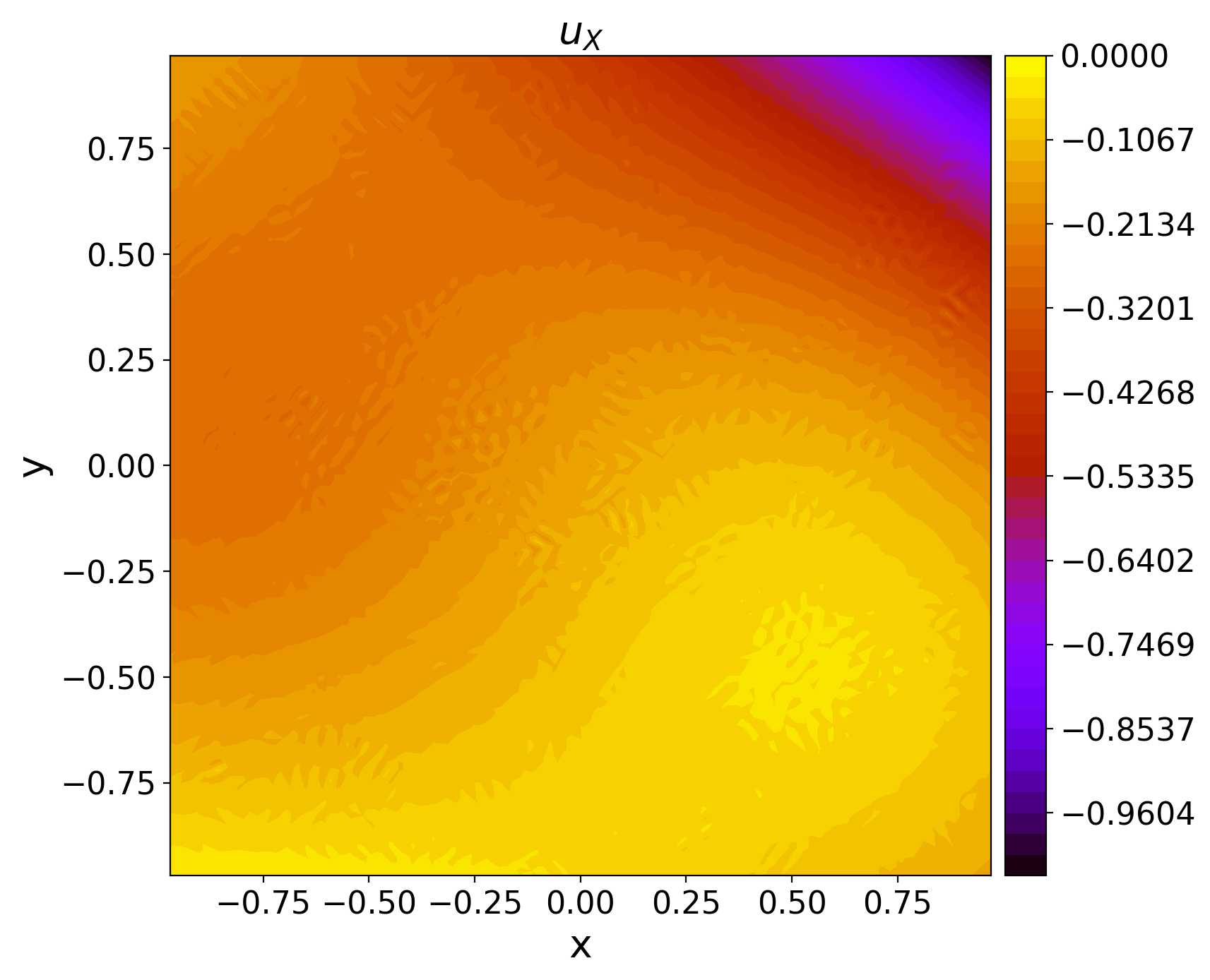}}
  \hspace*{0.2cm}
  \subfloat[][GMLS-PD]{\includegraphics[height=0.21\textwidth,trim={3cm 0 0 0 },clip]{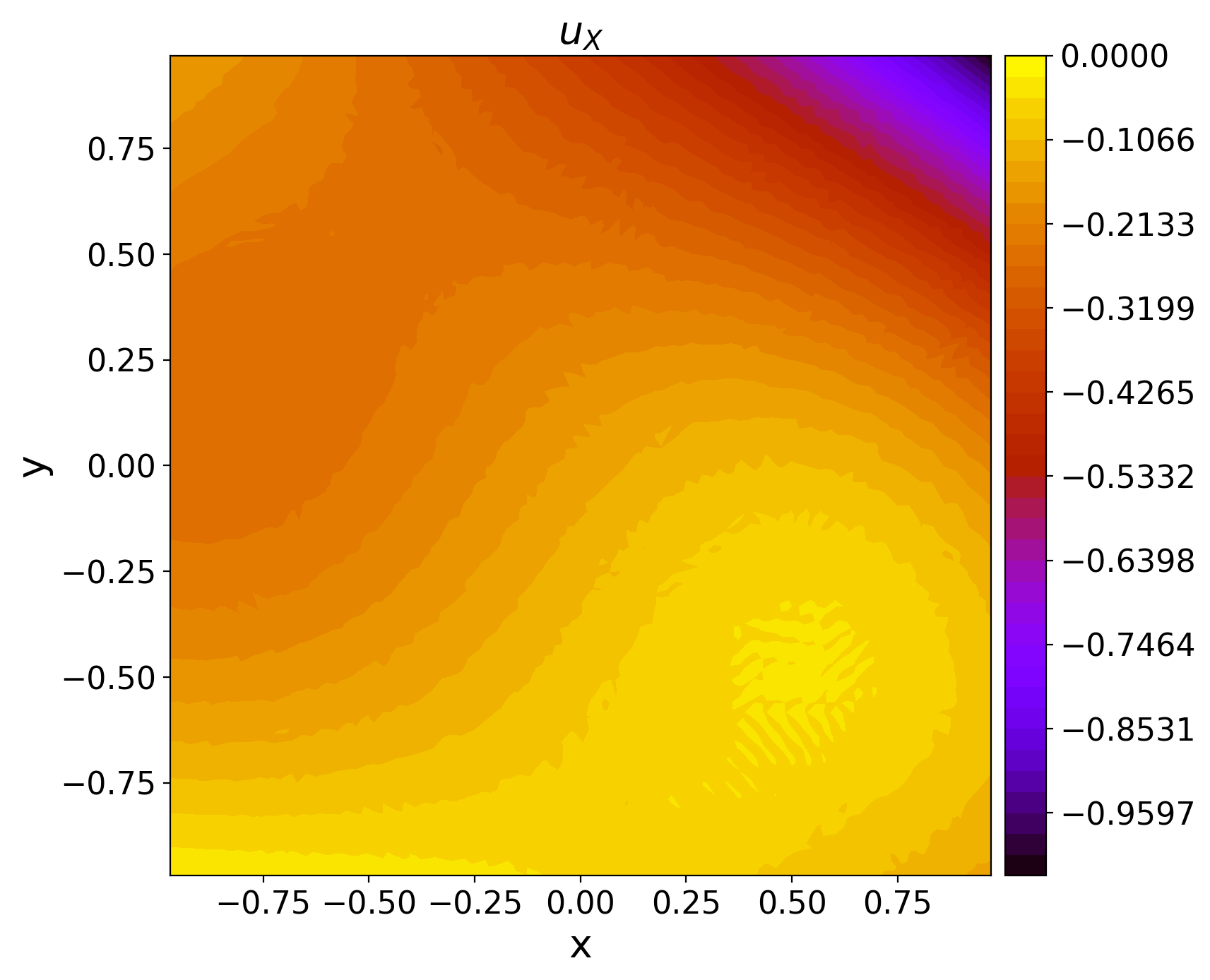}}
  \hspace*{0.2cm}
  \subfloat[][BA-RK-PD]{\includegraphics[height=0.21\textwidth,trim={3cm 0 0 0 },clip]{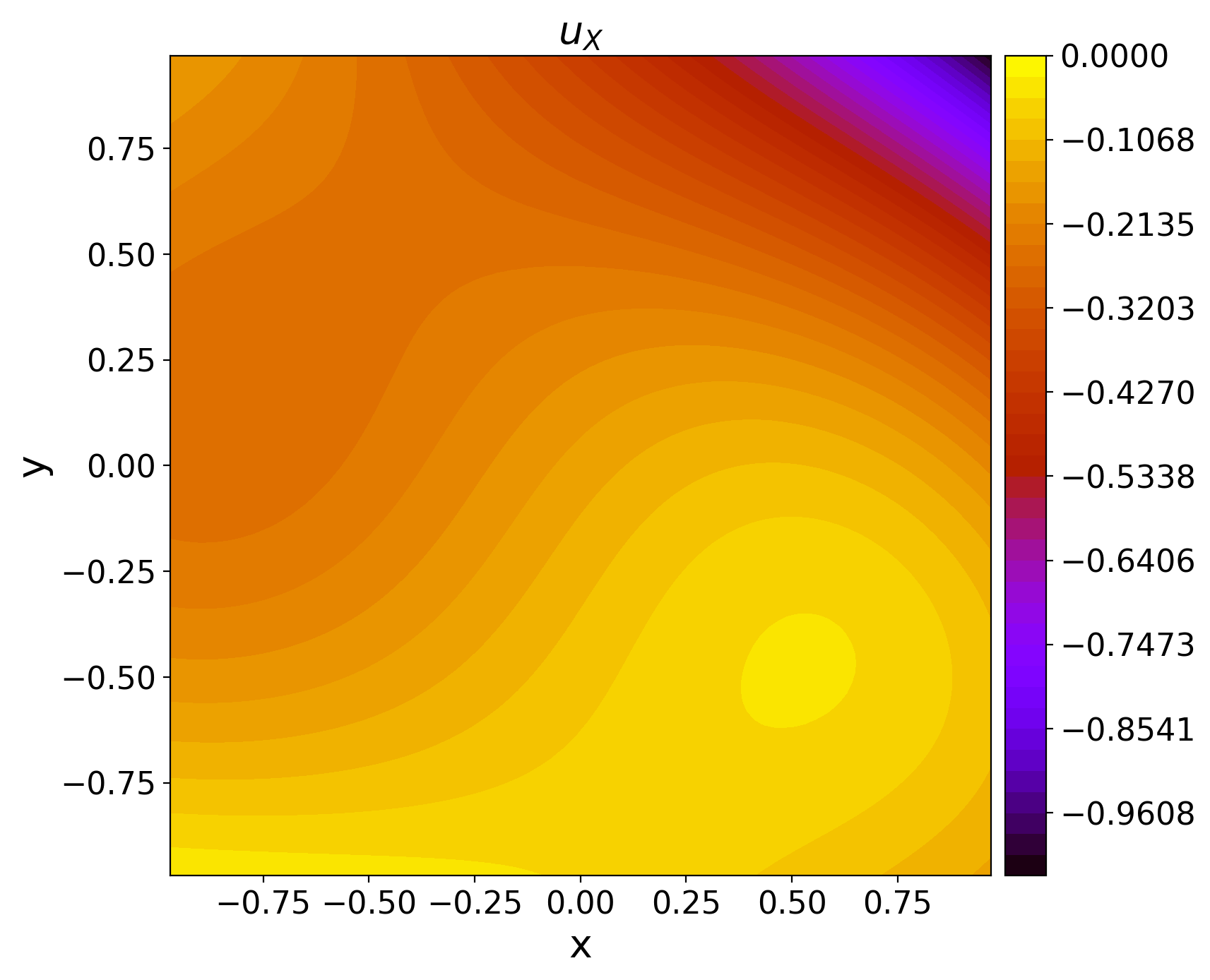}}
  \hspace*{0.2cm}
  \subfloat[][BA-GMLS-PD]{\includegraphics[height=0.21\textwidth,trim={3cm 0 0 0 },clip]{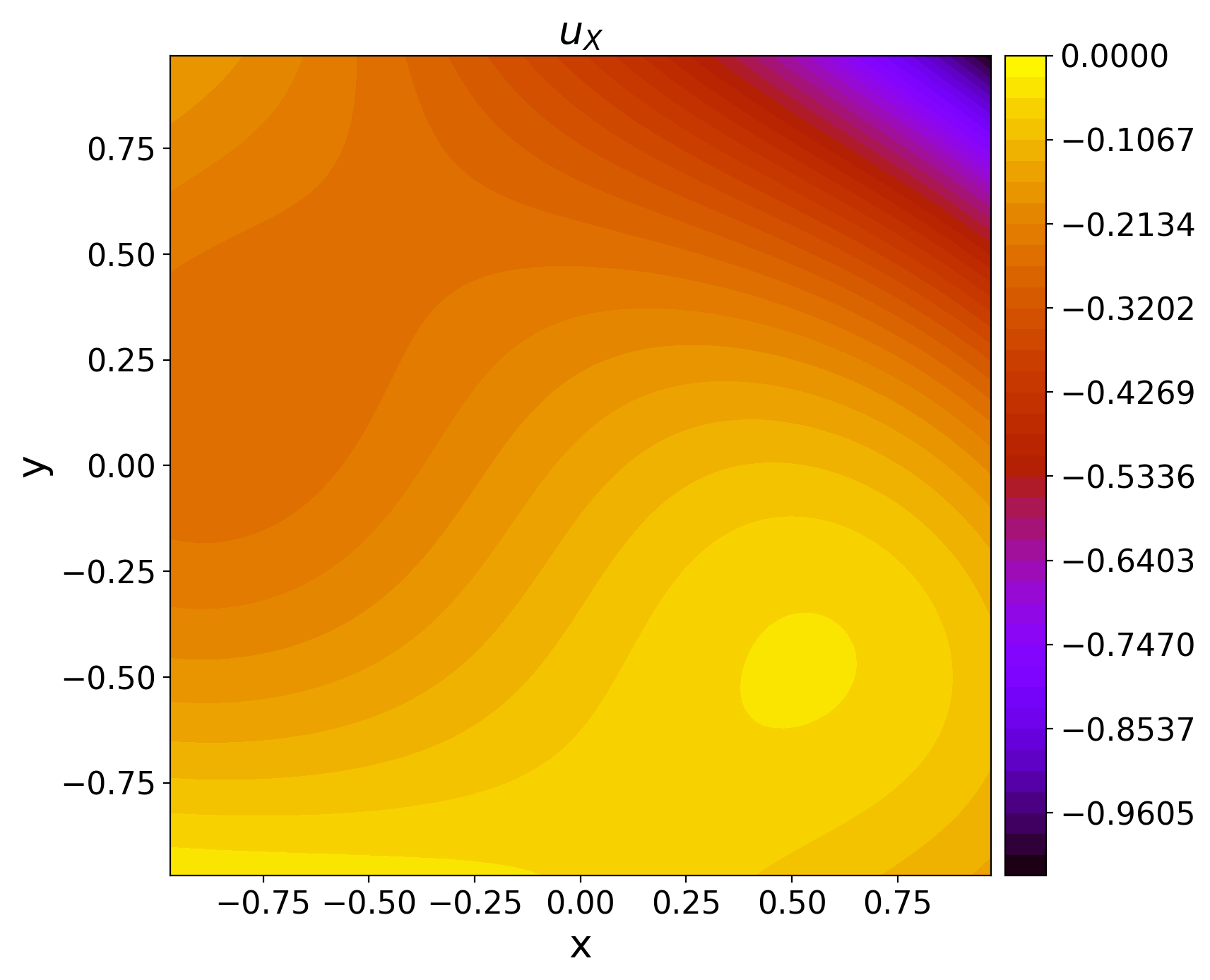}}
  \caption{Horizontal-displacement contours using quadratic formulations. Top and bottom rows correspond to the L1 and L3 non-uniform discretizations, respectively. While oscillations in the RK-PD and GMLS-PD fields are reduced with mesh refinement, they are not completely eliminated. On the other hand, there is no oscillation in the bond-associated solution fields, even on coarse grids.}
  \label{fig:square-contours}
\end{figure*}

\cref{fig:square-contours} demonstrates the instability issue in the underlying models. Contours of the horizontal displacements, solved by the quadratic formulations, are shown for two levels of non-uniform grids (denoted as L1 and L3 in \cref{fig:square-mesh}). While mesh refinement reduces the high degree of oscillations present on the coarse mesh for RK-PD and GMLS-PD, some spurious oscillations remain in the refined case. On the contrary, the bond-associated methods provide meaningful and stable solutions at all levels of mesh refinement.

\begin{figure*}[!ht]
  \centering
  \subfloat[][$\delta = 2.75 \, h$]{\includegraphics[height=0.29\textwidth]{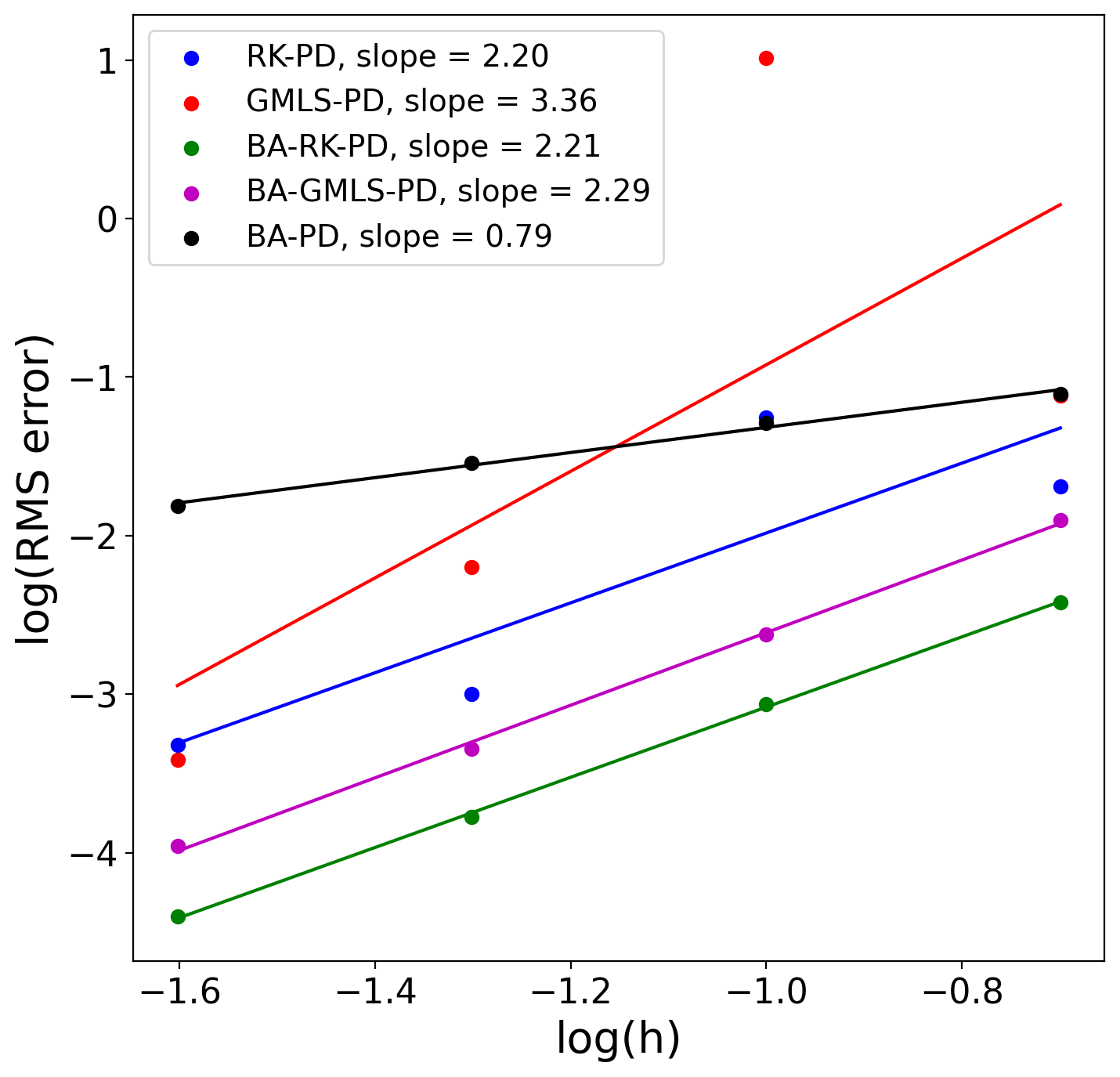}}
  \hspace*{0.3cm}
  \subfloat[][$\delta = 3.5 \, h$]{\includegraphics[height=0.29\textwidth]{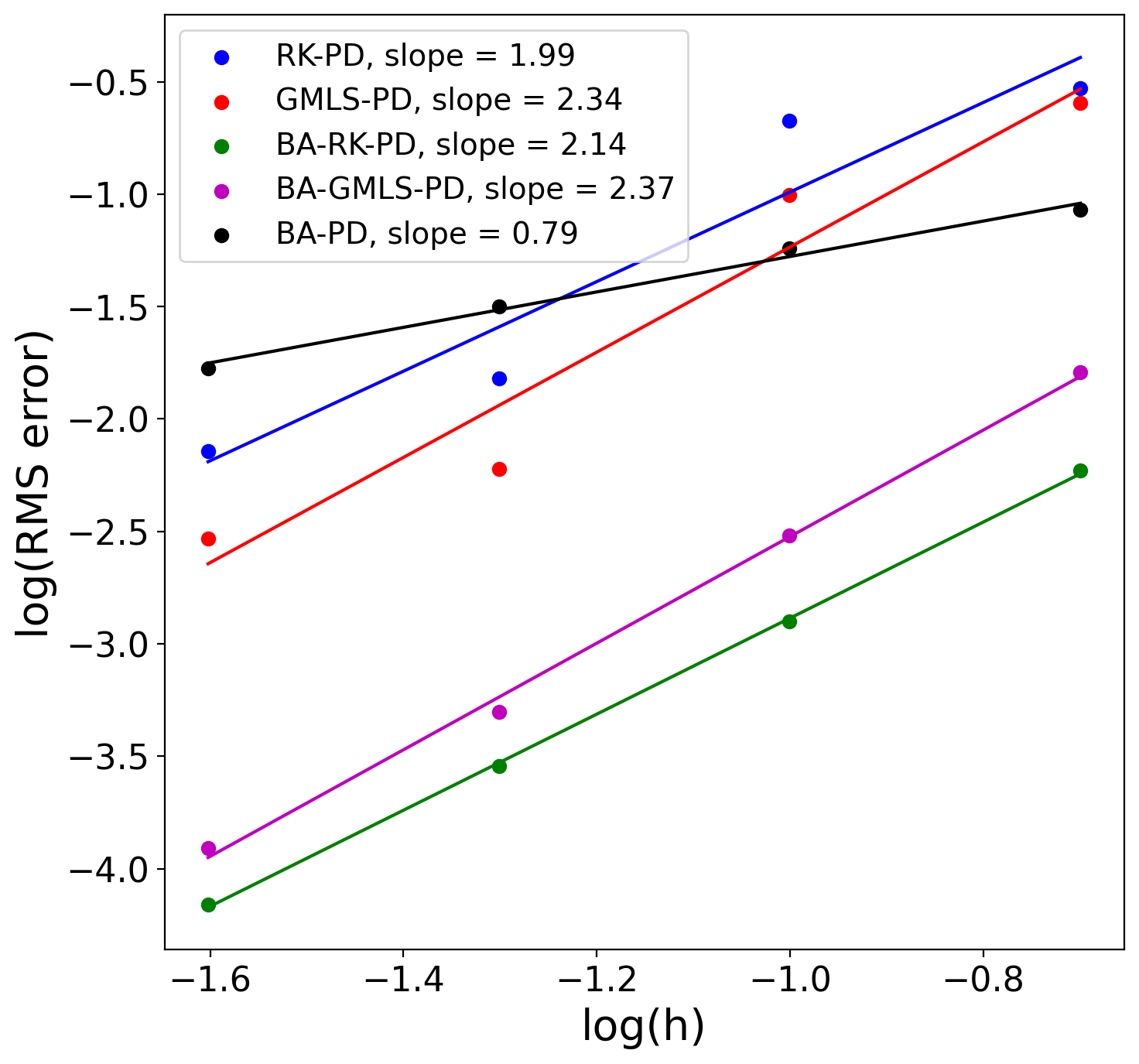}}
  \hspace*{0.3cm}
  \subfloat[][$\delta = 4.25 \, h$]{\includegraphics[height=0.29\textwidth]{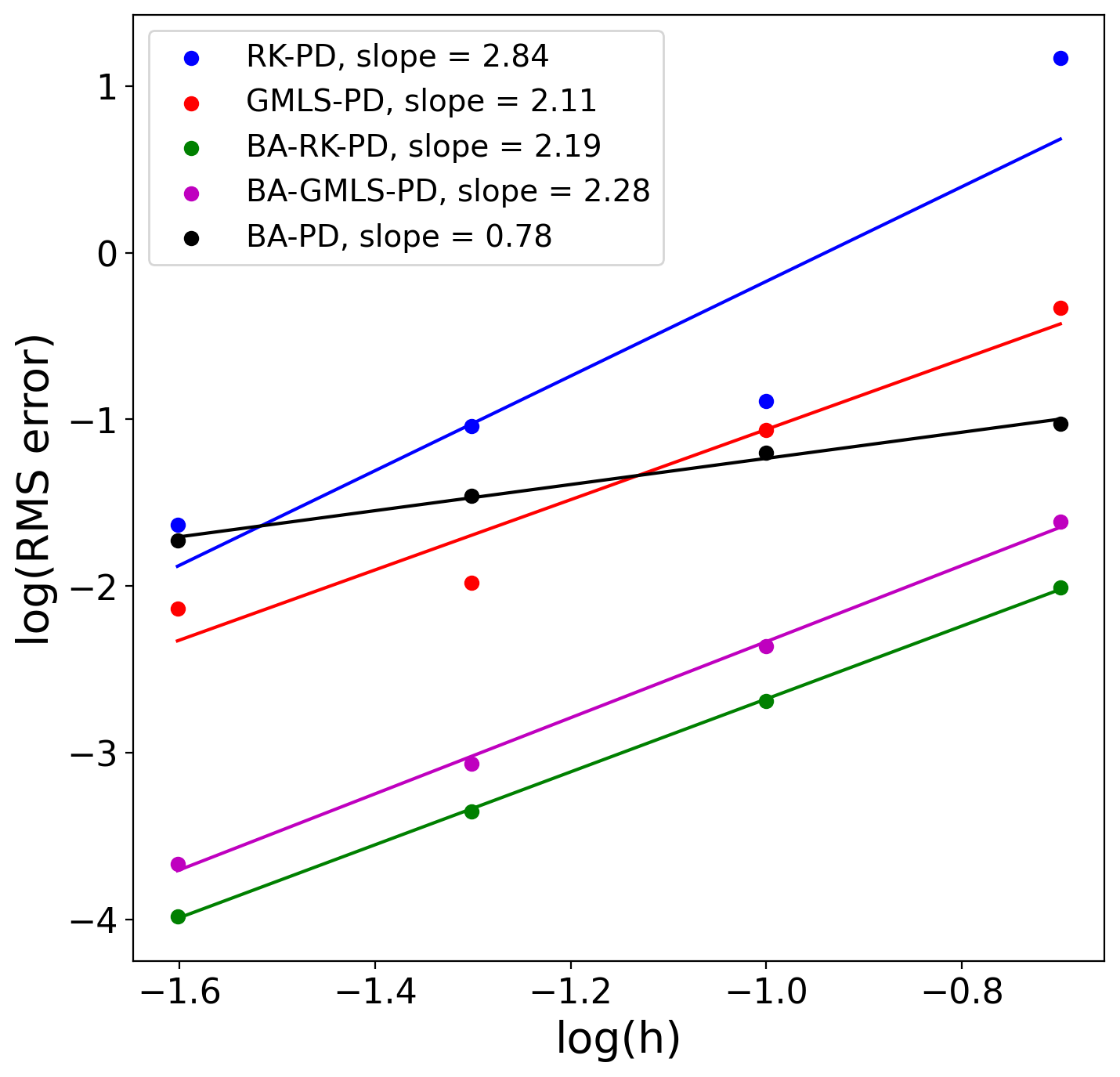}}
  \caption{Convergence of the quadratic models is compared for different horizon sizes. The horizon size has a large effect on RK-PD and GMLS-PD, while the bond-associated versions show robustness with respect to it. Plots show displacement RMS error values.}
  \label{fig:square-horizon}
\end{figure*}

We illustrate the effects of the neighborhood size in \cref{fig:square-horizon}, by considering three different horizons for the same nodal spacings. Quadratic models are used with non-uniform discretizations. In this test, for the non-stabilized models, i.e., RK-PD and GMLS-PD, error increases with $\delta$. On the other hand, the bond-associated methods demonstrate robustness with respect to the number of bonds in the neighborhood.

\subsubsection{Infinite plate with a circular hole in uniaxial tension}
\label{subsec:plate}

In this example, we consider a classical problem involving a curvilinear free surface, in which an infinite plate with a circular hole is subjected to far-field uni-axial tension. Under plane strain condition, the following solution can be obtained \citep[pp. 163--167]{gould1994introduction} using Airy stress functions \citep{michell1899direct}:
\begin{align*}
  & P_{rr}(r,\theta) = \frac{T_x}{2} \left[ 1 - \frac{a^2}{r^2} \right] + \frac{T_x}{2} \left[ 1 - 4 \frac{a^2}{r^2} + 3 \frac{a^4}{r^4} \right] \cos 2\theta , \notag \\
  & P_{\theta\theta}(r,\theta) = \frac{T_x}{2} \left[ 1 + \frac{a^2}{r^2} \right] - \frac{T_x}{2} \left[ 1 + 3 \frac{a^4}{r^4} \right] \cos 2\theta , \notag \\
  & P_{r\theta}(r,\theta) = P_{\theta r}(r,\theta) = - \frac{T_x}{2} \left[ 1 + 2 \frac{a^2}{r^2} - 3 \frac{a^4}{r^4} \right] \sin 2\theta , 
\end{align*}
where $a$ is the hole radius and $T_x$ is the magnitude of the applied far-field stress.

\begin{figure}[!htbp]
  \centering
  \includegraphics[width=0.45\textwidth]{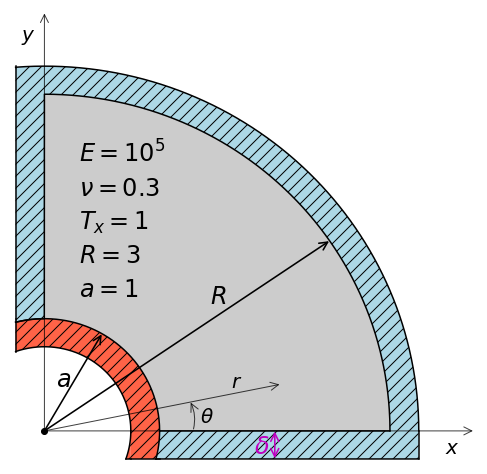}
  \caption{Infinite plate with a circular hole. Blue color indicates the natural boundary nodes, prescribed with the exact stress values. Red color denotes the free-surface region, where fictitious nodes are placed and given zero stress values. Nodes on the horizontal axis have fixed vertical positions and vice versa.}
  \label{fig:geometry}
\end{figure}

The infinite problem is modeled as a finite quarter plate with the exact stress values prescribed on the non-local collar surrounding the perimeter of the plate. To model the free-surface behavior at the hole region, a non-local, zero-stress, fictitious layer is utilized in that region. That is, as noted previously, the free-surface nodes are given zero stress values that directly contribute to the determination of $\nabla_h \cdot \mbt{P}$ for their material neighbors. Similar to the way the natural-bc nodes are employed, the free-surface nodes are not involved in the evaluation of the kinematic variable for their neighbors. Symmetry of the problem is used to apply kinematic (Dirichlet) boundary conditions by assigning zero vertical displacements to nodes on the horizontal axis and vice versa. Schematic of this modeling approach is provided in \cref{fig:geometry}.

\begin{figure*}[!ht]
  \centering
  \subfloat[][Polar - L0]{\includegraphics[height=0.25\textwidth]{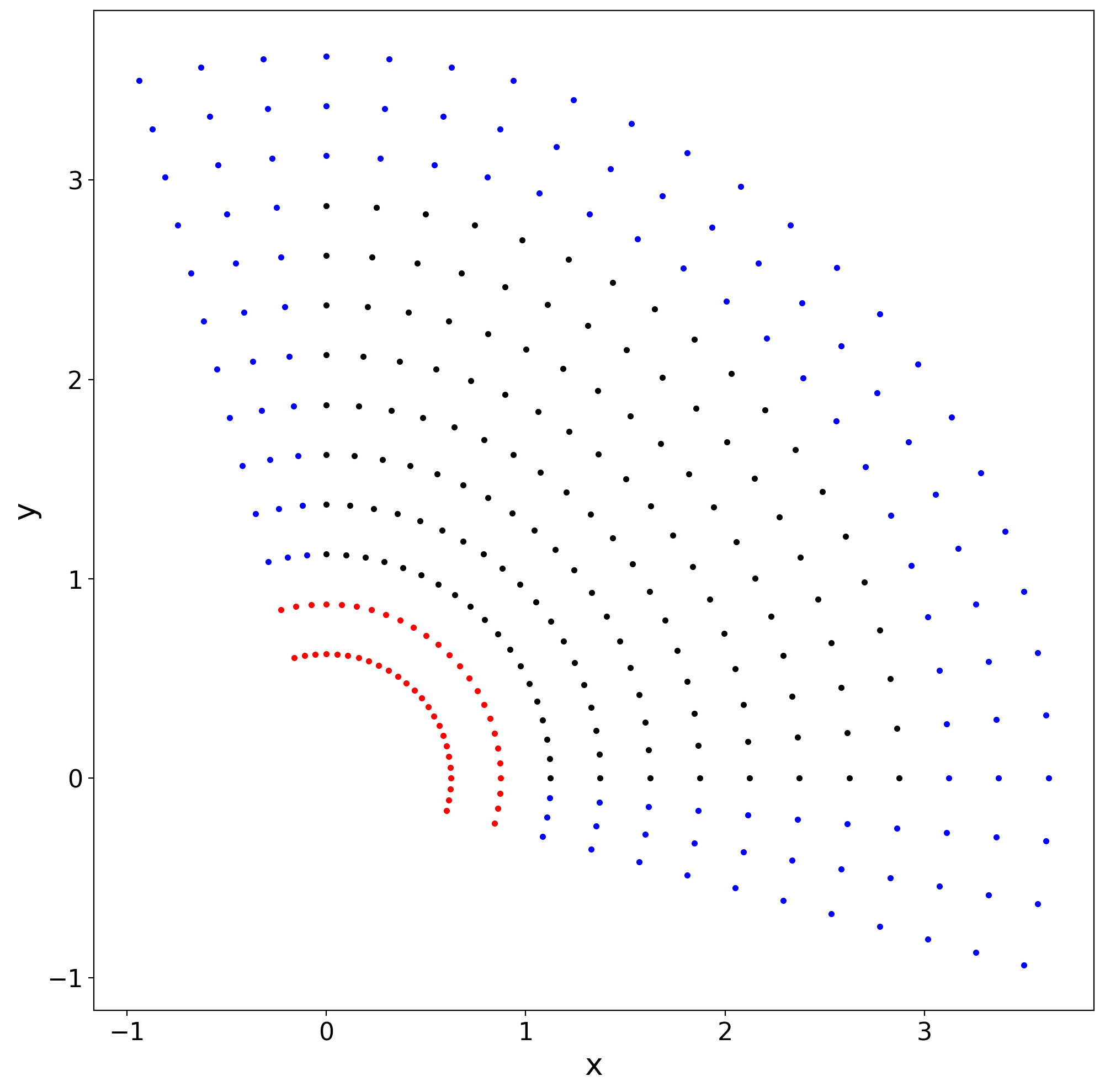}}
  \hspace*{0.1cm}
  \subfloat[][Polar - L1]{\includegraphics[height=0.25\textwidth,trim={2.5cm 0 0 0 },clip]{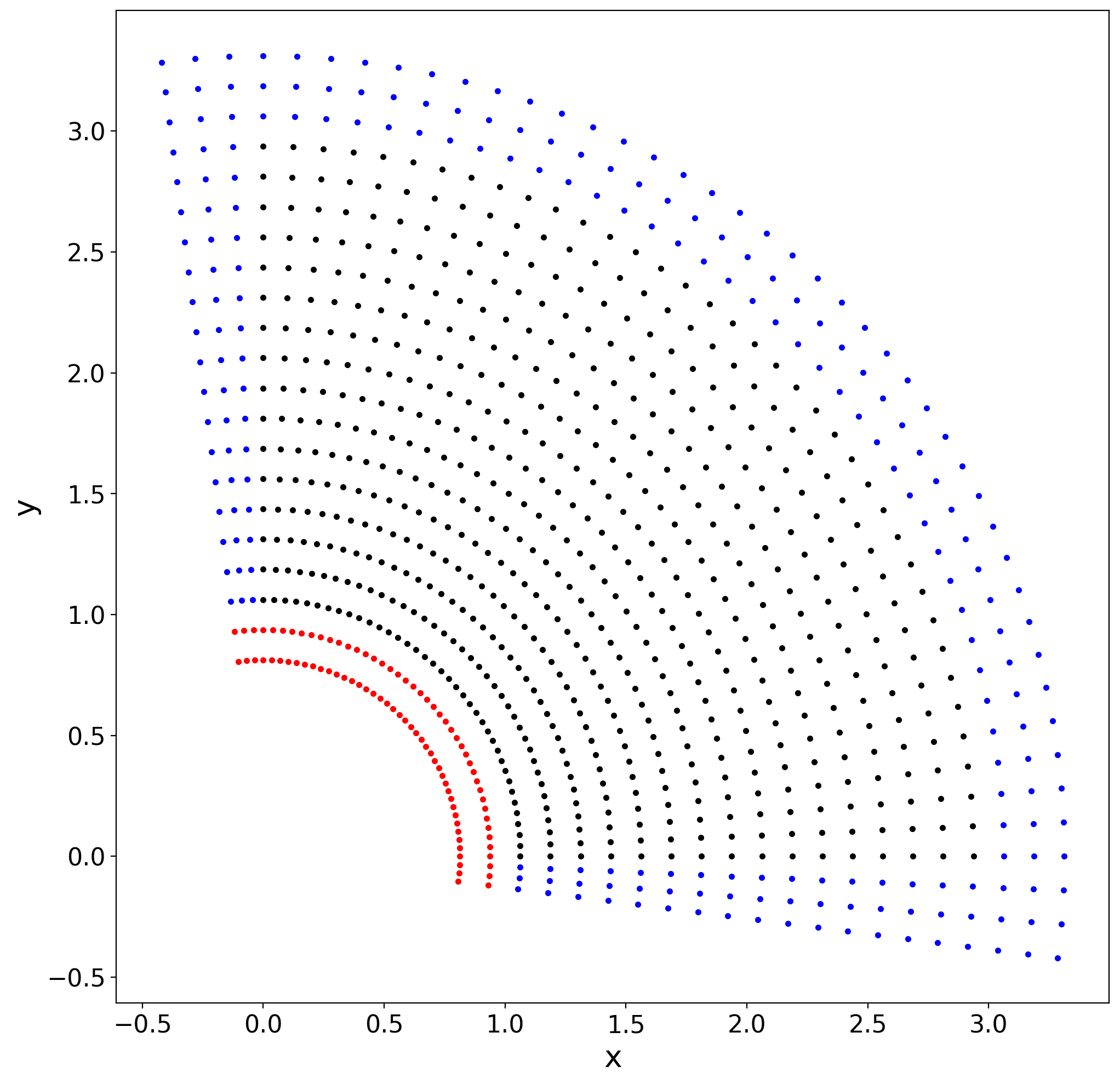}}
  \hspace*{0.1cm}
  \subfloat[][Polar - L2]{\includegraphics[height=0.25\textwidth,trim={2cm 0 0 0 },clip]{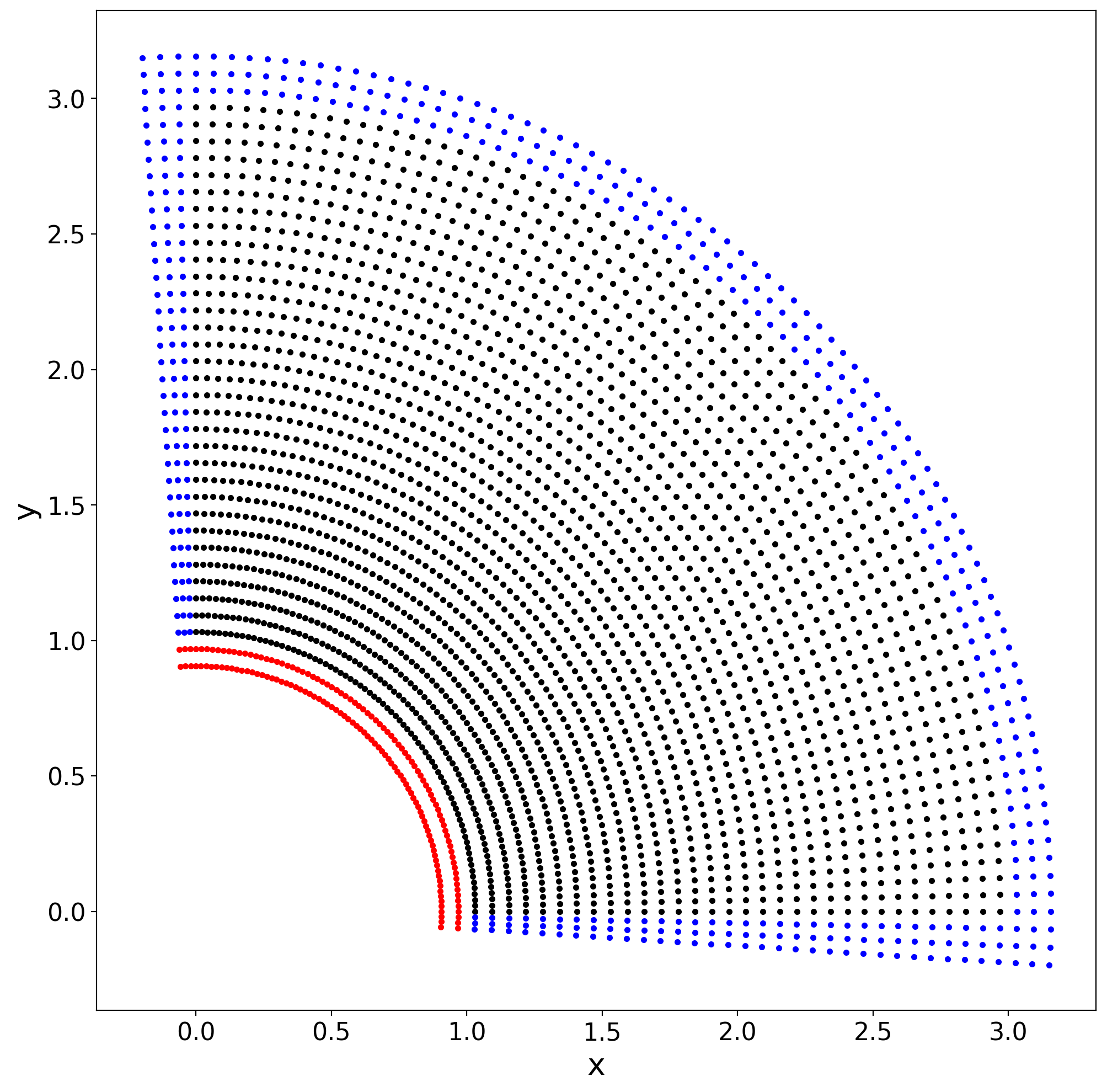}}
  \hspace*{0.1cm}
  \subfloat[][Polar - L3]{\includegraphics[height=0.25\textwidth,trim={2cm 0 0 0 },clip]{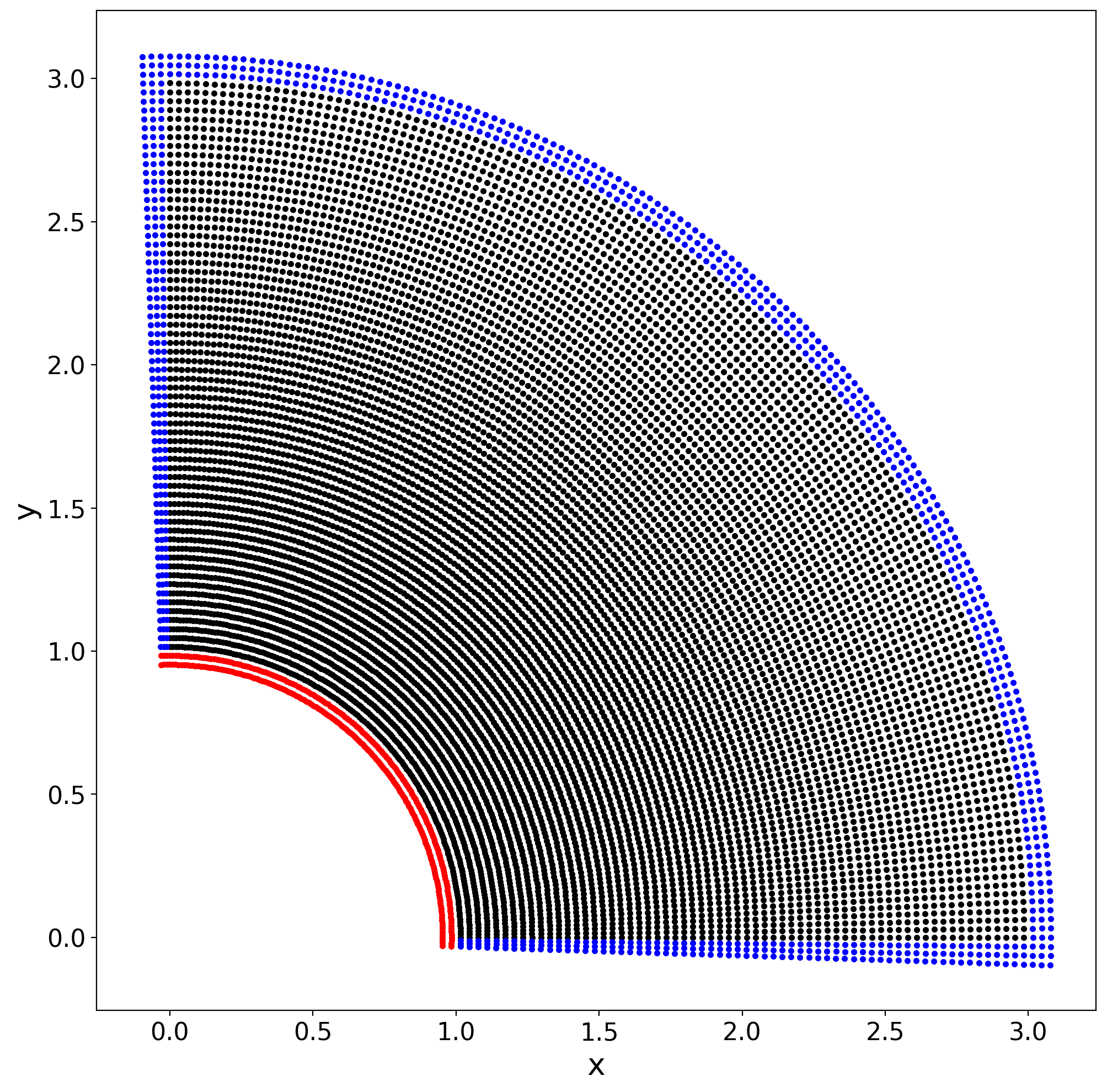}}

  \subfloat[][Triangular - L0]{\includegraphics[height=0.245\textwidth]{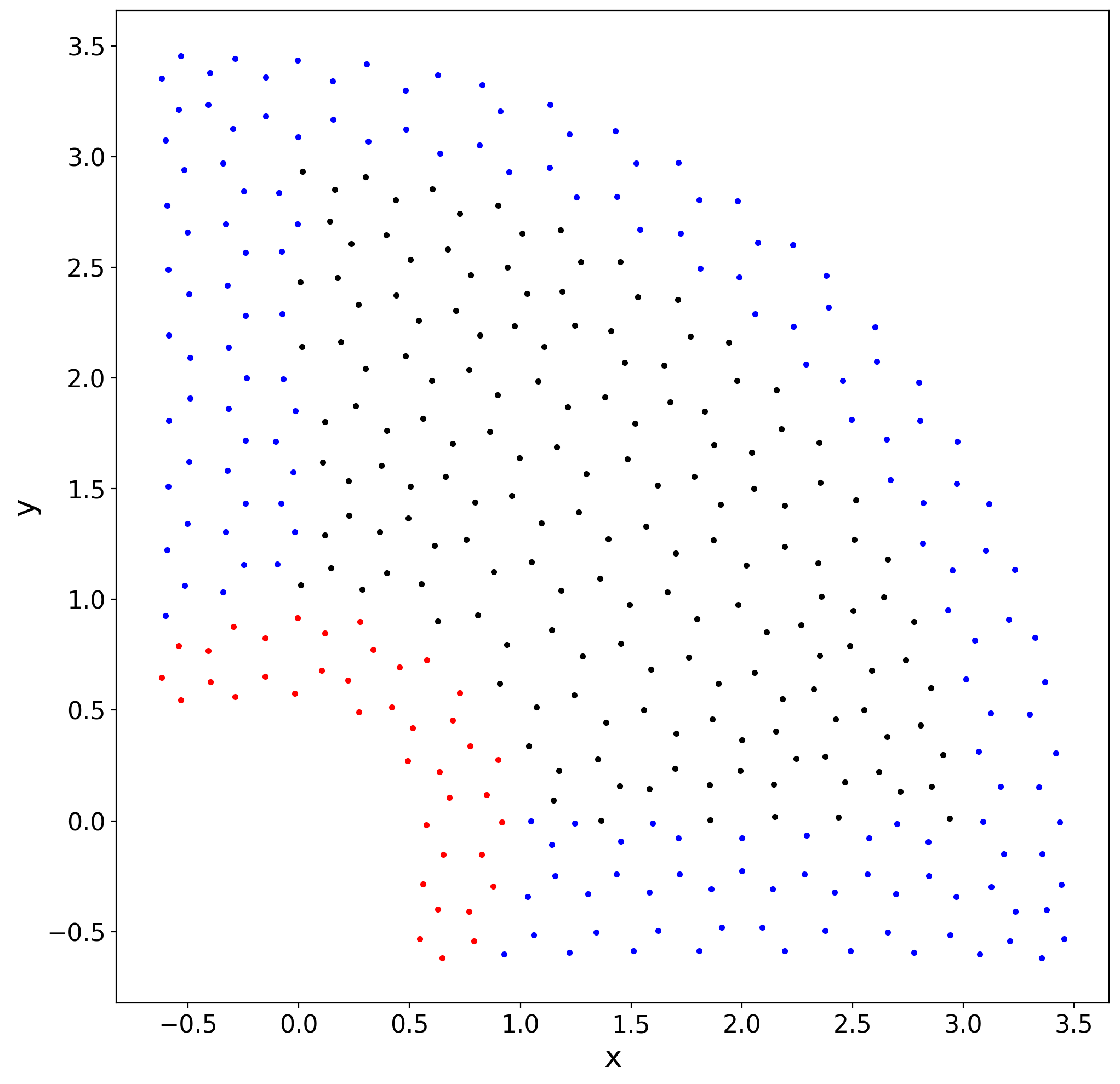}}
  \hspace*{0.1cm}
  \subfloat[][Triangular - L1]{\includegraphics[height=0.245\textwidth,trim={2cm 0 0 0 },clip]{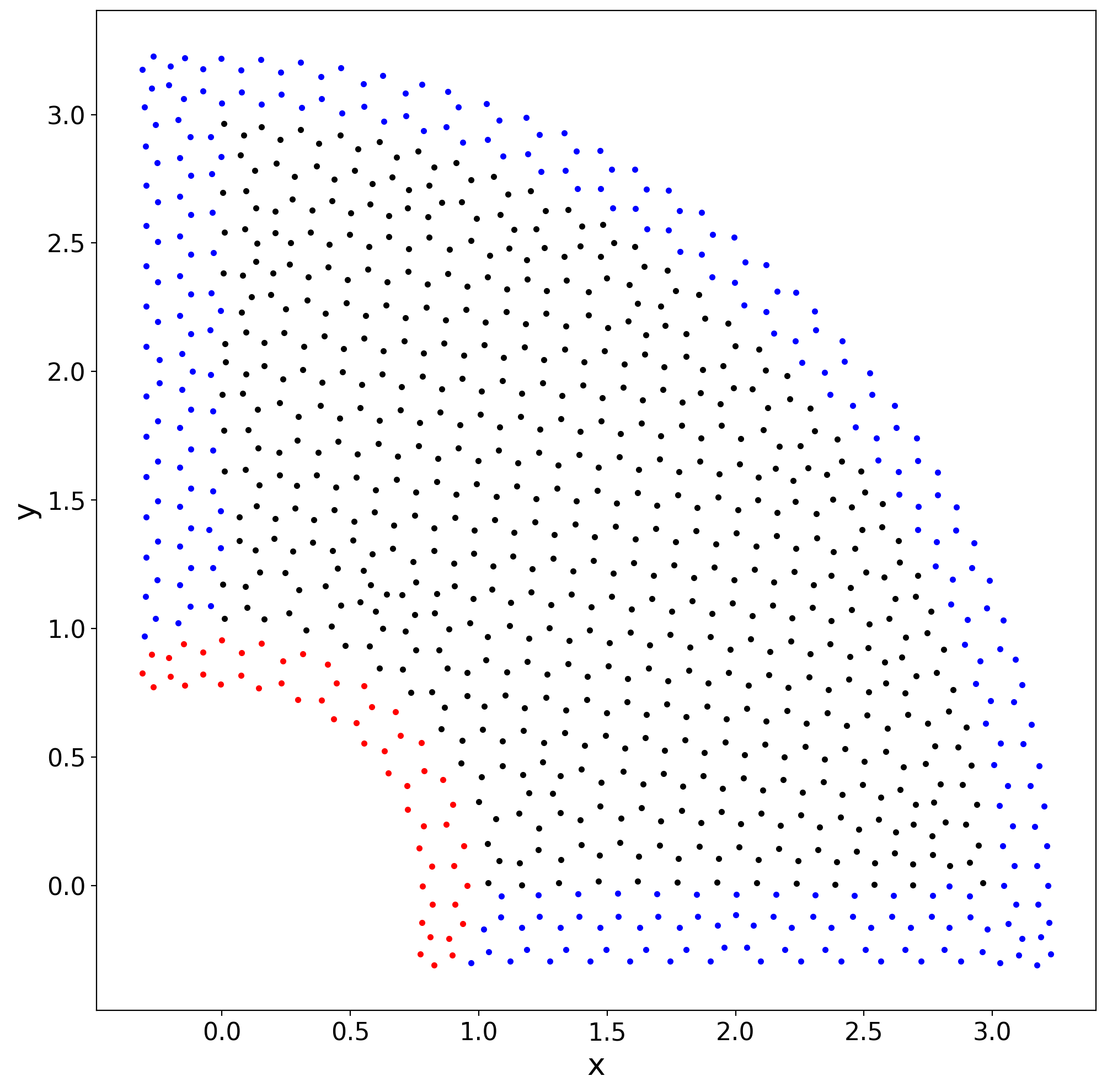}}
  \hspace*{0.1cm}
  \subfloat[][Triangular - L2]{\includegraphics[height=0.245\textwidth,trim={2cm 0 0 0 },clip]{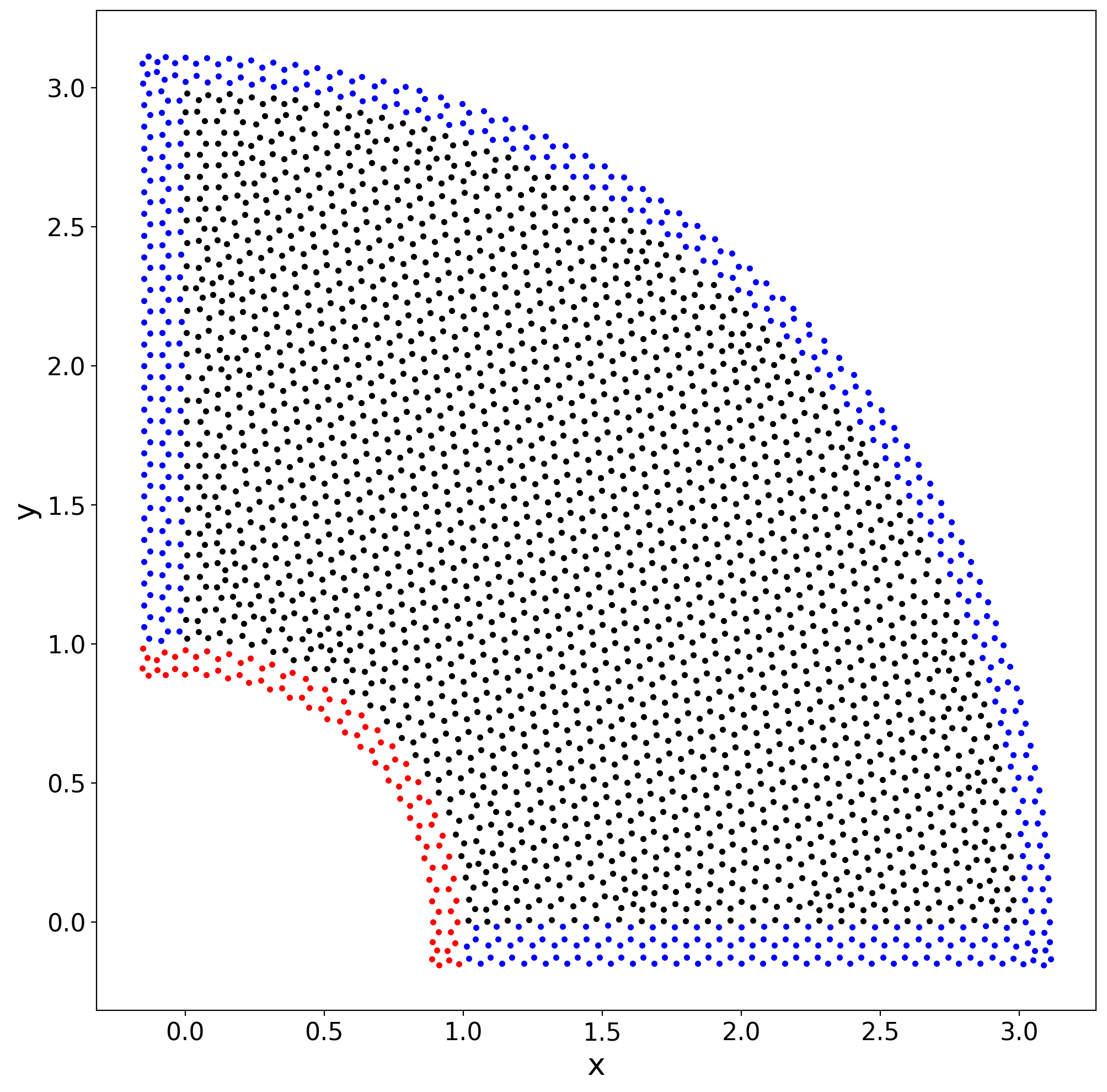}}
  \hspace*{0.1cm}
  \subfloat[][Triangular - L3]{\includegraphics[height=0.245\textwidth,trim={2cm 0 0 0 },clip]{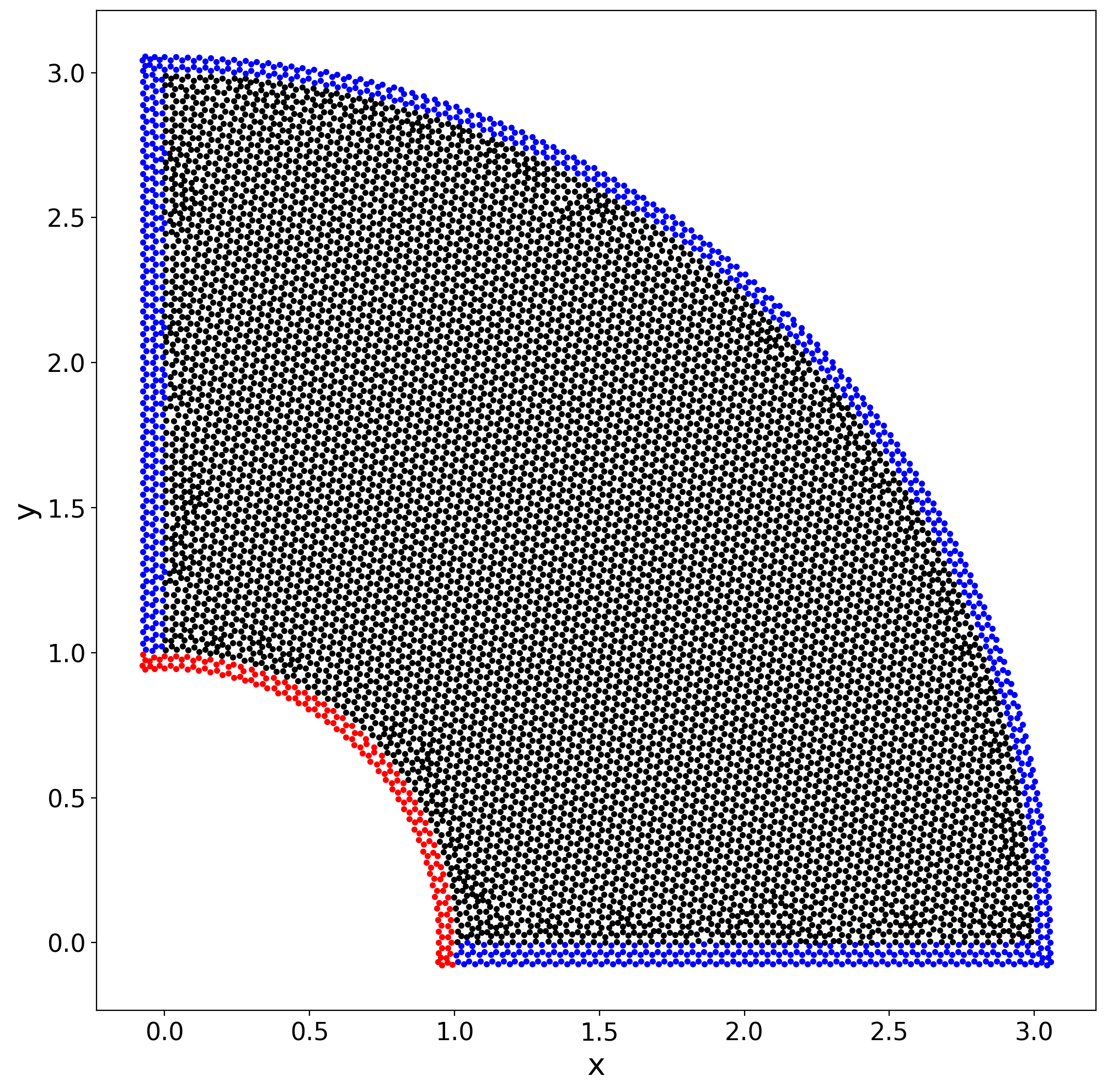}}
  \caption{Meshes for the plate with a circular hole problem. The meaning of red and blue colors is the same as in the previous figure.}
  \label{fig:plate-mesh}
\end{figure*}

In this example, for obtaining a nodal discretization, a finite element mesh is used, where each element is replaced by a meshfree node. The meshfree node is placed on the centroid of the element and is given its area (volume in 3-D). Two strategies are utilized to generate the mesh for this domain: (1) a structured meshing algorithm, where finite element nodes are uniformly distributed along a polar coordinate system to obtain quadrilateral elements, and (2) a semi-unstructured algorithm using the open-source Python library {\em pyGmsh}, which is based on Gmsh \citep{geuzaine2007gmsh}, where Delaunay triangular meshes are constructed. We denote these two schemes as {\em polar} and {\em triangular} discretizations, respectively, in the remainder of this section. Different mesh refinement levels are used to perform a convergence test. For each meshfree discretization, average node spacing is defined by $$h = \sqrt{\frac{\sum_i^N A_i}{N}},$$ where $A_i$ is element area. The approximate average spacing of $h\approx[0.2, \, 0.1, \, 0.05, \, 0.025]$ is considered for both mesh types. The different levels of nodal discretizations are shown in \cref{fig:plate-mesh}.

\begin{figure*}[!htbp]
  \centering
  \includegraphics[width=0.65\textwidth]{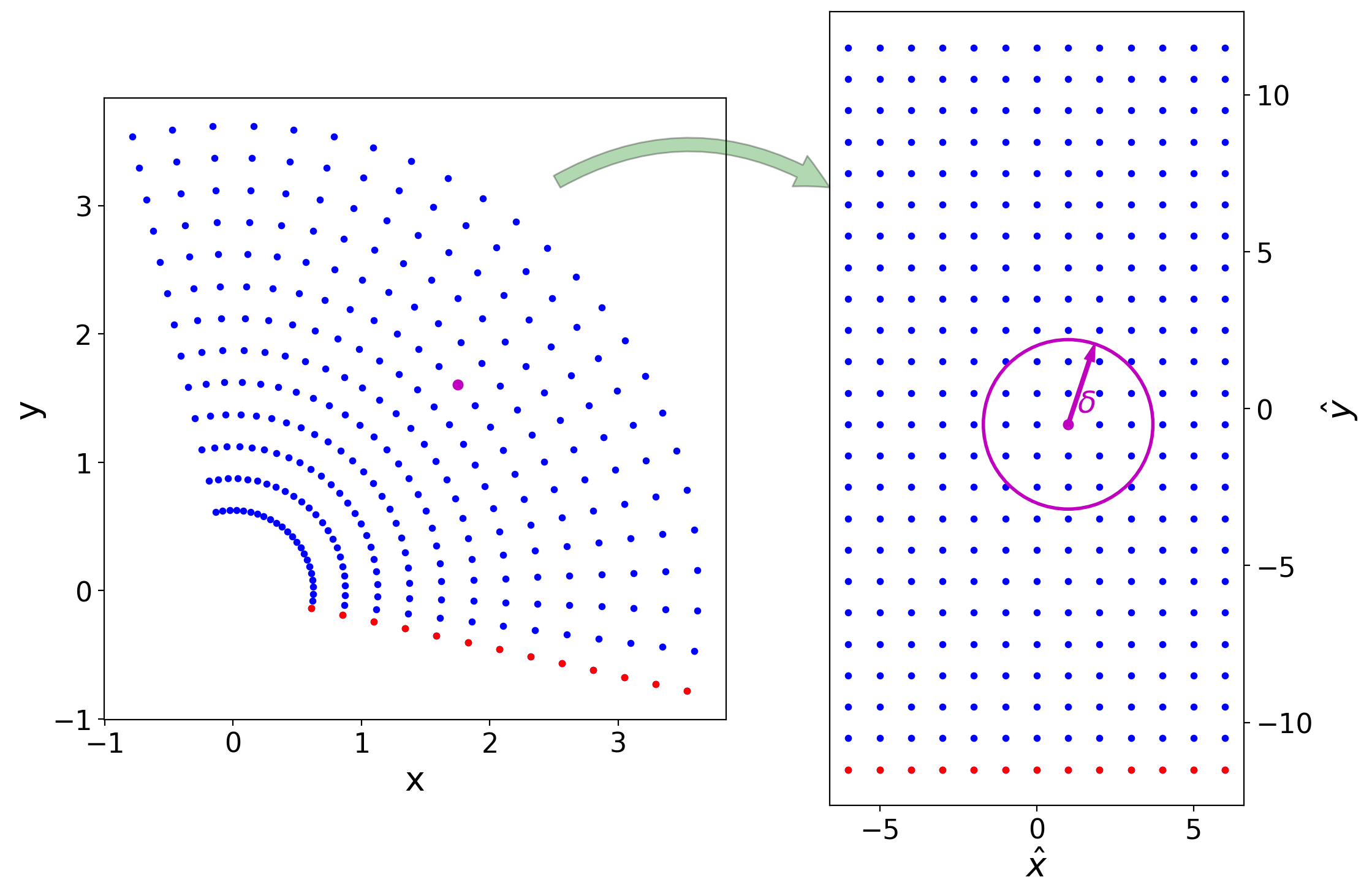}
  \caption{Mapping is used to define the neighbor sets in the parametric space, resulting in a uniform number of neighbors for each node (except near the boundaries).}
  \label{fig:mapping}
\end{figure*}

As explained in Part I of this paper, using a constant horizon size in the {\em physical} domain of the polar case to construct neighborhoods would result in a significantly larger neighbor sets near the hole, where more nodes are in place. That would reduce efficiency and potentially damage the robustness of the RK and GMLS algorithms as they depend on the number of neighbors (cf. \cref{fig:square-horizon}). To obtain neighbor sets with consistent number of neighbors regardless of the location of the neighborhood (excluding in the vicinity of boundaries), while maintaining consistency between neighbors (i.e., if a bond $IJ$ exists, $JI$ also exists), a mapping algorithm is used. That is, the nodes are transformed from the physical domain to a {\em parametric} space, where a uniform spacing of 1 is achieved between the adjacent nodes. The neighbor sets are constructed in the parametric space. This approach is illustrated in \cref{fig:mapping}. In the parametric configuration, the horizon size for the linear, quadratic, and cubic models is adopted as $\delta = 1.75$, $\delta = 2.75$, and $\delta = 3.75$, respectively. Note that although the influence functions are computed in the parametric space, the quadrature weights should be computed in the physical space, to obtain consistent gradient operators. 

The neighbor sets, for the triangular discretization scheme, are formed in the physical space. In this case, the horizon size is chosen as $\delta = 2.25 h$, $\delta = 3.25 h$, and $\delta = 4.25 h$ for the linear, quadratic, and cubic methods, respectively.

\begin{figure*}[!ht]
  \centering
  \subfloat[][Polar discretization - Linear]{\includegraphics[height=0.39\textwidth]{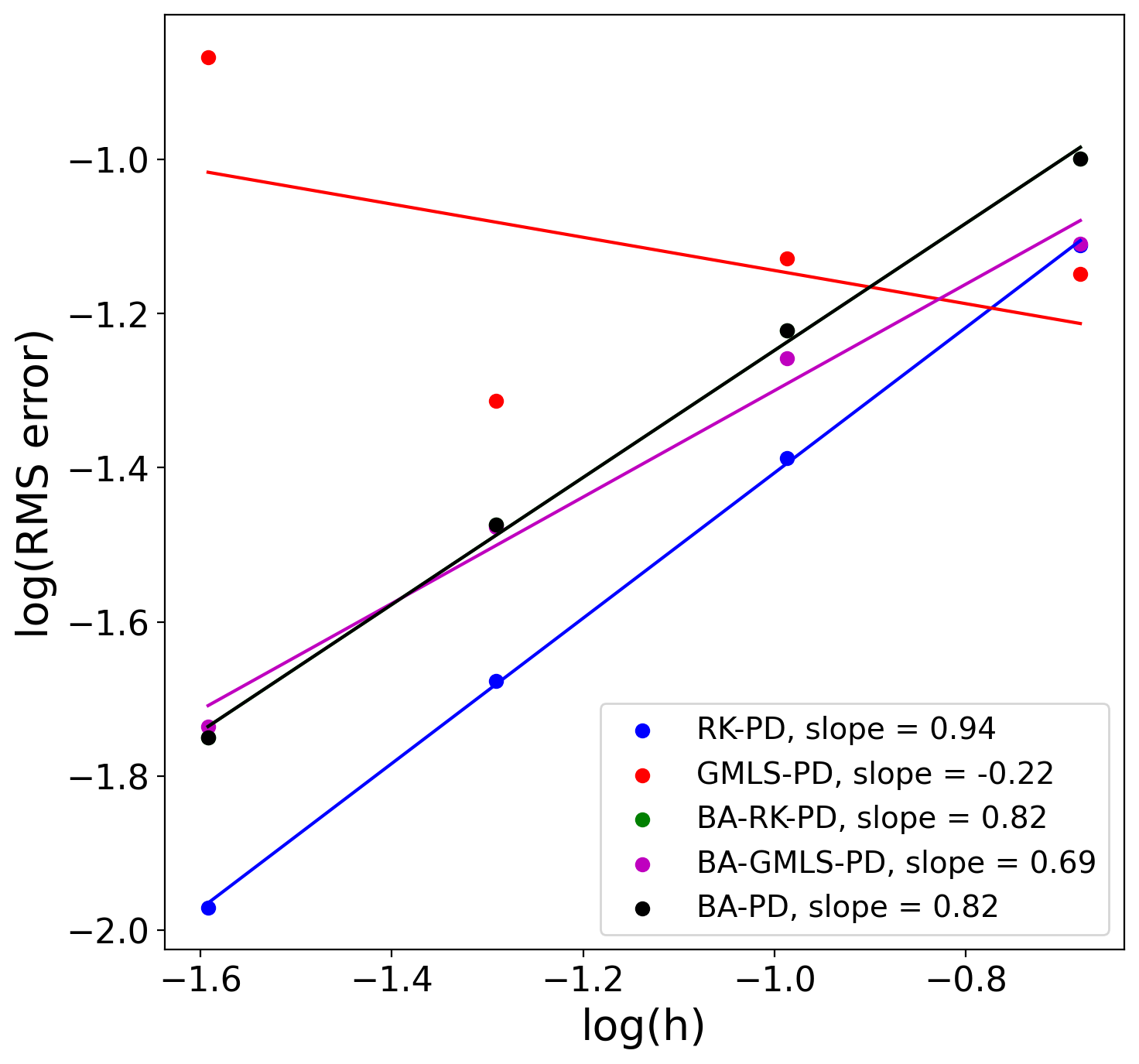}}
  \hspace*{1cm}
  \subfloat[][Triangular discretization - Linear]{\includegraphics[height=0.39\textwidth]{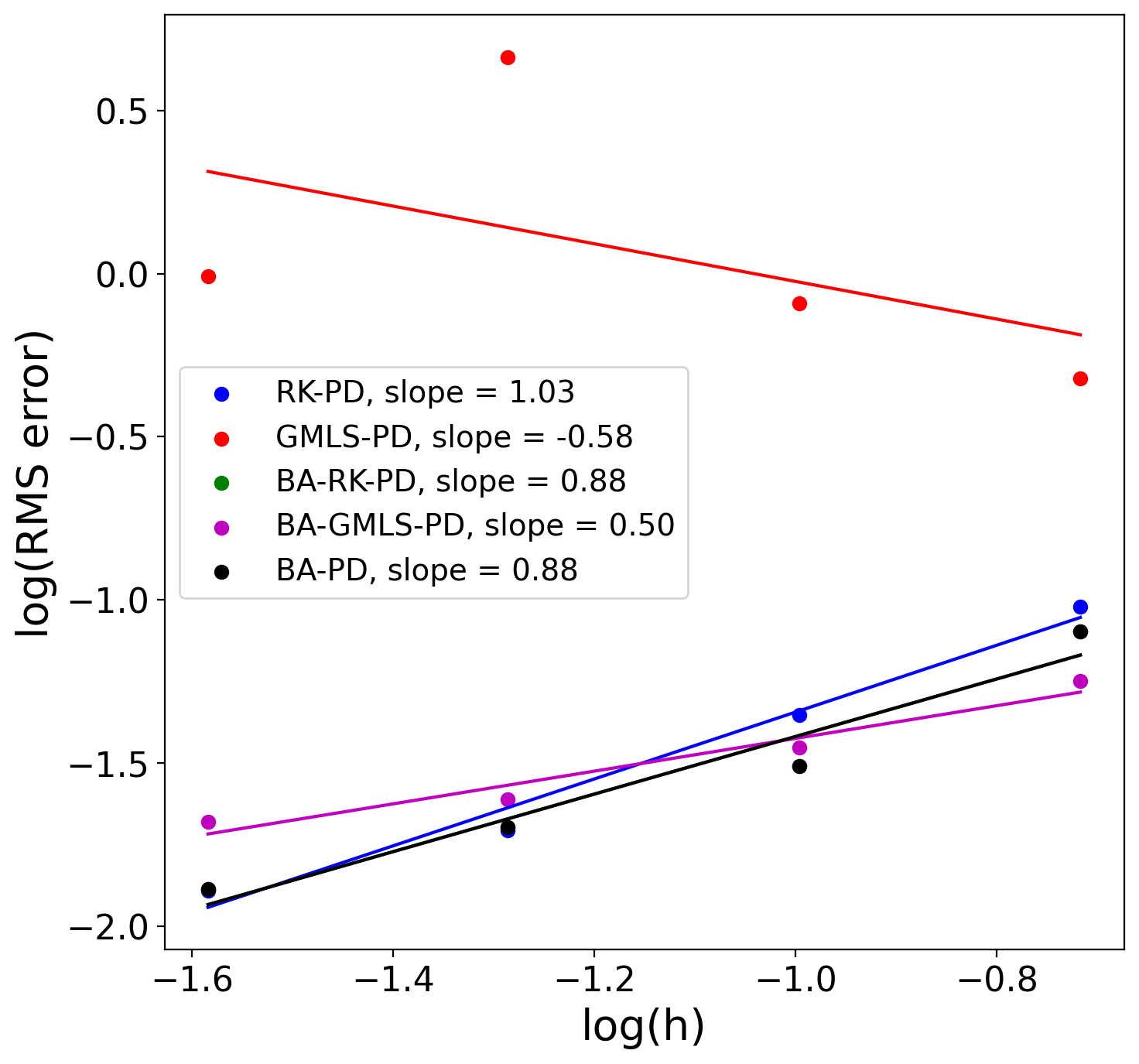}}

  \subfloat[][Polar discretization - Quadratic]{\includegraphics[height=0.39\textwidth]{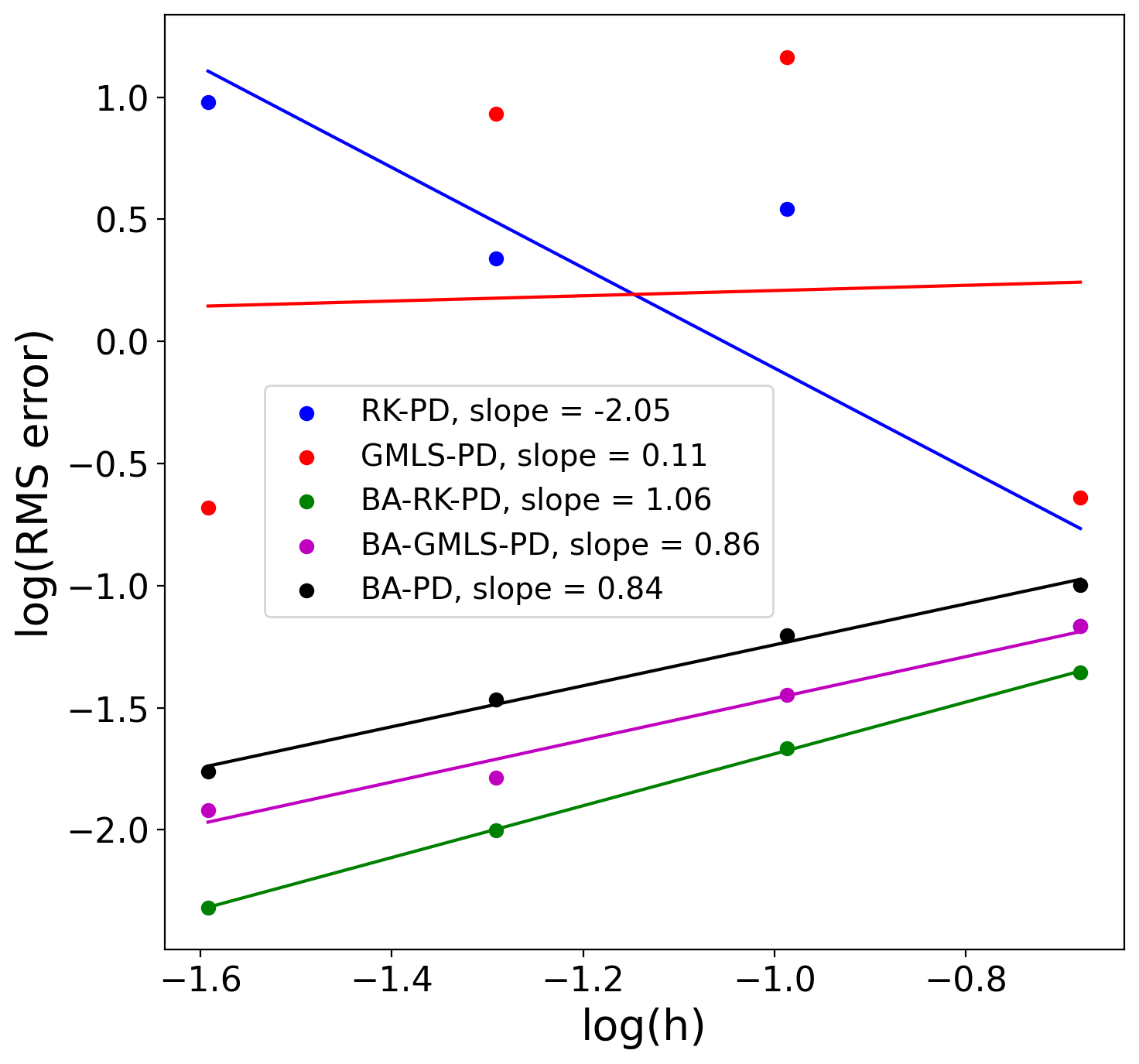}}
  \hspace*{1cm}
  \subfloat[][Triangular discretization - Quadratic]{\includegraphics[height=0.39\textwidth]{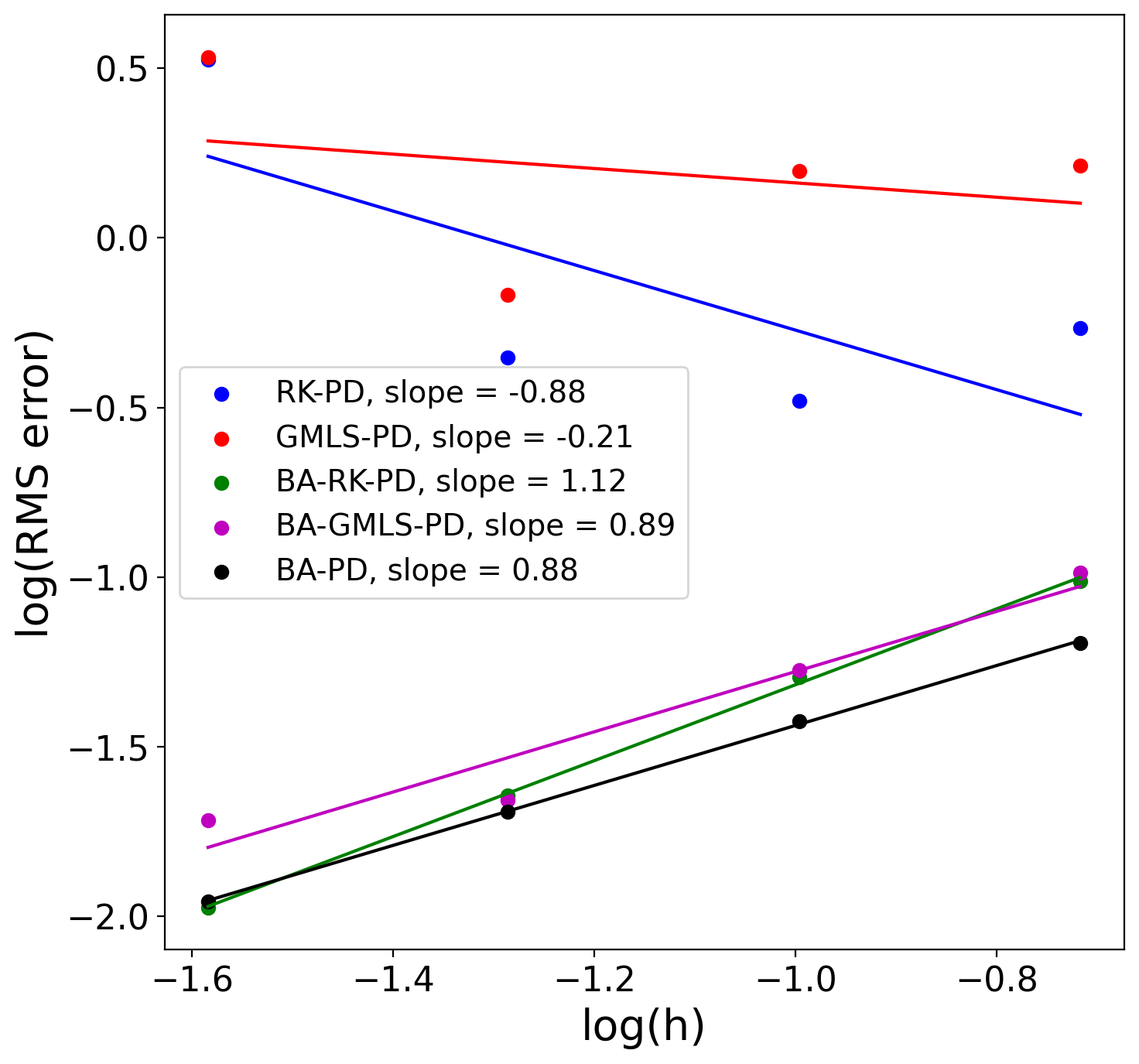}}

  \subfloat[][Polar discretization - Cubic]{\includegraphics[height=0.39\textwidth]{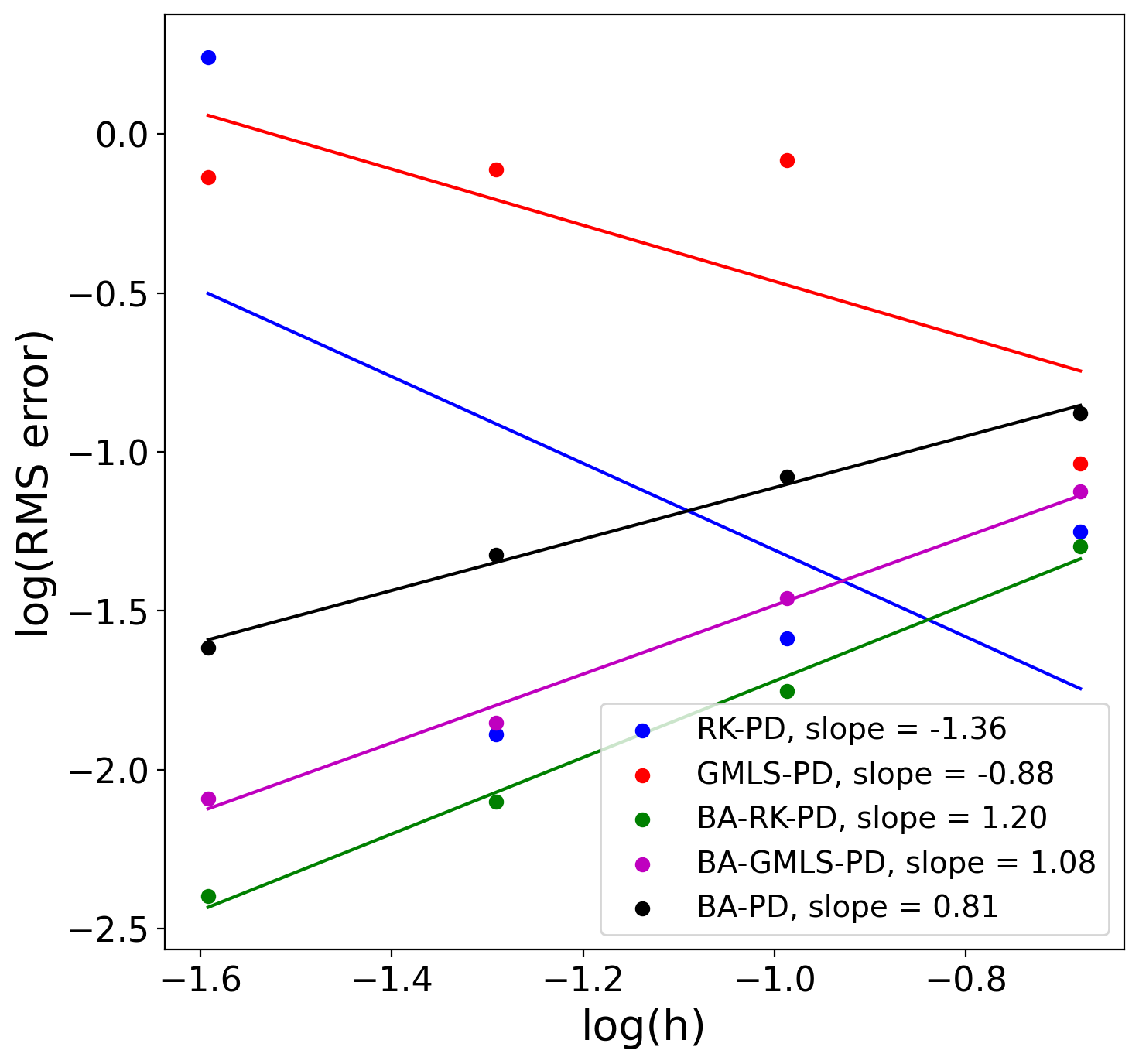}}
  \hspace*{1cm}
  \subfloat[][Triangular discretization - Cubic]{\includegraphics[height=0.39\textwidth]{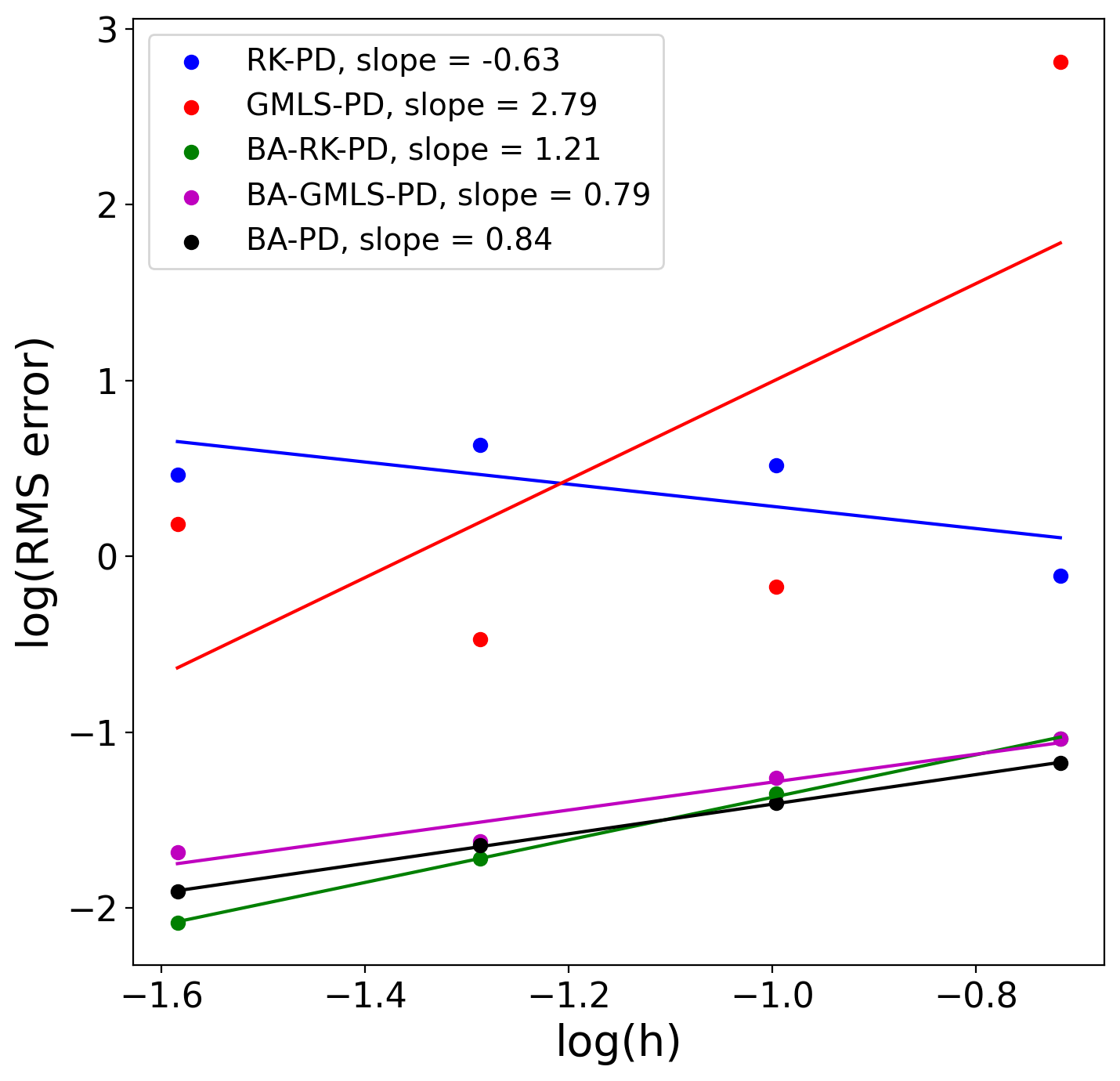}}
  \caption{Convergence of RK-PD, GMLS-PD, BA-RK-PD, BA-GMLS-PD, and BA-PD discretizations for the plate with a circular hole problem. Linear, quadratic, and cubic formulations are tested on the polar and triangular discretizations. BA-PD and BA-RK-PD overlap for (a--b). Plots show stress RMS errors.}
  \label{fig:plate-convergence}
\end{figure*}

The stress RMS errors are shown in \cref{fig:plate-convergence}. The bond-associated models obtained a near-linear convergence rate for quadratic and cubic cases (and a somewhat smaller rate for the linear case), which demonstrates asymptotically compatible convergence to the local solution in this problem. Note that the solution here is not smooth, involving a jump in stresses near the free surface, and only a first-order convergence rate with the non-local approach can be achieved. RK-PD linear case is asymptotically compatible; however, as soon as the number of neighbors increases (i.e., larger horizon), neither of the quadratic or cubic methods are able to achieve good results. The GMLS-PD results are associated with high degrees of error in all cases here. The difficulty of the base models in converging to the solution is attributed to the instability problem, as noted earlier. For the polar discretizations, the higher-order correction enhances the accuracy of the BA variants in this example.

\begin{figure*}[!ht]
  \centering
  \subfloat[][RK-PD]{\includegraphics[height=0.4\textwidth]{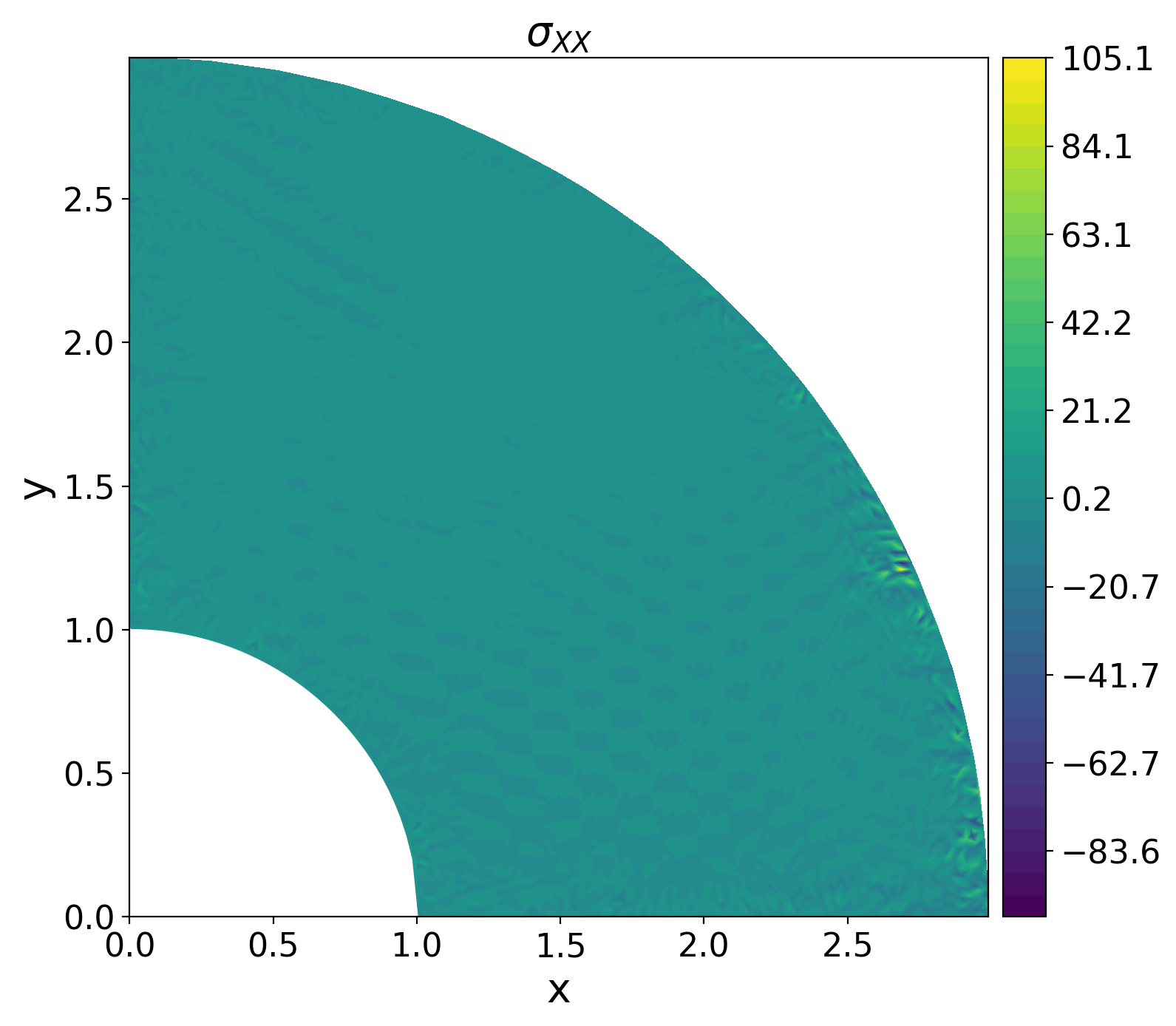}}
  \hspace*{0.5cm}
  \subfloat[][GMLS-PD]{\includegraphics[height=0.4\textwidth]{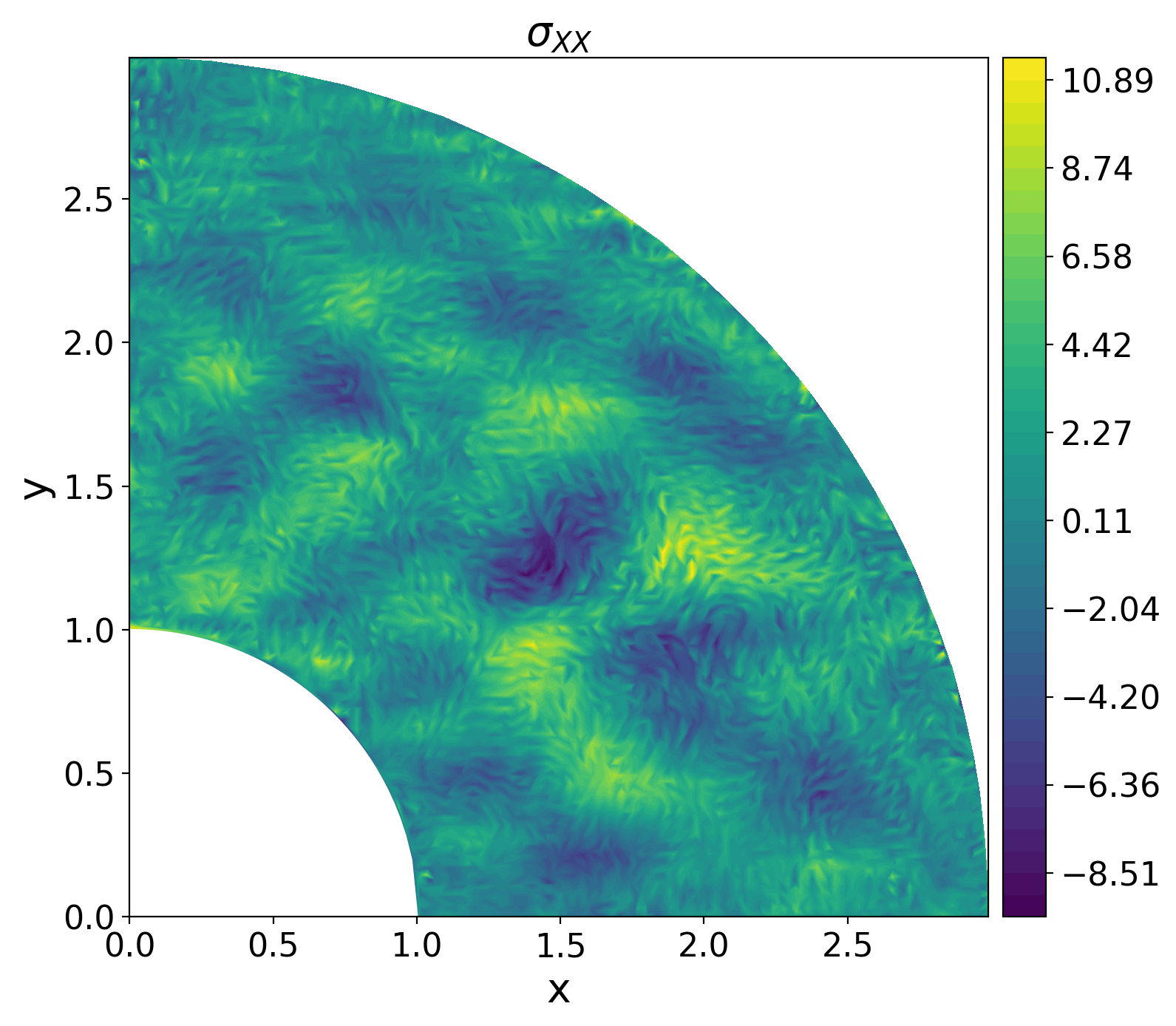}}

  \subfloat[][BA-RK-PD]{\includegraphics[height=0.4\textwidth]{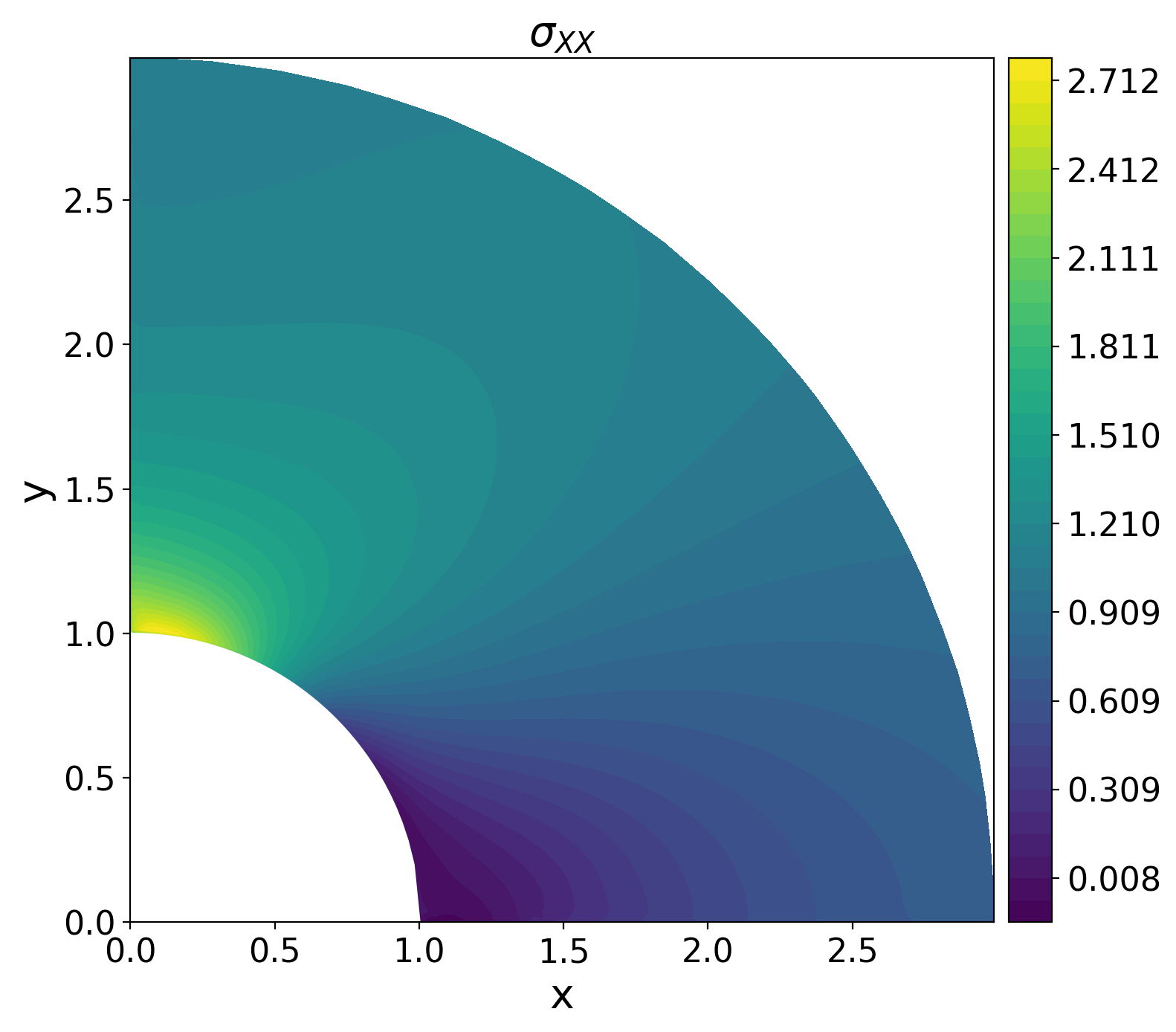}}
  \hspace*{0.5cm}
  \subfloat[][BA-GMLS-PD]{\includegraphics[height=0.4\textwidth]{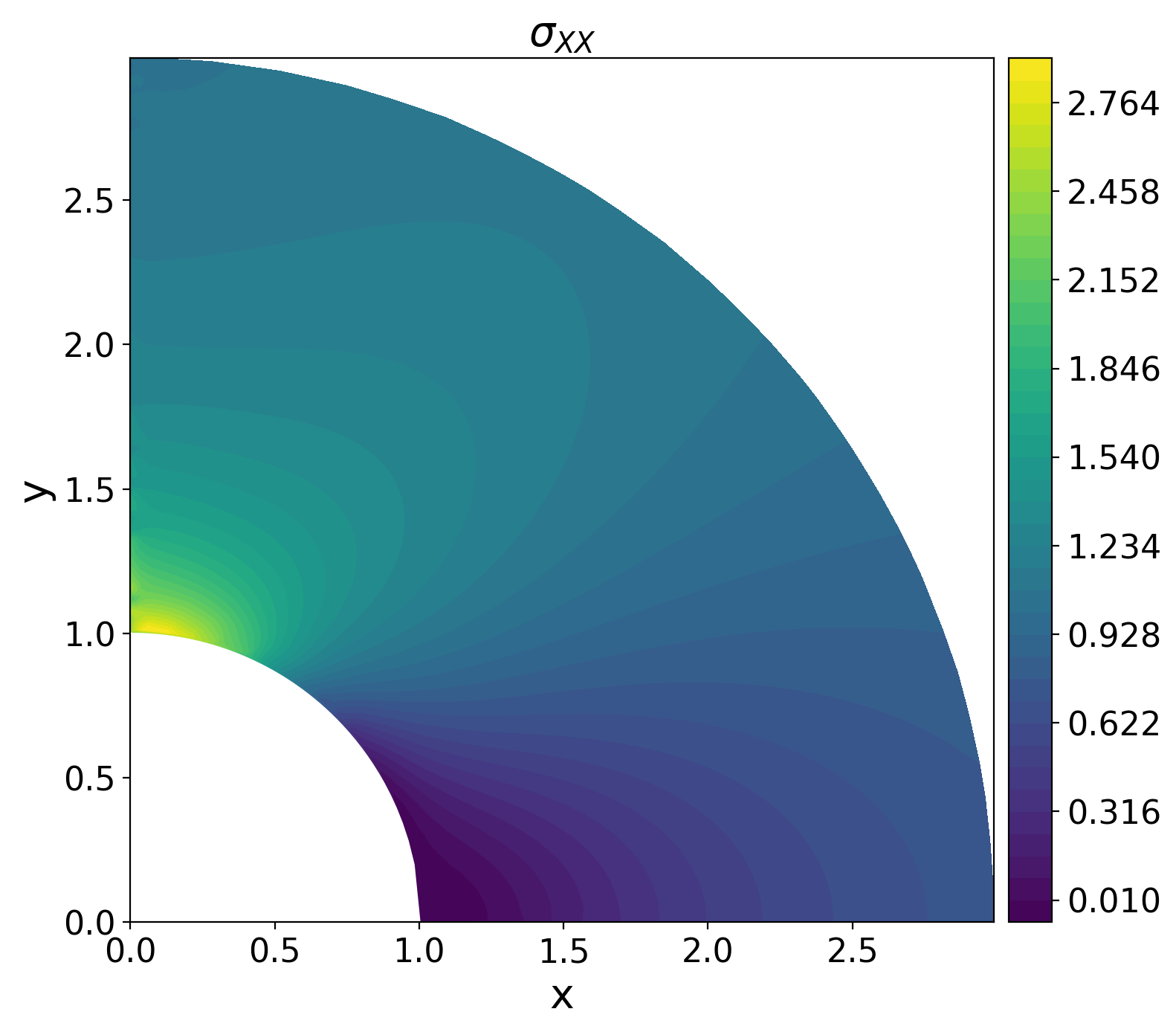}}
  \caption{Horizontal-stress contours for the plate with a hole problem for the quadratic formulations and L3 triangular discretization. Only the bond-associated results (c-d) produce a meaningful, stable solution, which is in good agreement with the exact solution (not shown here).}
  \label{fig:plate-contours}
\end{figure*}

The stability issue of the base models are illustrated in \cref{fig:plate-contours}, where the horizontal-stress contours, for the quadratic models on the most refined triangular mesh, are plotted. While the BA variants capture the exact solution in this problem, the instability issue ruins the solution for the base models.

%\begin{figure}[!ht]
  %\centering
  %\subfloat[][RK-PD]{\includegraphics[height=0.4\textwidth]{Hole_contour_dispX_RK-PD.png}}
  %\hspace*{0.5cm}
  %\subfloat[][GMLS-PD]{\includegraphics[height=0.4\textwidth]{Hole_contour_dispX_GMLS-PD.png}}

  %\subfloat[][BA-RK-PD]{\includegraphics[height=0.4\textwidth]{Hole_contour_dispX_BA-RK-PD.png}}
  %\hspace*{0.5cm}
  %\subfloat[][BA-GMLS-PD]{\includegraphics[height=0.4\textwidth]{Hole_contour_dispX_BA-GMLS-PD.png}}
  %\caption{Horizontal-displacement contours in the plate with hole problem using the quadratic formulations and the L3 triangular discretization. Only the bond-associated models are able to capture the correct behavior. The bond-associated solutions (c-d) are in good agreement with the exact solution (not shown here).}
  %\label{fig:plate-contours}
%\end{figure}

\begin{figure*}[!ht]
  \centering
  \subfloat[][Polar discretization - Quadratic]{\includegraphics[height=0.4\textwidth]{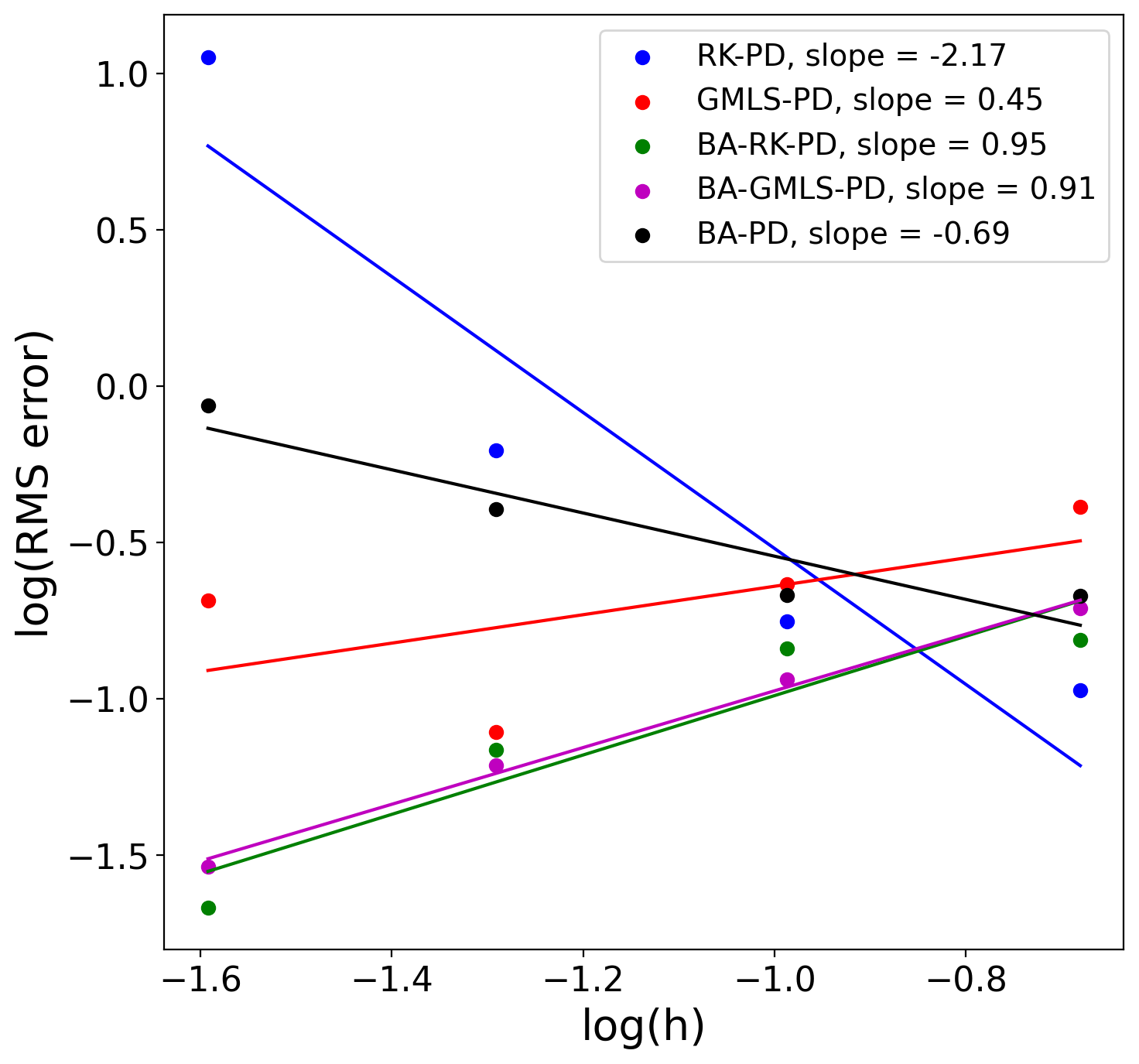}}
  \hspace*{1cm}
  \subfloat[][Triangular discretization - Quadratic]{\includegraphics[height=0.4\textwidth]{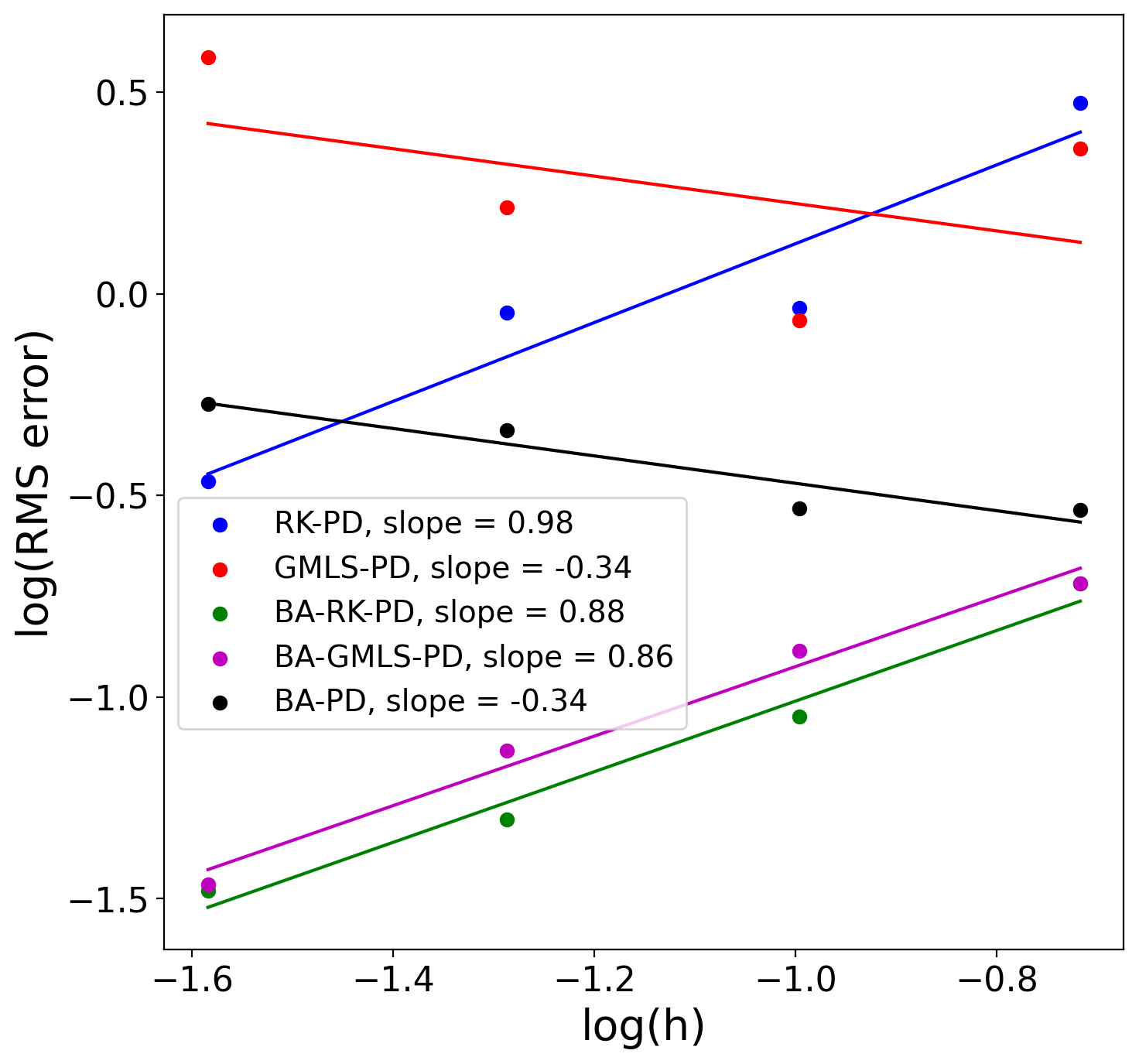}}
  \caption{Convergence for the plate with a hole problem in the mild near-incompressibility regime. Plots show stress RMS error values for the quadratic discretizations.}
  \label{fig:plate-incompressibility}
\end{figure*}

The methods considered here are also tested in the mild near-incompressibility regime ($\nu=0.495$) using the same example. \cref{fig:plate-incompressibility} shows the results for the quadratic formulations, where the higher-order, bond-associated variants are the only models that are able to converge to the analytical solution, obtaining a near-linear convergence rate. RK-PD, GMLS-PD, and BA-PD do not achieve convergence in this test. For higher levels of near-incompressibility, special treatment will be required to avoid volumetric locking (see, e.g.~\citep{moutsanidis2020treatment}).

\section{Conclusions}
\label{sec:conclusions}

The present work builds on Part I where a core formulation of PD in the framework of the correspondence modeling was developed based on the unification of RK-PD and GMLS approaches and with the addition of bond-associated stabilization. The core formulation exhibited a good combination of accuracy, stability, and robustness with respect to the choice of horizon size, attributes that are critical in engineering-scale applications. In the present effort, we studied the behavior of the proposed bond-associated models for the approximation of wave dispersion phenomena, and developed a new procedure for the imposition of natural boundary conditions for PD formulations in the strong form. Good performance in dynamic analysis and the ability to impose natural boundary conditions are likewise considered essential attributes of any framework that is to be applied for solving nonlinear, dynamic, and coupled multi-physics problems.

Wave propagation results presented here indicate that stability of numerical formulations that comes with bond-associated modeling is essential to accurately capture wave dispersion phenomena, especially at higher wave numbers. 

The main idea behind the enforcement of stress boundary conditions is to introduce the so-called natural-bc nodes on the boundaries where stress boundary conditions are enforced. Having introduced the natural-bc nodes, two neighborhood sets are defined. One such neighborhood set excludes the natural-bc nodes and is used for the computation of the non-local deformation gradient. This ensures that the formulation passes the patch test. The other neighborhood set includes the natural-bc nodes and is used for the computation of the internal force vector. This ensures that the stress information from the natural-bc nodes is appropriately accounted for in the in the internal force vector. The approach is similar to the modeling of damage in PD, where the broken bonds are not involved in the kinematic variable computation, and contribute only to the internal force evaluation. The proposed approach was successfully tested on a series of 2D linear elastostatics benchmark problems involving mixed boundary conditions, and also including mild near-incompressibility. It was shown that the bond-associated, higher-order (at least a second-order) model gives the best results. 

Future efforts will utilize the developed framework to solve dynamic problems involving material failure. We will apply the presented methodology to semi-Lagrangian PD \citep{behzadinasab2020semi}, which is well-suited to simulate very large deformation problems with fragmentation \citep{behzadinasab2020peridynamic}. This framework will be used to study fluid--structure interaction phenomena, where the surrounding fluid imposes stress boundary conditions on the solid structure that are naturally handled in the proposed framework.

\begin{acknowledgements}
  Y.~Bazilevs was partially supported through the Sandia contract No. 2111577. J.T.~Foster acknowledges funding under AFOSR MURI Center for Materials Failure Prediction through Peridynamics: project No. ONRBAA12-020 and Sandia National Laboratories contract No. 1885207. The authors also thank Nathaniel Trask for helpful discussions on the subject. 

\end{acknowledgements}

\bibliography{paper}

\end{document}